\documentclass[final]{siamart190516}
\usepackage[utf8]{inputenc}
\usepackage{xcolor}
\usepackage{mathtools}
\usepackage{graphicx}
\usepackage{amsfonts}
\usepackage{amssymb}
\usepackage{algorithmic}
\usepackage{subcaption}
\usepackage{tcolorbox}
\usepackage{scrextend}
\usepackage{listings}
\usepackage{bm}
\usepackage{float}
\usepackage{multirow}
\usepackage{afterpage}
\usepackage{etoolbox}
\usepackage{makecell}

\newcommand{\etal}{et al.\ }

\newcommand{\vv}[1]{ \boldsymbol{\mathbf{#1}} }

\newcommand*{\colorboxed}{}
\def\colorboxed#1#{%
  \colorboxedAux{#1}%
}
\newcommand*{\colorboxedAux}[3]{%
  \begingroup
    \colorlet{cb@saved}{.}%
    \color#1{#2}%
    \boxed{%
      \color{cb@saved}%
      #3%
    }%
  \endgroup
}

\title{A Deep Learning Approach for the Computation of Curvature in the Level-Set Method\thanks{Submitted to the SIAM Journal on Scientific Computing.\funding{This research acknowledges support from ONR N00014-17-1-2676.}}}

\author{Luis \'{A}ngel Larios-C\'{a}rdenas\thanks{Department of Computer Science, University of California, Santa Barbara, CA 93106, USA (\email{lal@cs.ucsb.edu}).}
\and Frederic Gibou\thanks{Department of Mechanical Engineering and Department of Computer Science, University of California, Santa Barbara, CA 93106, USA (\email{fgibou@ucsb.edu}).}}

\headers{Level-Set Curvature Estimation Using Deep Learning}{Luis \'{A}ngel Larios-C\'{a}rdenas and Frederic Gibou}

\begin{document}

\maketitle

\begin{abstract}
We propose a deep learning strategy to estimate the mean curvature of two-dimensional implicit interfaces in the level-set method.  Our approach is based on fitting feed-forward neural networks to synthetic data sets constructed from circular interfaces immersed in uniform grids of various resolutions.  These multilayer perceptrons process the level-set values from mesh points next to the free boundary and output the dimensionless curvature at their closest locations on the interface.  Accuracy analyses involving irregular interfaces, in both uniform and adaptive grids, show that our models are competitive with traditional numerical schemes in the $L^1$ and $L^2$ norms.  In particular, our neural networks approximate curvature with comparable precision in coarse resolutions, when the interface features steep curvature regions, and when the number of iterations to reinitialize the level-set function is small.  Although the conventional numerical approach is more robust than our framework, our results have unveiled the potential of machine learning for dealing with computational tasks where the level-set method is known to experience difficulties.  We also establish that an application-dependent map of local resolutions to neural models can be devised to estimate mean curvature more effectively than a universal neural network.
\end{abstract}




\section{Introduction}

Free boundary problems are a large class of models based on partial differential equations that arise from a variety of distant areas sharing the same mathematical structure \cite{CMGP16}.  Among the most common applications, we find models for simulating multiphase flows, heat conduction, propagation of fire fronts, image segmentation, solidification of multicomponent alloys, and morphogenesis.  There are two general numerical approaches for evolving interfaces subject to velocity fields: a Lagrangian or explicit representation and an Eulerian or implicit formulation.

In the Lagrangian strategy, we discretize the interface into a finite number of pieces and advect them by solving an elementary ordinary differential equation.  This method enjoys the advantage of simplicity for updating the positions of the elements and thus leads to accurate volume preservation \cite{Tryggvason;Bunner;Esmaeeli;etal:01:A-Front-Tracking-Met}.  This accuracy, however, deteriorates rather quickly when the velocity field causes pronounced topology deformations, and one does not resort to special procedures for maintaining interface smoothness and regularity \cite{Osher;Fedkiw:02:Level-Set-Methods-an}.  Also, explicit methods can become quite demanding when considering merging and splitting of moving fronts.  We collectively refer to the combination of these special techniques with the Lagrangian formulation to solve the evolution ODE as front-tracking methods.  Some representative developments on this branch may be found in \cite{Unverdi;Tryggvason:92:A-front-tracking-met, Juric;Tryggvason:96:A-Front-Tracking-Met, Juric;Tryggvason:98:Computations-of-Boil, Qian;Tryggvason;Law:98:A-Front-Tracking-Met}, and \cite{Tryggvason;Bunner;Esmaeeli;etal:01:A-Front-Tracking-Met} provides an excellent review on the subject.

\begin{figure}[t]
	\centering
	\includegraphics[width=0.65\textwidth]{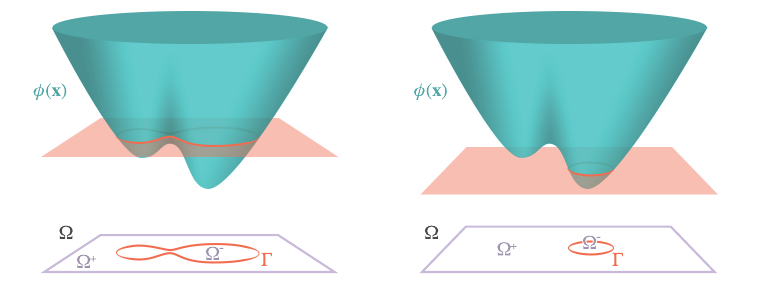}
	\caption{\small Time evolution of a level-set function, $\phi(\vv{x})$, its zero isocontour, $\Gamma$, and the computational domain, $\Omega$.}
	\label{fig.levelSet}
\end{figure}

Implicit methods circumvent the problems with boundary deformations and complicated topology maintenance schemes by using an implicit interface representation \cite{Osher;Fedkiw:02:Level-Set-Methods-an}.  Typical Eulerian formulations include the phase-field \cite{Braun;Murray:97:Adaptive-Phase-Field, ALS16, QB10, Elder;Grant;Provatas;etal:01:Sharp-interface-limi, Jeong;Goldenfeld;Dantzig:01:Phase-Field-Model-fo, Karma;Rappel:96:Phase-Field-Modeling, Karma;Rappel:97:Quantitative-Phase-F, Karma:01:Phase-Field-Formulat, Nestler;Danilov;Galenko:05:Crystal-growth-of-pu, Provatas;Goldenfeld;Dantzig:98:Efficient-Computatio, Provatas;Goldenfeld;Dantzig:99:Adaptive-Mesh-Refine}, the volume-of-fluid \cite{DeBar:74:Fundamentals-of-the-, Noh;Woodward:76:SLIC-simple-line-int, Hirt;Nichols:81:Volume-of-Fluid-VOF-, Youngs:84:An-Interface-Trackin, Benson:92:Computational-method, Sussman;Puckett:00:A-Coupled-Level-Set-, Renardy;Renardy:02:PROST:-A-Parabolic-R, Benson:02:Volume-of-Fluid-Inte, Aulisa;Manservisi;Scardovelli;etal:03:A-geometrical-area-p, Dyadechko;Shashkov:06:Moment-of-Fluid-Inte, Yang;James;Lowengrub;etal:06:An-adaptive-coupled-, Wang;Yang;Stern:12:A-new-volume-of-flui, CurvatureML19, VOFCurvature3DML19, DespresJourdren;MLDesignOfVOF;20}, and the level-set methods \cite{Osher;Sethian:88:Fronts-Propagating-w, MKZ98, Sussman;Fatemi;Smereka;etal:98:An-Improved-Level-Se, Sethian:99:Level-set-methods-an, Sussman;Puckett:00:A-Coupled-Level-Set-, Osher;Fedkiw:01:Level-Set-Methods:-A, Osher;Fedkiw:02:Level-Set-Methods-an, Enright;Fedkiw;Ferziger;etal:02:A-Hybrid-Particle-Le, Yang;James;Lowengrub;etal:06:An-adaptive-coupled-, Chen;Min;Gibou:07:A-Supra-Convergent-F, Min;Gibou:07:A-second-order-accur, Min;Gibou:07:Geometric-integratio, Min;Gibou:08:Robust-second-order-, Chen;Min;Gibou:09:A-numerical-scheme-f, Mirzadeh;Guittet;Burstedde;etal:16:Parallel-level-set-m, CLSVOF17, GFO18, Theillard;Gibou;Pollock:14:A-Sharp-Computationa, Rycroft;Gibou:12:Simulations-of-a-str}.  The main advantage of these frameworks is their natural ability to handle complex topological changes, as depicted in Figure \ref{fig.levelSet}.

Phase-field models are implicit methods that have been used extensively in the materials sciences to treat phase transitions or decompositions with complicated morphologies.  This method, based on the numerical solution of the nonlinear Cahn--Hilliard or Cahn--Allen equation \cite{ALS16}, introduces a diffuse profile where an auxiliary continuous order parameter identifies the phase \cite{Braun;Murray:97:Adaptive-Phase-Field}.  The phase-field approach has been considerably popular in the investigation of microstructural pattern formation in alloys and in the simulation of solidification phenomena, such as dendritic crystallization (see, for example, \cite{Karma;Rappel:97:Quantitative-Phase-F, Karma:01:Phase-Field-Formulat, Nestler;Danilov;Galenko:05:Crystal-growth-of-pu}).  Phase-field models, however, experience important difficulties when handling large domains.  Qin and Bhadeshia \cite{QB10} have stated that the possibility of adjusting the interface width to unrealistic physical values results in a loss of detail to represent sharp transitions.  In turn, these issues translate into severe time step restrictions and low accuracy.  Also, it is unclear how the parameters of the phase-field representation relate to physical quantities.

The volume-of-fluid (VOF) approach, pioneered by \cite{DeBar:74:Fundamentals-of-the-, Noh;Woodward:76:SLIC-simple-line-int, Hirt;Nichols:81:Volume-of-Fluid-VOF-}, is a popular method in computational fluid dynamics that uses volume fractions of fluid as index function values at each computational cell \cite{CurvatureML19, VOFCurvature3DML19}.  In two-phase flows, for example, a unit value corresponds to a cell full of fluid, while zero is associated with an empty cell.  Any cell with an intermediate value contains part of a free boundary or material interface, which is assembled through a piecewise linear or parabolic reconstruction \cite{Hirt;Nichols:81:Volume-of-Fluid-VOF-, Benson:02:Volume-of-Fluid-Inte, Renardy;Renardy:02:PROST:-A-Parabolic-R}.  VOF methods approximately conserve volume and mass locally, treat coalescing free boundaries automatically, require a minimum of stored information, and can be easily taken to three dimensions \cite{Hirt;Nichols:81:Volume-of-Fluid-VOF-}.  However, extending the algorithm to handling more materials is nontrivial \cite{Benson:92:Computational-method}.  Also, it is challenging to use only the volume fractions to compute the interface curvature, normal vectors, and forces because of the discontinuous nature of the VOF representation.  

An emerging line of research in VOF methods has been exploring the latest machine learning technologies to tackle the difficulties with interface reconstruction, advection, and estimation of geometric quantities.  Recently, Despr\'{e}s and Jourdren have addressed interface reconstruction and standard VOF procedures by proposing a VOF machine learning scheme adapted to bimaterial compressible Euler computations \cite{DespresJourdren;MLDesignOfVOF;20}.  These authors have trained a machine learning flux function with straight lines, arcs, and corners, aiming for an accurate transport/remap of reconstructed interfaces of various shapes.  Likewise, Qi \etal \cite{CurvatureML19} and Patel \etal \cite{VOFCurvature3DML19} have proposed data-driven strategies to improve on curvature estimation accuracy.  In \cite{CurvatureML19}, a feed-forward neural network processes the volume fractions of a nine-cell stencil and outputs the dimensionless curvature at the center cell.  This neural model was trained in uniform grids with a large synthetic data set, spanning a wide range of curvatures and orientations from circular interfaces.  Qi \etal have shown that their neural network produces satisfactory results in experiments with a static sinusoidal wave and with a moving free boundary in a flow solver.  Similarly, Patel \etal \cite{VOFCurvature3DML19} have devised supervised multilayer perceptrons to establish functional relationships between volume fractions in three-dimensional stencils and local interface curvatures.  They have presented a systematic approach for generating well-balanced data sets using spherical patches at different configurations.  Patel \etal have also tested their models on several shapes and coupled them with a multiphase solver to verify their effectiveness in standard simulations with bubbles.  When compared to conventional curvature computation methods, these models can easily outperform the convolution scheme and match the correctness of the height function method.

The level-set method is a powerful framework for evolving arbitrary interfaces, which are captured as the zero isocontour of high-dimensional functions (see Section \ref{sec.levelSetMethod}).  The interface is advected under a velocity field by solving a Hamilton--Jacobi equation using finite-difference, finite-volume, or finite-element schemes.  As a result, the dynamics of the moving boundary implies no need for asymptotic analysis, and only standard time step restrictions for stability and consistency are required \cite{Min;Gibou:07:A-second-order-accur}.  Because of its numerical nature, however, the most important level-set method difficulties are (1) maintaining the smoothness of the underlying level-set function and (2) conserving mass when the interface undergoes severe stretching or tearing \cite{MKZ98, Sussman;Puckett:00:A-Coupled-Level-Set-}.  Hybrid methods have been introduced to improve on mass conservation.  In the coupled level-set and VOF method by Sussman and Puckett \cite{Sussman;Puckett:00:A-Coupled-Level-Set-}, the VOF index function complements the level-set function to transport the implicit boundary.  There, the authors use a piecewise linear interface calculation (PLIC) to correct for and preserve the enclosed fluid mass or volume.  At the same time, a smooth level-set function is used to compute surface normal vectors and curvature \cite{CLSVOF17}.  More recently, Yang \etal \cite{Yang;James;Lowengrub;etal:06:An-adaptive-coupled-} have presented the adaptive coupled level-set and VOF volume tracking method for unstructured triangular grids.  Compared to previous developments, their adaptive-mesh algorithm is founded on an analytical PLIC for triangular grids.  Their method accurately resolves complex topological changes, areas of steep curvature, and near-contact regions of colliding fronts.  Another hybrid approach that has addressed mass loss is the particle-level-set method proposed by Enright \etal \cite{Enright;Fedkiw;Ferziger;etal:02:A-Hybrid-Particle-Le}.  In this work, marker particles are randomly seeded near the interface and are passively advected with a flow field.  These Lagrangian marker particles are subsequently used to rebuild the level-set function in underresolved regions.

The mass loss problem in the level-set method is especially worse when one uses coarse grids for discretization  \cite{Yang;James;Lowengrub;etal:06:An-adaptive-coupled-}.  Thus, multiple studies have focused on applying mesh refinement techniques \cite{Berger;Oliger:84:Adaptive-mesh-refine, Strain1999} to handle large numbers of discrete elements near evolving fronts.  On this line of research, Chen, Min, and Gibou \cite{Chen;Min;Gibou:07:A-Supra-Convergent-F, Chen;Min;Gibou:09:A-numerical-scheme-f} and Min and Gibou \cite{Min;Gibou;Ceniceros:06:A-supra-convergent-f, Min;Gibou:06:A-second-order-accur, Min;Gibou:07:Geometric-integratio, Min;Gibou:08:Robust-second-order-} have developed node-based, robust level-set tools and extrapolation and reinitialization schemes that have been incorporated into adaptive solvers (see \cite{GFO18} and the references therein).  The level-set methods on quadtrees (in two dimensions) and octrees (in three dimensions) proposed by Min and Gibou \cite{Min;Gibou:07:A-second-order-accur} are archetypical of applying local grid refinement to reduce the computational costs associated with high resolutions.  Mirzadeh \etal \cite{Mirzadeh;Guittet;Burstedde;etal:16:Parallel-level-set-m} have extended these level-set technologies on quad and octrees to parallel algorithms for distributed memory machines.  Mirzadeh and coauthors presented their algorithms while solving the Stefan problem with a combination of \verb|MPI| \cite{MPI14} and the \verb|p4est| \cite{Burstedde;Wilcox;Ghattas:11:p4est:-Scalable-Algo} library.  Local level-set methods have also been considered for limiting the band of grid points processed around the interface \cite{Peng;Merriman;Osher;etal:99:A-PDE-based-fast-loc} or for using efficient data structures that only record those mesh points \cite{Brun;Guittet;Gibou:12:A-local-level-set-me, Adalsteinsson;Sethian:95:A-Fast-Level-Set-Met, Nielsen;Museth:06:Dynamic-Tubular-Grid}.  The authors in \cite{Brun;Guittet;Gibou:12:A-local-level-set-me}, for example, have noted that only recording the set of adjacent vertices leads to the materialization of numerical noise when the level-set function is reinitialized.  They have also shown that level-set methods based on octrees have superior performance.

Curvature is one of the most important derived interface attributes for its relation to surface tension in physics and its regularization property in optimization.  It is thus crucial to compute mean curvature accurately from any free boundary representation.  The recent review by Popinet \cite{Popinet;NumModelsOfSurfTension;18} highlights the distinct classes of methods routinely used for curvature estimation in Eulerian surface tension models.  In the level-set method, specifically, the correctness of such an approximation depends on the smoothness of the level-set function.  One can enforce both level-set smoothness and regularity by reinitializing the underlying function using iterative procedures.  In this regard, and based on the observations in \cite{G.Russo;P.Smereka:00:A-level-set-method-f}, du Ch\'{e}n\'{e}, Min, and Gibou \cite{Chene;Min;Gibou:08:Second-order-accurat} have introduced a high-order accurate discretization of the reinitialization equation that produces second-order accurate curvature computations in the $L^\infty$ norm.  One should be cautioned, however, that high-order reinitialization procedures are costly and cannot always be realized on nonuniform grids.  Furthermore, although alternative, fast redistancing approaches exist (e.g., the fast marching method \cite{Sethian:96:A-Fast-Marching-Leve, Chacon;Vladimirsky:15:A-Parallel-Two-Scale, Chacon;Vladimirsky:12:Fast-two-scale-metho} and the fast sweeping method \cite{Zhao:04:A-Fast-Sweeping-Meth, Zhao:07:Parallel-Implementat, Detrixhe;Gibou;Min:13:A-parallel-fast-swee, Detrixhe;Gibou:14:Hybrid-Massively-Par}), these do not produce smooth enough level-set functions to guarantee satisfactory results.

Motivated by the work in \cite{CurvatureML19, VOFCurvature3DML19}, we introduce a deep learning strategy to estimate the mean curvature of two-dimensional interfaces in the level-set method.  Our approach is based on fitting feed-forward neural networks to synthetic data sets constructed from circular interfaces immersed in uniform grids of various resolutions.  These multilayer perceptrons process the level-set values from vertices next to the free boundary and output the dimensionless curvature at their closest locations on the interface.  Accuracy analyses with irregular interfaces show that our neural models are competitive, on average, with traditional numerical schemes.  In particular, our neural networks approximate mean curvature with comparable precision when we deal with relatively coarse resolutions, when the interface features steep curvature regions, and when the number of iterations for level-set reinitialization is small.  Unlike the machine learning advancements in VOF \cite{CurvatureML19, VOFCurvature3DML19}, our case studies reveal that there is still room for improvement given that the level-set numerical schemes are very robust.  Specifically, we have noticed that a numerical two-step procedure has superior performance in the absence of sharp turns, as the resolution increases, and as one uses more operations for level-set redistancing.

We organize the contents of our report as follows: we present the essentials of the level-set method in section \ref{sec.levelSetMethod} before describing the deep learning curvature approximation in section \ref{sec:DNN}.  We devote section \ref{sec:results} to numerical experiments and analyses, and section \ref{sec:conclusions} draws some conclusions and provides avenues for future work.


\section{The level-set method}
\label{sec.levelSetMethod}

The level-set method \cite{Osher;Sethian:88:Fronts-Propagating-w} is an Eulerian formulation of the numerical evolution of an implicit interface, captured as the zero isocontour of the level-set function $\phi(\vv{x}): \mathbb{R}^n \rightarrow \mathbb{R}$.  Given the computational domain $\Omega \subseteq \mathbb{R}^n$ and an interface $\Gamma$, $\phi$ is defined as the signed distance function
\begin{equation}
\phi(\vv{x}) = 
	\begin{cases}
	-d, & \vv{x} \in \Omega^-, \\
	+d, & \vv{x} \in \Omega^+, \\
	 0, & \vv{x} \in \Gamma,
	\end{cases}
\label{eq.levelSetFunction}
\end{equation}
where $d$ is the Euclidean distance from $\vv{x}$ to $\Gamma$, and $\Omega^-$ and $\Omega^+$ are resulting inside and outside regions in the partitioned domain.

Assuming that $\vv{v}(\vv{x})$ is a velocity field defined for all $\vv{x} \in \Omega$, the evolution of the level-set function $\phi(\vv{x})$ satisfies
\begin{equation}
\phi_t + \vv{v} \cdot \nabla\phi = 0,
\label{eq.levelSetEquation}
\end{equation}
which is known as the level-set equation.

\begin{figure}[h]
	\centering
	\includegraphics[width=0.4\textwidth]{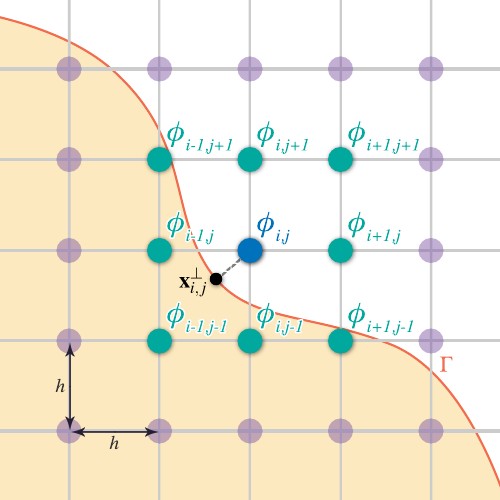}
	\caption{\small Level-set values used in the standard computation of the curvature at a grid node $(i,j)$.  The point $\vv{x}_{i,j}^\perp$ denotes the normal projection of $\vv{x}_{i,j}$ onto the interface, $\Gamma$.}
	\label{fig.stencil}
\end{figure}

Figure \ref{fig.stencil} depicts the zero level set, $\Gamma$, of a bivariate function, $\phi(x,y): \mathbb{R}^2 \rightarrow \mathbb{R}$, traveling across a portion of a discretized domain, $\Omega$.  It also shows the nine $\phi$ values employed in the numerical estimation of mean curvature, $\kappa$, at node $(i,j)$.  The formula to approximate $\kappa$ is given by
\begin{equation}
\kappa = \nabla \cdot \frac{\nabla \phi}{|\nabla \phi|} = \frac{\phi_x^2 \phi_{yy} - 2\phi_x\phi_y\phi_{xy} + \phi_y^2 \phi_{xx}}{\left( \phi_x^2 + \phi_y^2 \right)^{3/2}},
\label{eq.numericalCurvature}
\end{equation}
where the partial derivatives of $\phi$ at vertex $(i,j)$ can be approximated using second-order accurate central differences involving $\phi_{p,q}$, for $i-1 \leq p \leq i+1$ and $j-1 \leq q \leq j+1$, in the case of a uniform Cartesian grid.

In practice, one opts for a signed distance function, $\phi(x,y)$, because it produces robust numerical results \cite{Sussman;Smereka;Osher:94:A-Level-Set-Approach}.  A signed distance function is also uniquely determined as the viscosity solution of the Eikonal equation, and it satisfies the uniform gradient condition $|\nabla\phi| = 1$ \cite{Min:10:On-reinitializing-le}.  However, as $\Gamma$ evolves, $\phi$ deteriorates numerically and drifts away from its desired signed distance form, developing noisy features that get amplified when used in the approximation of partial derivatives in \eqref{eq.numericalCurvature}.  Consequently, it is standard to reinitialize $\phi$ into a signed distance function periodically. The equation for level-set reinitialization was first introduced by \cite{Sussman;Smereka;Osher:94:A-Level-Set-Approach}:
\begin{equation}
\phi_\tau + S(\phi^0)(|\nabla\phi| - 1) = 0,
\label{eq.reinitialization}
\end{equation}
where $\tau$ represents pseudo time, $\phi^0$ is the level-set function before reshaping, and $S(\phi^0)$ is a smoothed-out sign function.

To evolve the pseudotransient partial differential equation \eqref{eq.reinitialization} to steady-state ($\phi_\tau = 0$ and $|\nabla\phi| = 1$), one uses a TVD Runge--Kutta scheme in time and a Godunov spatial discretization of the Hamiltonian $H(\nabla\phi) = S(\phi^0)(|\nabla\phi| - 1)$ in a combination of Euler steps.  It is typical to use between 5 and 20 iterations to reinitialize a level-set function, with more iterations translating into smoother $\phi$ values and thus a more accurate approximation of $\kappa$.  This number of iterations is not a prescribed quantity, and it usually depends on how close $\phi^0$ is to a signed distance function.  In addition, one should keep in mind that the computational cost associated with redistancing is significant when the number of iterations is large.


\section{Deep learning and neural network curvature approximation}
\label{sec:DNN}

\subsection{Fundamentals}

An artificial neural network is a computational graph of elementary units that lie interconnected in some particular way to compute a function of its inputs by using connection weights as intermediate parameters \cite{A18}.  Neural networks are powerful supervised learning methods \cite{Mehta19}.  They learn a target function by successively adapting their weights to minimize a loss metric given their current parameter configuration and a set, $\mathcal{D}$, of training data pairs.

\begin{figure}[t]
	\centering
	\begin{subfigure}[b]{0.4\textwidth}
		\includegraphics[width=\textwidth]{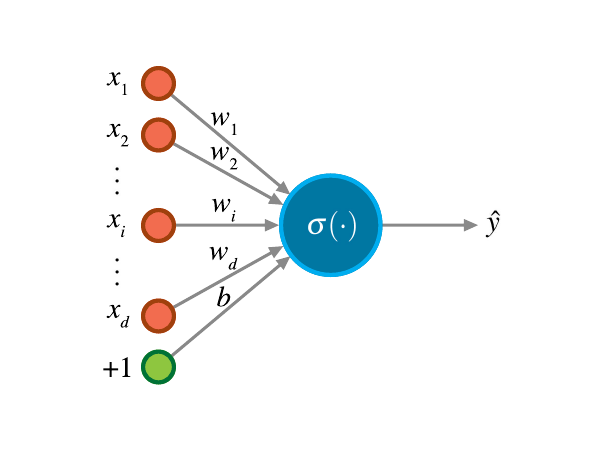}
		\caption{Perceptron}
		\label{fig.perceptron}
	\end{subfigure}\hfill%
	\begin{subfigure}[b]{0.4\textwidth}
		\includegraphics[width=\textwidth]{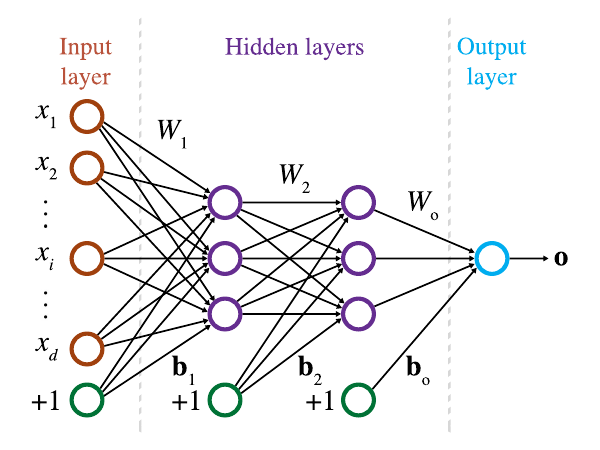}
		\caption{Multilayer perceptron}
		\label{fig.mlp}
	\end{subfigure}
	\caption{\small The artificial neural network basic architectures.}
	\label{fig.nnet}
\end{figure}

The simplest neural network architecture is known as a \textit{perceptron} and consists of a single input layer and one nonlinear output neuron (see Figure \ref{fig.perceptron}).  Given an input layer with $d$ linear units, $\vv{x} \in \mathbb{R}^d$, a vector of edge parameters, $\vv{w} \in \mathbb{R}^d$, and an offset bias, $b \in \mathbb{R}$, the neural prediction is computed as

\begin{equation}
\hat{y} = \sigma( \vv{w} \cdot \vv{x} + b ) = \sigma\left( \sum_{i=1}^d w_i x_i + b \right),
\label{eq.perceptron}
\end{equation}
where $\sigma(\cdot)$ plays the role of an activation function.  Traditional $\sigma(\cdot)$ options include the step function, sigmoids, and hyperbolic tangent, although more recent nonlinearities like the rectified linear units (ReLUs), leaky rectified linear units, and exponential linear units have become commonplace \cite{Mehta19}.

The training algorithm of a perceptron works by feeding each data pair $(\vv{x}_p, y_p) \in \mathcal{D}$ into the network and outputting a prediction $\hat{y}_p$.  The weights are then updated based on the error $e(\vv{x}_p) = y_p - \hat{y}_p$ as follows:

\begin{equation}
\vv{w} \leftarrow \vv{w} + \alpha e(\vv{x}_p) \vv{x}_p \qquad \text{or } \qquad \vv{w} \leftarrow \vv{w} + \alpha \sum_{(\vv{x}_p, y_p) \in \mathcal{B}} e(\vv{x}_p) \vv{x}_p,
\label{eq.perceptronUpdate}
\end{equation}
where $\alpha$ is the learning rate and $\mathcal{B} \subset \mathcal{D}$ is a batch of random samples.

The perceptron learning process is iteratively executed until convergence.  Every sample $(\vv{x}_p, y_p) \in \mathcal{D}$ may be cycled through many times, and each such iteration is known as an epoch \cite{A18}.  The perceptron learning algorithm based on equation \eqref{eq.perceptronUpdate} is an (mini-batch) stochastic gradient-descent method.  Its goal is to minimize the prediction squared error by performing gradient-descent updates with respect to randomly selected (batches of) training samples \cite{A18}.

\textit{Multilayer neural networks} consist of several (perceptron) units arranged in multiple layers.  The input layer transmits data into the network, the hidden layers perform some additional (non)linear calculations, and the output layer yields the predictions.  The output from each intermediate layer is carried over as the (weighted) input to the next level in the architecture.  The output layer is usually dedicated to simple classification (i.e., discrete result) or linear regression (i.e., continuous result) \cite{Mehta19}.  Our discussion centers around neural networks that predict a continuous output.  Figure \ref{fig.mlp} depicts a multilayer, feed-forward, and fully connected neural network with 2 hidden layers and a one-unit output layer.

The neural network parameters, $\vv{\Theta}$, are contained in bias vectors, $\vv{b}_m$, and weight matrices, $W_m$, for $1 \leq m \leq M + 1$, where $M$ is the number of hidden layers.  The weight matrix for the connections between the input with $d$ units and the first hidden layer with $m_1$ neurons is denoted as $W_1 \in \mathbb{R}^{d \times m_1}$; the weights between the $r$th and the $(r+1)$th hidden layer are provided in $W_{r+1} \in \mathbb{R}^{m_r \times m_{r+1}}$; and the matrix that holds the weights between the last hidden layer and the output layer with $o$ units is $W_o \in \mathbb{R}^{m_M \times o}$ \cite{A18}.  The layerwise recurrence equations that transform an input $\vv{x}$ into the predicted output $\vv{o}$ are expressed as

\begin{equation*}
\begin{cases}
	\vv{h}_1 = \sigma(W_1^T \vv{x} + \vv{b}_1) & \text{from input to first hidden layer}, \\
	\vv{h}_{r+1} = \sigma(W_{r+1}^T \vv{h}_r + \vv{b}_{r+1}), 1 \leq r \leq M - 1 & \text{from } r\text{th to } (r+1)\text{th hidden layers}, \\
	\vv{o} = \sigma(W_o^T \vv{h}_M + \vv{b}_o) & \text{from hidden to output layer}, 
\end{cases}
\end{equation*}
where $\sigma(\cdot)$ is an activation function.

Training a neural network involves computing a loss function and using gradient descent to find the best set of weights and biases \cite{Mehta19}.  For continuous output data, common loss functions are the mean squared error (MSE)

\begin{equation}
e(\vv{\Theta}) = \frac{1}{N}\sum_{p=1}^{N}\left( y_p - \hat{y}_p(\vv{\Theta}) \right)^2
\label{eq.mse}
\end{equation}
and the mean absolute error (MAE)

\begin{equation}
e(\vv{\Theta}) = \frac{1}{N}\sum_{p=1}^{N}\left| y_p - \hat{y}_p(\vv{\Theta}) \right|,
\label{eq.mae}
\end{equation}
where $N = |\mathcal{D}|$.  It is also not uncommon to add regularization terms in \eqref{eq.mse} and \eqref{eq.mae} to reduce the effect of intrinsic noise and encourage generalization \cite{Mehta19}.

Some optimizers often employed for training include stochastic gradient descent, Nesterov, RMSProp, Adadelta, Adagrad, and Adam.  We refer the interested reader to \cite{Mehta19, A18, Keras15} for a complete exposition of gradient descent and cost function optimizers.

The gradients for feed-forward networks are calculated through backpropagation \cite{Mehta19}.  This algorithm computes the gradient of a composition of functions by a direct application of dynamic programming.  The process comprises two phases: a forward phase, in which the network outputs and the local derivatives are calculated, and a backward phase, where the products of these local values are accumulated and the gradient of the loss function with respect to the different weights is used to update $\vv{\Theta}$ \cite{A18}.  The reader may consult \cite{A18} and \cite{Mehta19} for more details about backpropagation.

\subsection{Training a neural network to approximate curvature}
\label{sec.buildingNNets}

Recent advancements in deep learning have shown that multilayer perceptrons with several hidden layers possess the mathematical capacity to approximate any complex function \cite{Mehta19}.  Armed with this principle, we seek to build a neural network to accurately estimate the mean curvature of two-dimensional interfaces in the level-set method.  Following the approaches in \cite{CurvatureML19, VOFCurvature3DML19}, our goal is to learn a relation,
\begin{equation}
h\kappa_{i,j} = f\begin{pmatrix}
	\phi_{i-1,j+1}, & \phi_{i,j+1}, & \phi_{i+1,j+1}, \\
	\phi_{i-1,j}, & \phi_{i,j}, & \phi_{i+1,j}, \\
	\phi_{i-1,j-1}, & \phi_{i,j-1}, & \phi_{i+1,j-1}
\end{pmatrix},
\label{eq.neuralHKappa}
\end{equation}
where $h\kappa_{i,j}$ is the dimensionless curvature at a location on the zero isocontour of the level-set function $\phi(\vv{x}): \mathbb{R}^2 \rightarrow \mathbb{R}$ that is the closest to node $(i,j)$ (see Figure \ref{fig.stencil}).

Our level-set curvature neural network is a supervised model that employs training samples $(\vv{\phi}_p, o_p) \in \mathcal{D}$, $\vv{\phi}_p \in \mathbb{R}^9$, $o_p \in \mathbb{R}$, to build the association between $o_p = h\kappa$ and the nodal $\phi$ values.  As in \cite{CurvatureML19}, we train and test our model with synthetic data sets constructed from circular interfaces of varying radii.  These interfaces correspond to the zero level sets of the functions:

\begin{equation}
\phi_{cs}(x,y) = \sqrt{(x - x_0)^2 + (y - y_0)^2} - r
\label{eq.implicitCircleDist}
\end{equation}
and
\begin{equation}
\phi_{cn}(x,y) = (x - x_0)^2 + (y - y_0)^2 - r^2,
\label{eq.implicitCircleNonDist}
\end{equation}
where $(x_0, y_0)$ are the coordinates of the center of the interface, $\Gamma$, and $r$ is the circle radius.  Unlike \eqref{eq.implicitCircleDist}, \eqref{eq.implicitCircleNonDist} is not a signed distance function and therefore needs to be reinitialized via the numerical techniques introduced in section \ref{sec.levelSetMethod}.  Following \cite{CurvatureML19}, we opt for circular interfaces because once we know their radii, the target dimensionless curvatures can be calculated as $h/r$.

The neural network spatial domain is defined as the nonperiodic unit square $\Omega \equiv [0,1] \times [0,1]$, which is discretized with $\rho$ equally spaced grid points along the Cartesian directions.  Given a mesh size, $h = \Delta x = \Delta y = 1 / ( \rho - 1 )$, we obtain the number of circular interfaces with

\begin{equation}
\nu(\rho) = \left\lceil \frac{\rho - 8.2}{2} \right\rceil + 1,
\label{eq.nCircles}
\end{equation}
which guarantees that any resolution gets sufficient coverage with training interfaces.

To make sure that the circles are fully contained in $\Omega$, we randomly set their center coordinates $(x_0, y_0) \in [0.5 - h/2, 0.5 + h/2]^2$ and restrict their minimum and maximum radii to $1.6h$ and $1/2 - 2h$.  Furthermore, we generate data pairs for each radius, $r_\iota$, $1 \leq \iota \leq \nu(\rho)$, five times to introduce randomness into the learning stage.  These policies, together with \eqref{eq.nCircles}, guarantee that (1) at least four grid points lie inside the smallest circle, (2) any vertex next to $\Gamma$ has a well-defined nine-point stencil (see Figure \ref{fig.stencil}), and (3) model generalization is encouraged by repeated sampling from circles with the same radius \cite{Mehta19, A18}.

The learning set $\mathcal{D}$ is then a collection of $N$ samples of the form

\small
\begin{multline}
\left( \vv{\phi}_p, o_p \right) =\\
\big(\left[\phi_{i-1,j+1}, \phi_{i,j+1}, \phi_{i+1,j+1}, \phi_{i-1,j}, \phi_{i,j}, \phi_{i+1,j}, \phi_{i-1,j-1}, \phi_{i,j-1}, \phi_{i+1,j-1}\right], h\kappa_{i,j}\big),
\label{eq.sample}
\end{multline}\normalsize
where a sample is generated for any vertex that either sits on $\Gamma$ or has an outgoing vertical or horizontal edge crossed by the free boundary.  In particular, the input values, $\vv{\phi}_p$, in \eqref{eq.sample} are obtained from the corresponding circular interfaces in \eqref{eq.implicitCircleDist} and \eqref{eq.implicitCircleNonDist}.  In the case of $\phi_{cn}(x,y)$, we employ 5, 10, 15, and 20 iterations to reinitialize the level-set function.  According to preliminary experiments, redistancing not only helps to produce robust numerical results (as alluded to in section \ref{sec.levelSetMethod}) but also helps to remove unstructured noise that can compromise the models' complexity and learning outcome.  Furthermore, by using several degrees of reinitialization our neural networks are capable of building functional connections between different input versions of the same stencil.  In like manner, this allows our models to accommodate a decent range of reinitialization operations that practitioners often use to regularize their level-set functions.  Also, for all pairs $(\vv{\phi}_p, o_p)$, we collect their negated version $(-\vv{\phi}_p, -o_p)$, so that the neural network can account for concave regions.

We have used TensorFlow \cite{Tensorflow15} and Keras \cite{Keras15} to assemble a dictionary of three neural networks.  Each model was designed for a space resolution, $\rho$, of 256, 266, and 276 grid points per unit length.  Then, they were trained to solve for $h\kappa$ given their own data sets, $\mathcal{D}$.  We have followed the literature convention \cite{Mehta19, A18} and split these data sets into three parts: training (70\%), validation (15\%), and test (15\%).  Table \ref{tbl.sampleStats} gives the number of samples collected with this procedure.  We point out that the data pairs in the learning stage and in the upcoming analyses were generated in C++ with our in-house implementation of the parallel adaptive level-set method of \cite{Mirzadeh;Guittet;Burstedde;etal:16:Parallel-level-set-m}.

\begin{table}[!b]
	\centering
	\footnotesize
	\bgroup
	\def\arraystretch{1.1}%
	\begin{tabular}{|c|c|c|c|}
		\hline
		$\rho$ & Training & Test & Validation \\
		\hline \hline
		256 & 3'145,410 & 674,017 & 674,017 \\ \hline
		266 & 3'399,948 & 728,560 & 728,560 \\ \hline
 		276 & 3'664,188 & 785,184 & 785,184 \\
		\hline
	\end{tabular}
	\egroup
	\caption{\small Number of samples collected for each neural network.}
	\label{tbl.sampleStats}
\end{table}

Our proposed level-set curvature neural networks have input layers with nine linear units and a single linear neuron in their output layers.  To stabilize the learning process, we scaled $\vv{\phi}_p$ in \eqref{eq.sample} and used z-scores during training and testing.  We also included only ReLU neurons in all the intermediate layers because of their beneficial nonsaturating property that ameliorates the problem of vanishing gradients \cite{Mehta19}.  For training, we considered batches of size 32, selected the MSE loss function, and employed the Adam optimizer with a learning rate of 0.00015.  Additionally, we exercised the early-stopping technique, which monitored the models' generalization performance through the MAE of the validation sets.  For this, we set a patience parameter of 30 epochs.  In the end, the learning stage for each of the networks never surpassed 200 epochs.  Table \ref{tbl.nnetStats} shows the optimal configurations for our three-element neural dictionary alongside some statistics from their learning stage.

\begin{table}[!b]
	\centering
	\footnotesize
	\bgroup
	\def\arraystretch{1.1}%
	\begin{tabular}{|c|c|c|c|c|c|c|}
		\hline
		$\rho$ & \makecell{Training\\epochs} & \makecell{Hidden\\layers} & \makecell{Units per\\hidden layer} & Test MSE & Test MAE & Max AE \\
		\hline \hline
		256 & 103 & 4 & [128, 128, 128, 128] & $3.86\times 10^{-7}$ & $2.91\times 10^{-4}$ & 0.154 \\ \hline
		266 & 72 & 4 & [140, 140, 140, 140] & $3.64\times 10^{-7}$ & $2.70\times 10^{-4}$ & 0.122 \\ \hline
 		276 & 33& 4 & [140, 140, 140, 140] & $5.42\times 10^{-7}$ & $3.01\times 10^{-4}$ & 0.164 \\
		\hline
	\end{tabular}
	\egroup
	\caption{\small Best neural network configurations, their test set statistics, and their maximum absolute error (Max AE) over $\mathcal{D}$.}
	\label{tbl.nnetStats}
\end{table}

\begin{figure}[!t]
	\centering
	\begin{subfigure}[b]{0.35\textwidth}
		\includegraphics[width=\textwidth]{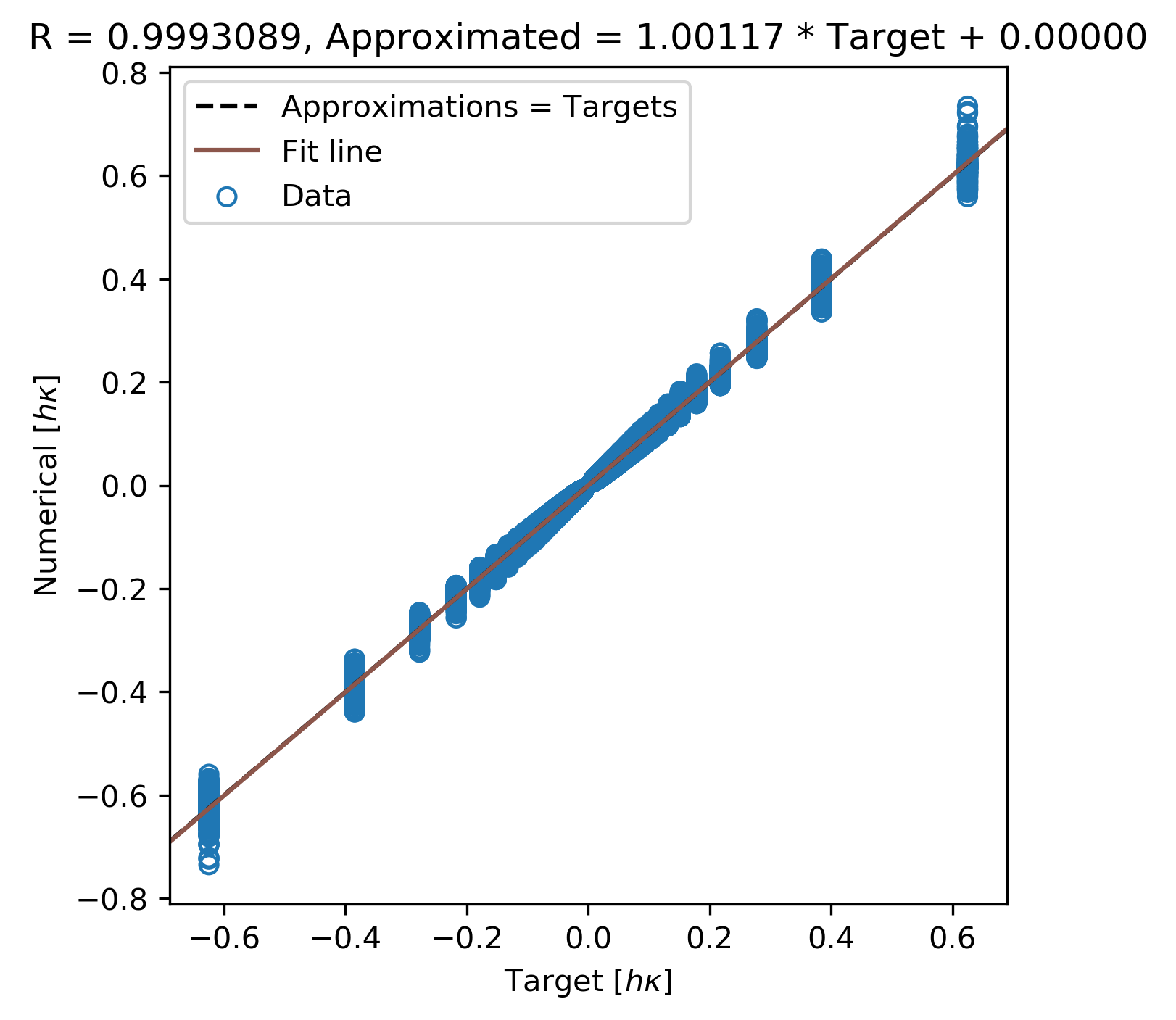}
        \caption{Numerical method}
        \label{fig.learningResults.266.numerics}
    \end{subfigure}
    ~
	\begin{subfigure}[b]{0.35\textwidth}
		\includegraphics[width=\textwidth]{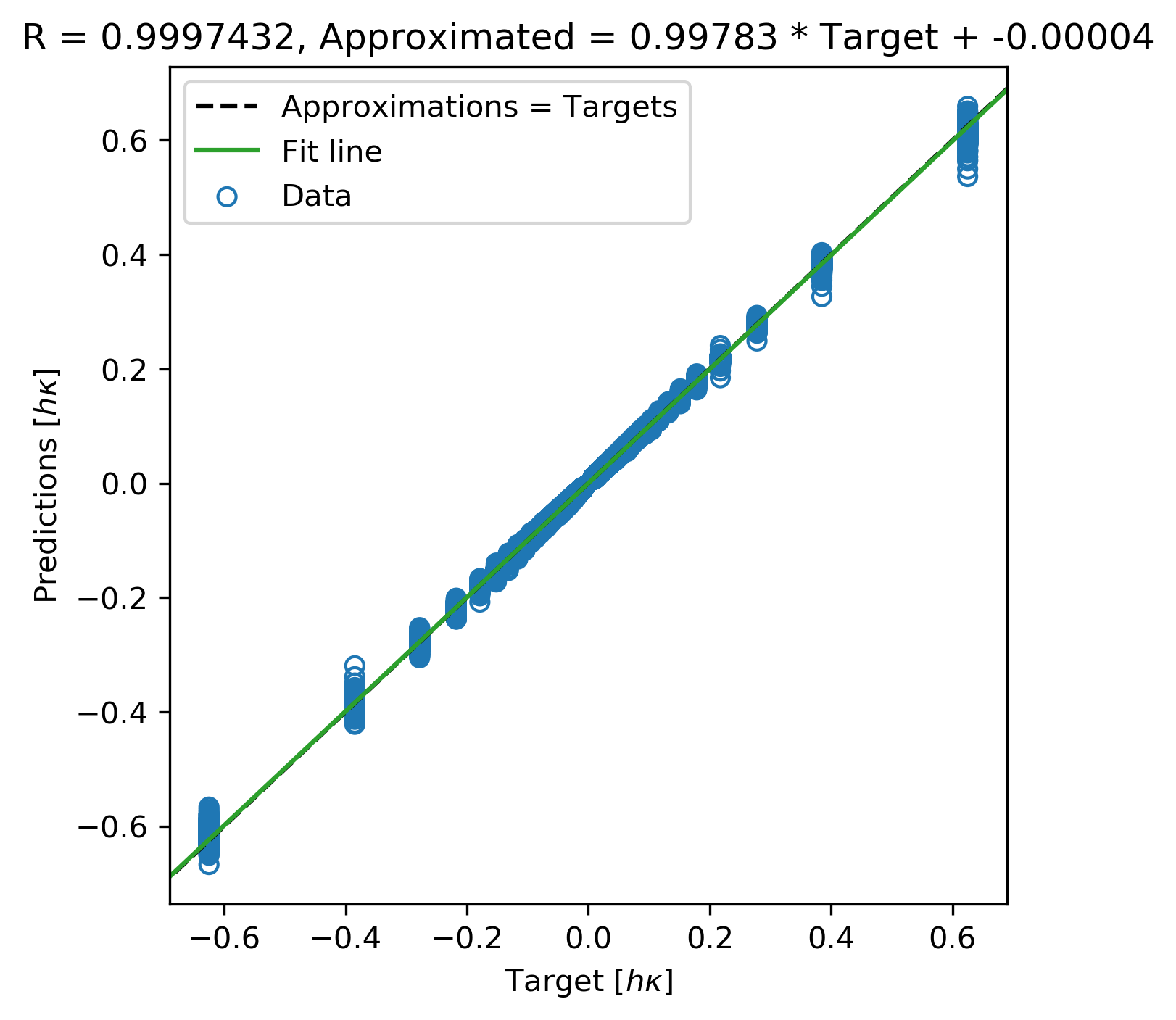}
        \caption{Training set}
        \label{fig.learningResults.266.train}
    \end{subfigure}
    
	\begin{subfigure}[b]{0.35\textwidth}
		\includegraphics[width=\textwidth]{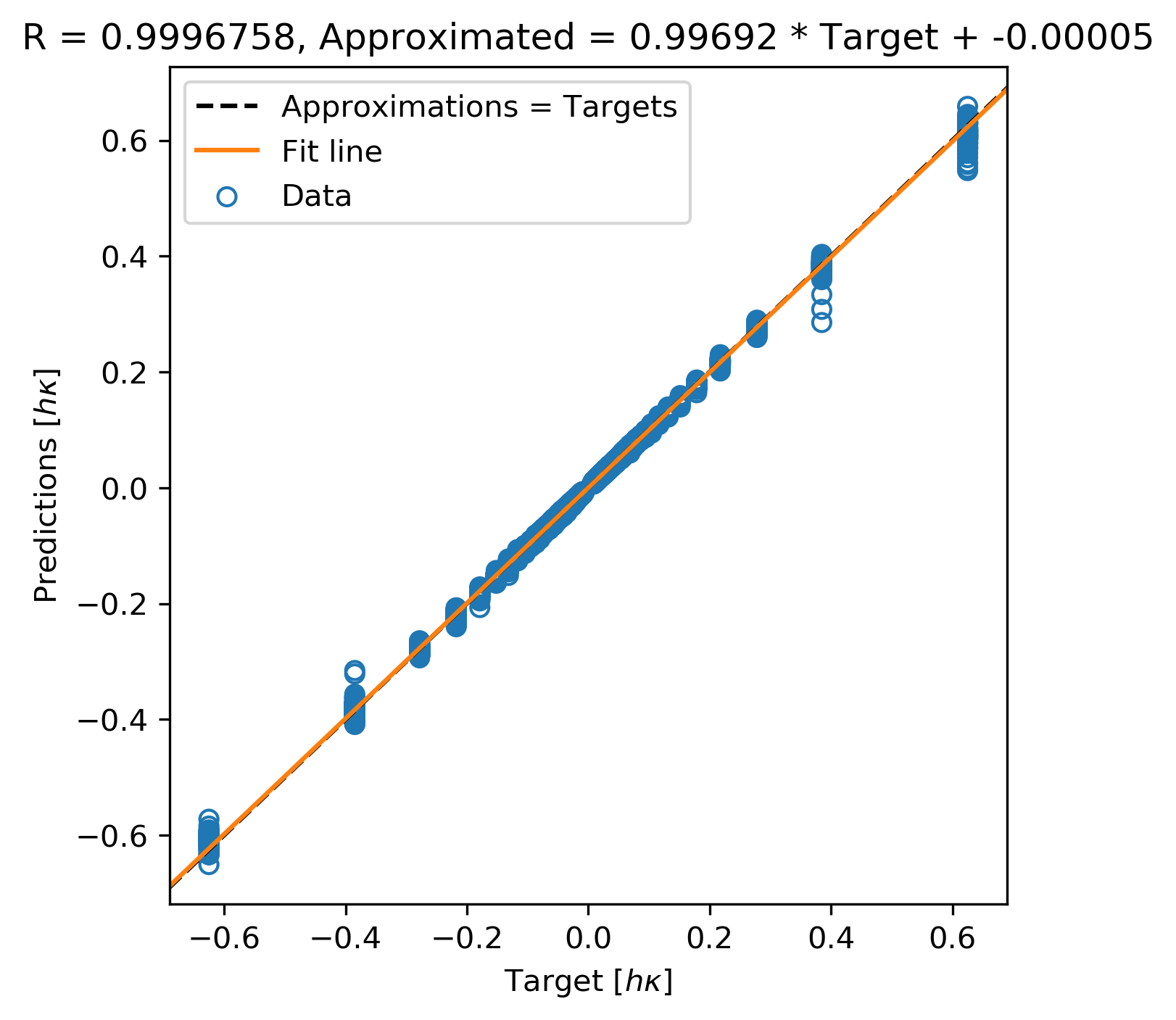}
        \caption{Testing set}
        \label{fig.learningResults.266.test}
    \end{subfigure}
    ~
	\begin{subfigure}[b]{0.35\textwidth}
		\includegraphics[width=\textwidth]{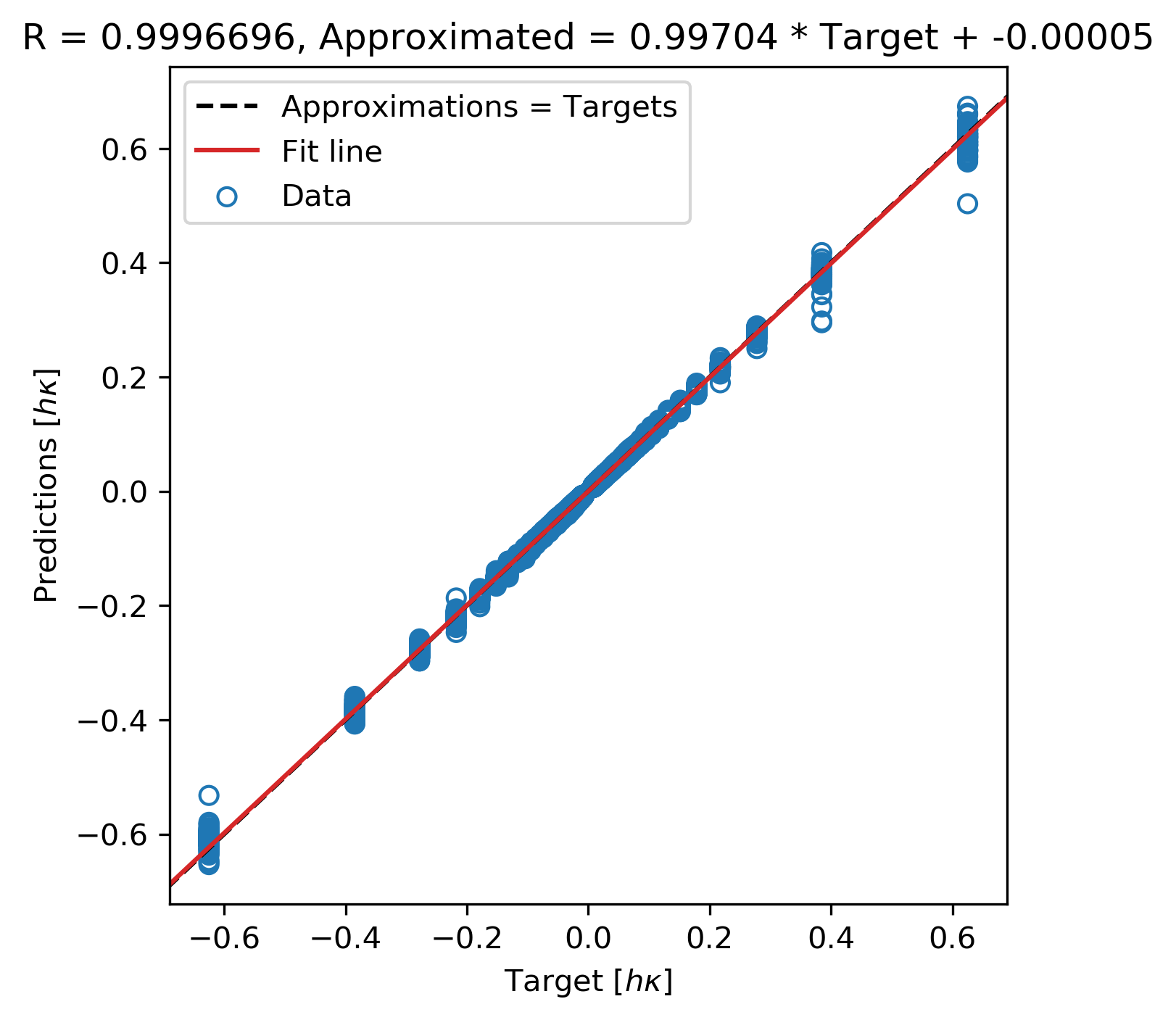}
        \caption{Validation set}
        \label{fig.learningResults.266.val}
    \end{subfigure}
	\caption{\small Correlation between expected and approximated $h\kappa$ using (a) the numerical method on the entire set $\mathcal{D}$ and (b)--(d) the neural network for the $266 \times 266$ grid resolution.}
	\label{fig.learningResults.266}
\end{figure}

Figure \ref{fig.learningResults.266} illustrates the quality of the fit achieved at the end of training by the neural network adapted for $\rho = 266$.  Each plot corresponds to the learning subsets extracted from $\mathcal{D}$.  To save space, we omit the results for the other two models because their outcomes were similar.  As a prelude to section \ref{sec:results}, Figure \ref{fig.learningResults.266} reveals that our neural model is competitive, if not marginally better in this case, than the numerical method.  This can be verified in the correlation factors and contrasted in the curvature estimations for circles with small radii (as $|h\kappa| \rightarrow 0.625)$.

To conclude this section, we remark that our level-set curvature neural networks and the methodology described above have been the products of diligent experimentation.  Regarding architectural design, for example, we tinkered with noise layers and layerwise pretraining via stacked denoising autoencoders \cite{A18} with no major success.  In addition, we tried to build a ``universal'' neural model to approximate curvature independently of the input domain resolution.  For this, we scaled the data pairs as $\frac{1}{h}(\vv{\phi}_p, o_p) \in \mathcal{D}$ and discarded normalization.  The resulting universal model, however, did not yield satisfactory $h\kappa$ estimations when the cell's width differed considerably from the training resolution.  This observation validated a well-known limitation of neural networks when they are used for extrapolation.  That is, because of their interpolative nature, functions can only be approximated within a span of the sample data used for training \cite{BNK;2020;ML-for-FD, Gibou:2019aa}.  For this reason, we hypothesize that a map of domain resolutions to neural networks has greater potential to infer the dimensionless curvature of implicit interfaces more effectively.  Our current approach has proven to deliver similar precision to conventional numerical schemes by breaking the inference problem into more tractable components.  Hence, we establish that our strategy remains valid (1) when the training resolutions are chosen appropriately and (2) when one extrapolates $h\kappa$ from a local resolution that is close to the training mesh size.  We have accounted for these constraints in the numerical experiments of section \ref{sec:results}.


\section{Experiments and results}
\label{sec:results}

\begin{figure}[!b]
	\centering
	\begin{subfigure}[b]{0.4\textwidth}
		\includegraphics[width=\textwidth]{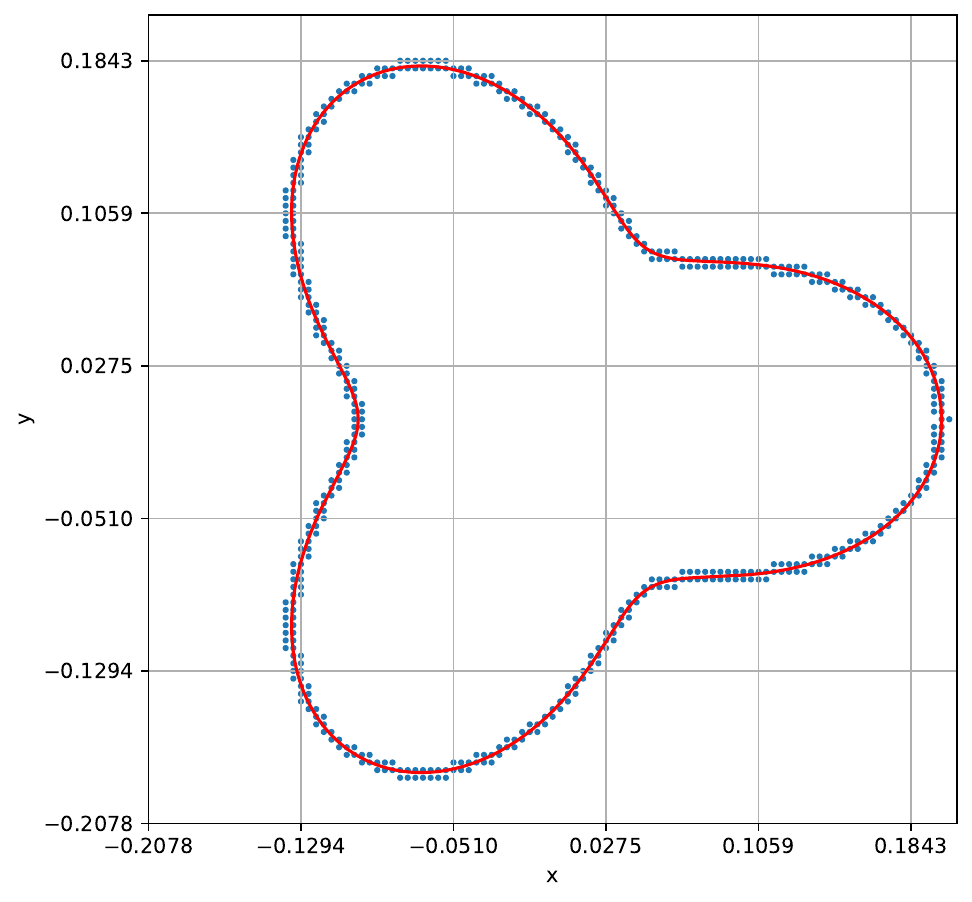}
		\caption{Smooth flower}
		\label{fig.flower.u.smooth}
	\end{subfigure}\hfill%
	\begin{subfigure}[b]{0.4\textwidth}
		\includegraphics[width=\textwidth]{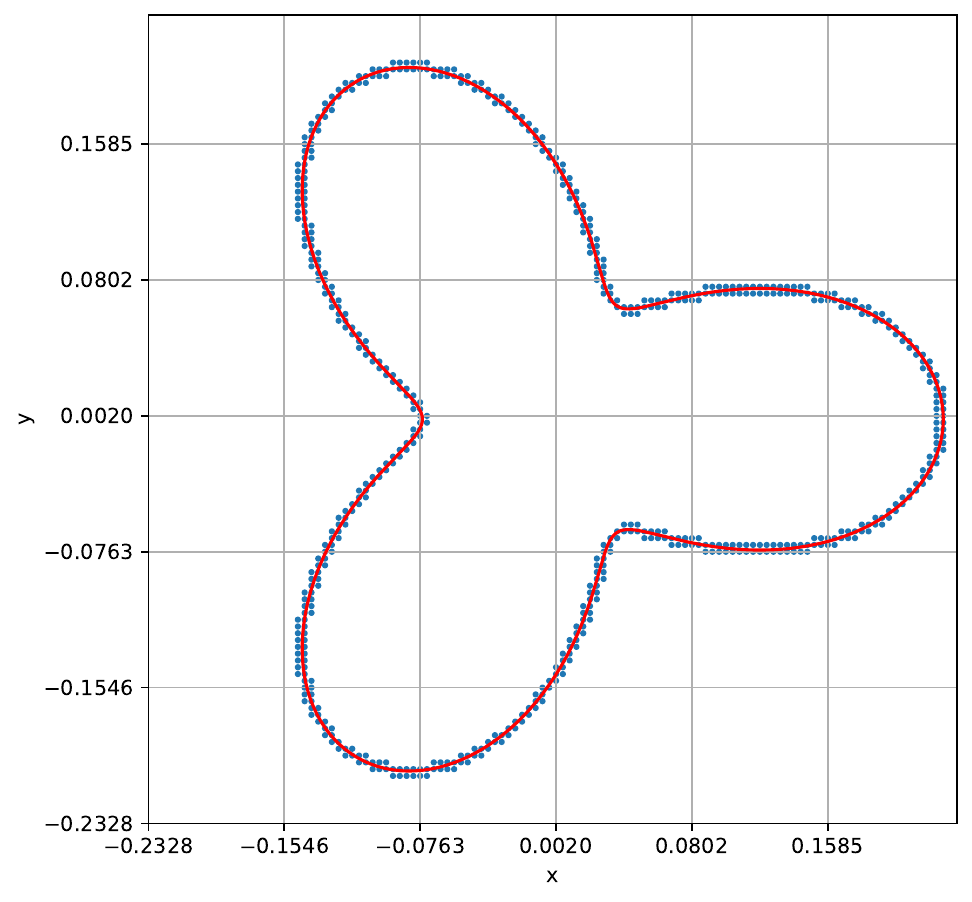}
		\caption{Acute flower}
		\label{fig.flower.u.acute}
	\end{subfigure}
	\caption{\small Flower-shaped interface embedded in low-resolution regular grids.}
	\label{fig.flower.u}
\end{figure}

We evaluate the accuracy of our level-set curvature neural networks on a two-dimensional irregular interface in both uniform and adaptive grids at different resolutions (see Figures \ref{fig.flower.u}, \ref{fig.flower.a.smooth}, and \ref{fig.flower.a.acute}).  The test interface corresponds to the zero level set of the bivariate function given by

\begin{equation}
\phi(x,y) = r(\theta) - a\cos{\left(p\theta\right)} - b,
\label{eq.implicitFlower}
\end{equation} 
where $r(\theta) = \sqrt{x^2 + y^2}$, $\theta \in [0, 2\pi)$ is the angle of $(x,y) \in \Omega$ with respect to the horizontal, and $a, b, p \in \mathbb{R}$ are shape parameters. For this interface, we can find the mean curvature analytically with the expression

\begin{equation}
\kappa(\theta) = \frac{r^2(\theta) + 2\left(r'(\theta)\right)^2 - r(\theta)r''(\theta)}{\left(r^2(\theta) + \left(r'(\theta)\right)^2\right)^{3/2}}.
\label{eq.kappaFlower}
\end{equation}

We have set $a = 0.05$ and $a = 0.075$ (keeping $b = 0.15$ and $p = 3$ constant) to generate a three-petaled interface with two steepness variations.  We refer to the first configuration as the \textit{smooth flower}, $\Gamma_s$, and to the second one as the \textit{acute flower}, $\Gamma_a$.  

Given a node with coordinates $\vv{x}_{i,j} = (x_i, y_j)$, we have trained our neural networks to infer the dimensionless curvature at the \emph{closest point}, $\vv{x}_{i,j}^\perp = (x_i^\perp, y_j^\perp)$, on the interface (see Figure \ref{fig.stencil}).  To perform an adequate comparison with the numerical approach, we first employ second-order finite differences to compute $h\kappa$ with \eqref{eq.numericalCurvature} everywhere on $\Omega$.  Then, we use $\phi$ to interpolate (bilinearly) the mean curvature to the location $\vv{x}_{i,j}^* = \vv{x}_{i,j} - \phi(\vv{x}_{i,j})\frac{\nabla \phi(\vv{x}_{i,j})}{|\nabla \phi(\vv{x}_{i,j})|}$.  The latter is the best approximation to $\vv{x}_{i,j}^\perp$ in view of the information supplied by the reinitialized level-set function.  In what follows, we refer to the finite-difference approximation of \eqref{eq.numericalCurvature} succeeded by bilinear interpolation as the \textit{compound numerical method}.  

The next series of analyses supports the thesis that we can build a resolution-based dictionary of neural models to estimate the mean curvature of implicit interfaces.  In our experiments, we generated samples by using 5, 10, and 20 iterations to reinitialize the level-set function \eqref{eq.implicitFlower}.  These samples were extracted from grids whose resolutions were equivalent to uniform mesh sizes of $h \sim 1/255$, $h \sim 1/265$, and $h \sim 1/275$.  We refer to these as \textit{low}, \textit{medium}, and \textit{high} resolutions, respectively.  They are intuitively matched to each of the trained neural networks described in section \ref{sec.buildingNNets}.

\subsection{Low-resolution regular grid}

We begin the accuracy analyses of our deep learning approach with the smooth flower interface, $\Gamma_s$, in a low-resolution regular grid.  We define the domain $\Omega \equiv [-0.207843, 0.207843]^2$ and discretize it into 107 equally spaced nodes along the Cartesian directions (see Figure \ref{fig.flower.u.smooth}).  This yields $h = 3.921569\times 10^{-3}$, which is equivalent to a $256 \times 256$ grid-point unit square.  Thus, we can use the trained neural network for $\rho = 256$.  Then, we collect 528 samples whose target $h\kappa$ values are calculated with \eqref{eq.kappaFlower}, where $\theta = \theta(x_i^\perp, y_j^\perp)$ and $(x_i^\perp, y_j^\perp)$ is the normal projection of node $(i,j)$ onto $r(\theta)$.

\begin{figure}[t]
	\centering
	\begin{subfigure}[b]{0.3\textwidth}
		\includegraphics[width=\textwidth]{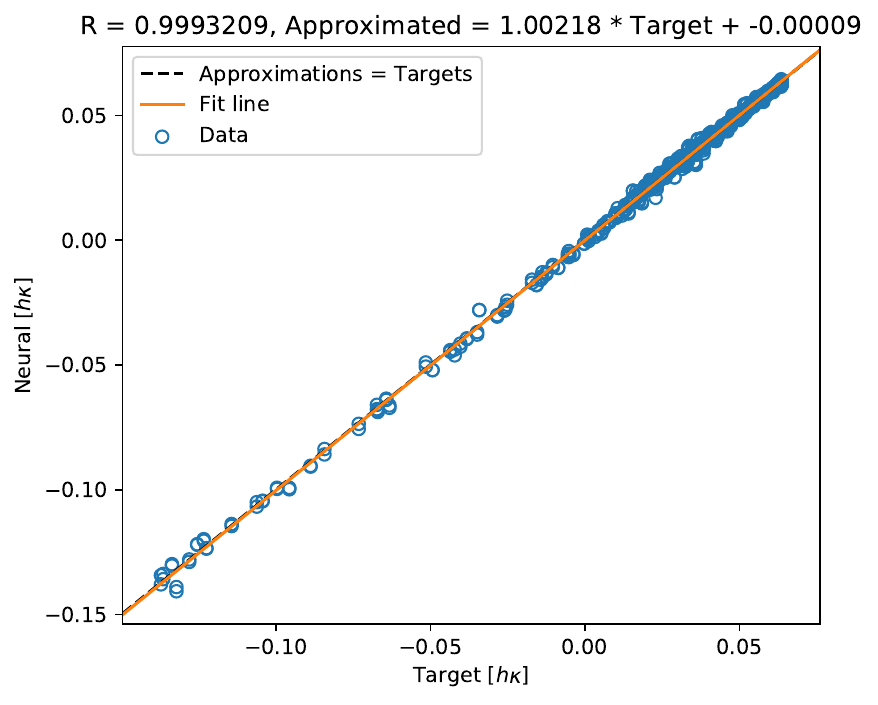}
        \caption{Neural, 5 iterations}
        \label{fig.flower.u.smooth.lowerRes.nnet.iter5}
    \end{subfigure}
	\begin{subfigure}[b]{0.3\textwidth}
		\includegraphics[width=\textwidth]{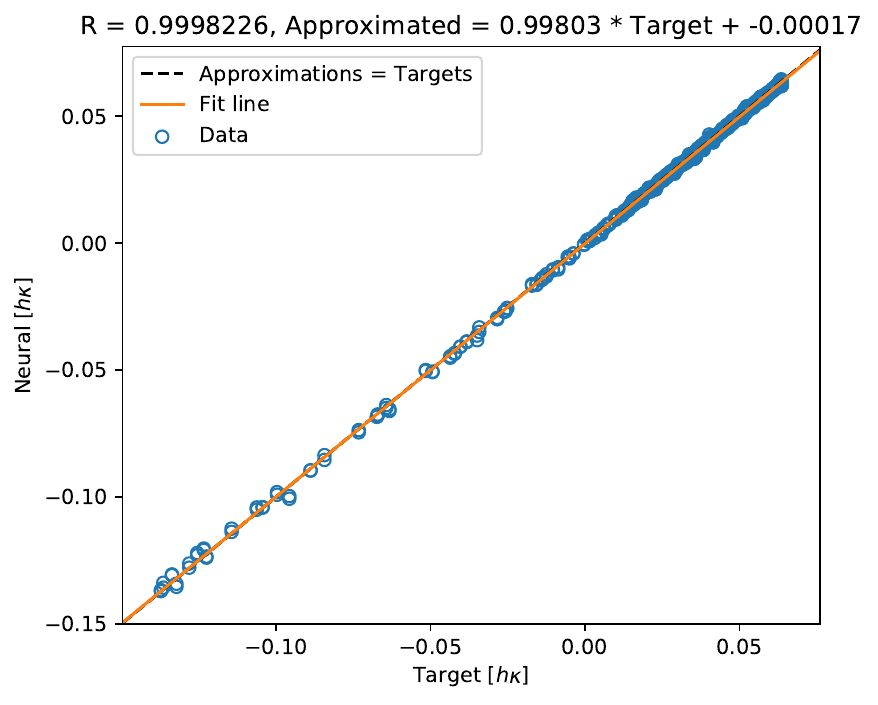}
		\caption{Neural, 10 iterations}
		\label{fig.flower.u.smooth.lowerRes.nnet.iter10}
	\end{subfigure}
	\begin{subfigure}[b]{0.3\textwidth}
		\includegraphics[width=\textwidth]{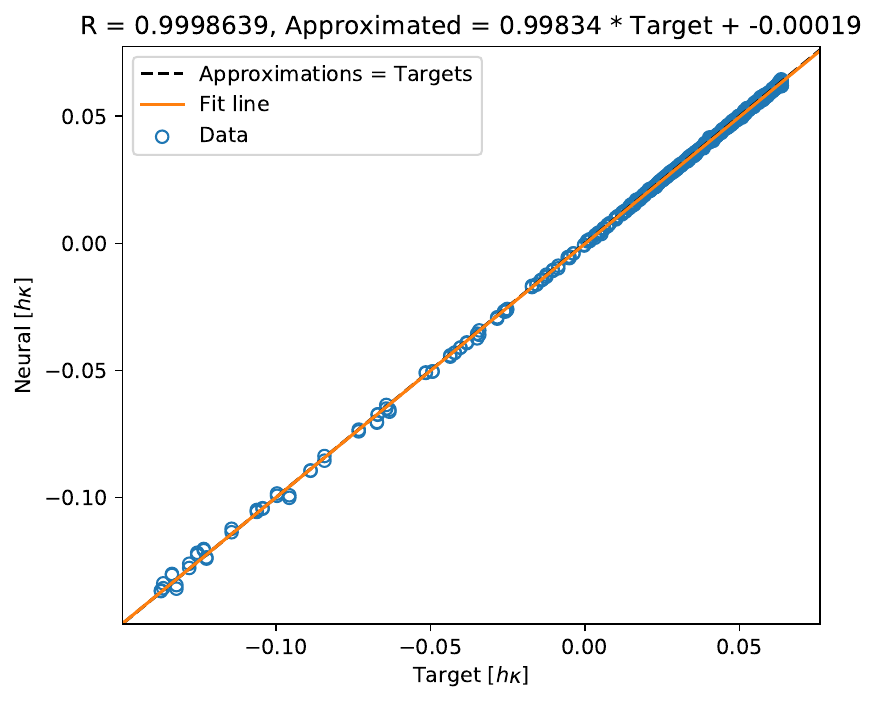}
		\caption{Neural, 20 iterations}
		\label{fig.flower.u.smooth.lowerRes.nnet.iter20}
	\end{subfigure}
   
	\begin{subfigure}[b]{0.3\textwidth}
		\includegraphics[width=\textwidth]{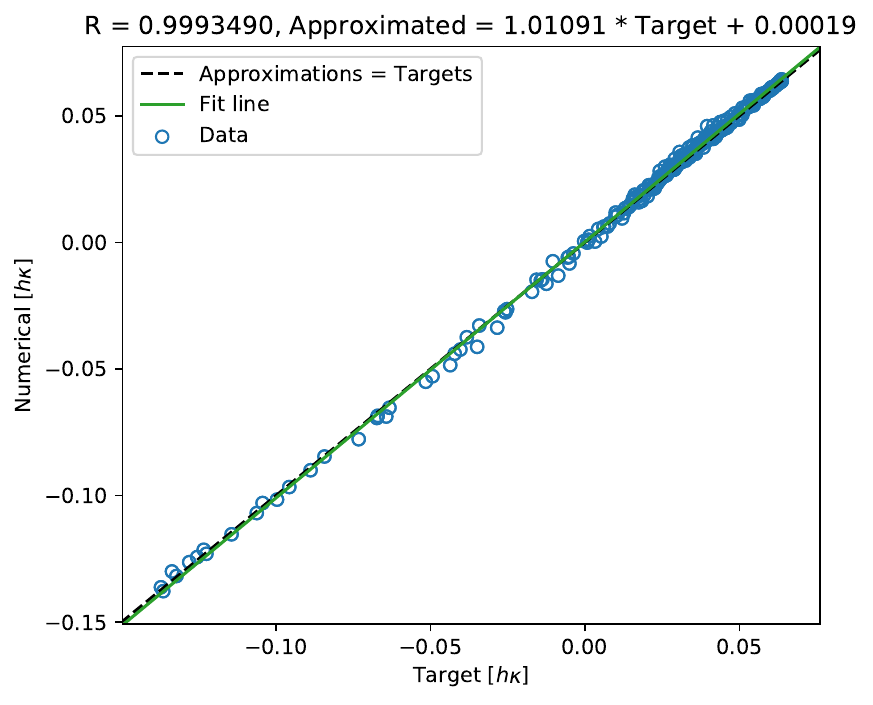}
		\caption{Numerical, 5 iterations}
		\label{fig.flower.u.smooth.lowerRes.numerics.iter5}
	\end{subfigure}
	\begin{subfigure}[b]{0.3\textwidth}
		\includegraphics[width=\textwidth]{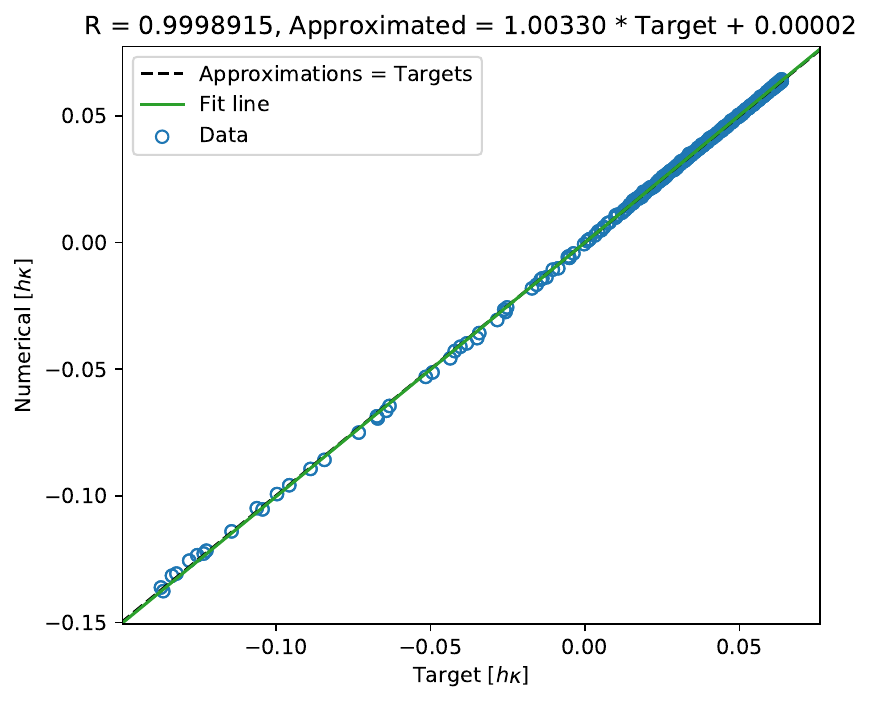}
		\caption{Numerical, 10 iterations}
		\label{fig.flower.u.smooth.lowerRes.numerics.iter10}
	\end{subfigure}
	\begin{subfigure}[b]{0.3\textwidth}
		\includegraphics[width=\textwidth]{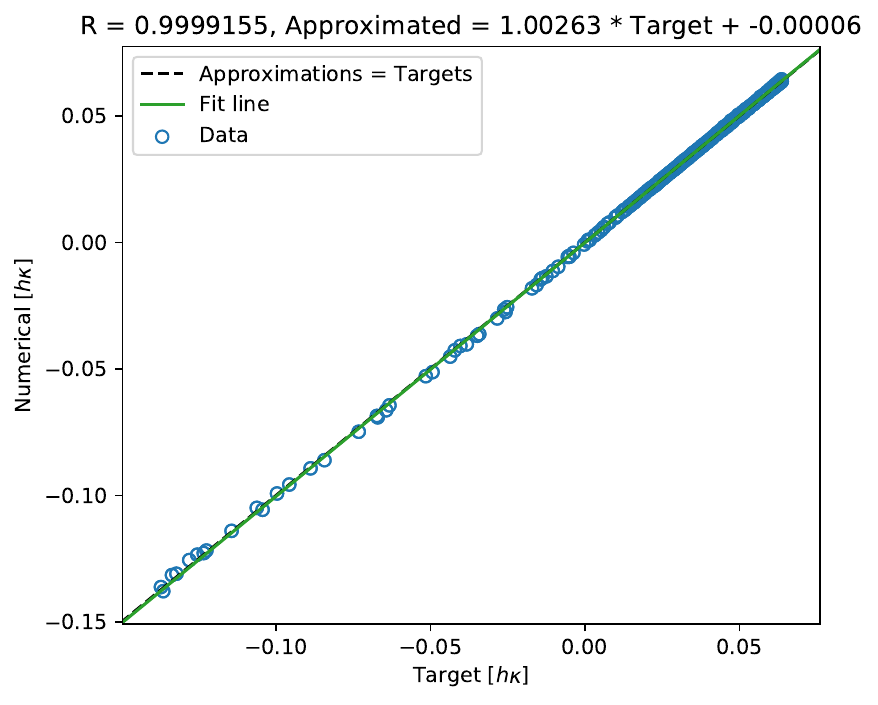}
		\caption{Numerical, 20 iterations}
		\label{fig.flower.u.smooth.lowerRes.numerics.iter20}
	\end{subfigure}
	\caption{\small Correlation between expected and approximated $h\kappa$ using the neural network and the numerical method for the smooth flower interface in a regular grid of $107 \times 107$ nodes.}
	\label{fig.flower.u.smooth.lowRes.correlation}
\end{figure}

Figure \ref{fig.flower.u.smooth.lowRes.correlation} illustrates the quality of the neural and numerical $h\kappa$ approximations as one varies the number of iterations to reinitialize \eqref{eq.implicitFlower}.  These results visually confirm that our deep learning estimations are in good agreement with those computed with conventional numerical schemes.  Unlike the numerical method, however, we remark that our neural network does not require any subsequent interpolation to calculate $h\kappa$ at the interface.  

Table \ref{tbl.flower.u.smooth.lowRes.errors} provides some statistics for the current analysis of $\Gamma_s$.  These results support our previous statement that the deep learning approach performs, on average, as well as the numerical technique.  Specifically, the MAE is not significantly distinct between both methodologies.  When one uses 5 iterations to reinitialize the level-set function, for example, the neural network performs better than the numerical approach in the error $L^1$ and $L^2$ norms.  In terms of the $L^\infty$ norm, however, the finite-difference schemes and subsequent bilinear interpolation are more robust and improve up to 43\% of the corresponding error incurred by the neural network.

\begin{table}[!b]
	\centering
	\footnotesize
	\bgroup
	\def\arraystretch{1.1}%
	\begin{tabular}{|l|l|r|r|r|}
		\hline
		Iterations & Method & MAE & Max AE & MSE \\
		\hline \hline
		\multirow{2}{*}{5} & Neural & $1.142733\times 10^{-3}$ & $8.479774\times 10^{-3}$ & $2.589772\times 10^{-6}$ \\
 		& Numerical & $1.205349\times 10^{-3}$ & $6.373532\times 10^{-3}$ & $2.920415\times 10^{-6}$ \\
		\hline \hline
		\multirow{2}{*}{10} & Neural & $5.850180\times 10^{-4}$ & $5.003561\times 10^{-3}$ & $7.187668\times 10^{-7}$ \\
 		& Numerical & $4.324966\times 10^{-4}$ & $2.860623\times 10^{-3}$ & $4.412982\times 10^{-7}$ \\
		\hline \hline
		\multirow{2}{*}{20} & Neural & $4.772196\times 10^{-4}$ & $4.472834\times 10^{-3}$ & $5.695442\times 10^{-7}$ \\
 		& Numerical & $3.166947\times 10^{-4}$ & $2.578525\times 10^{-3}$ & $3.341286\times 10^{-7}$ \\
		\hline
	\end{tabular}
	\egroup
	\caption{\small Error analysis for the smooth flower interface in a regular grid of $107 \times 107$ nodes.}
	\label{tbl.flower.u.smooth.lowRes.errors}
\end{table}

\begin{figure}[t]
	\centering
	\begin{subfigure}[b]{0.3\textwidth}
		\includegraphics[width=\textwidth]{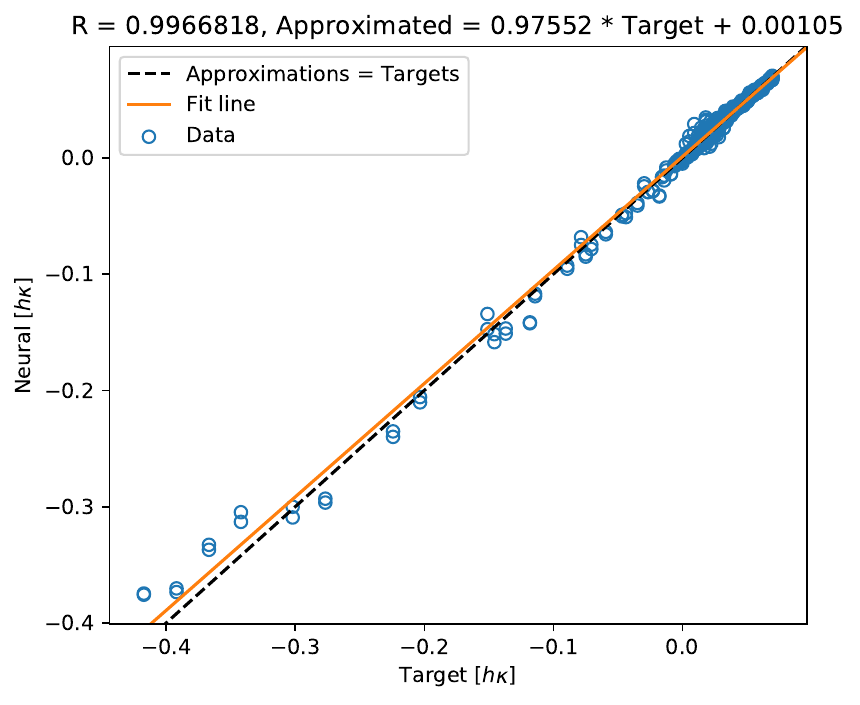}
        \caption{Neural, 5 iterations}
        \label{fig.flower.u.acute.lowerRes.nnet.iter5}
    \end{subfigure}
	\begin{subfigure}[b]{0.3\textwidth}
		\includegraphics[width=\textwidth]{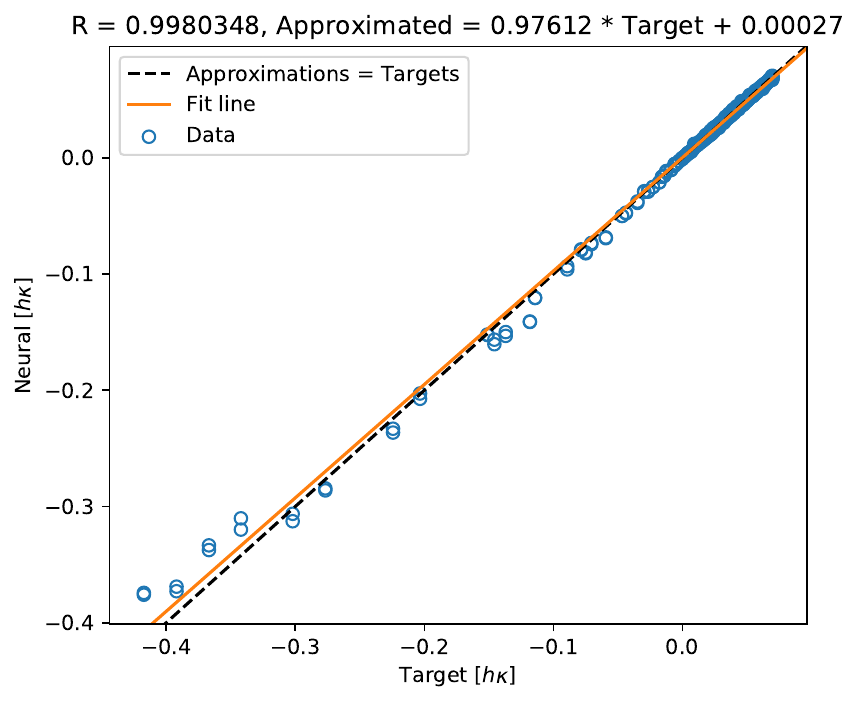}
		\caption{Neural, 10 iterations}
		\label{fig.flower.u.acute.lowerRes.nnet.iter10}
	\end{subfigure}
	\begin{subfigure}[b]{0.3\textwidth}
		\includegraphics[width=\textwidth]{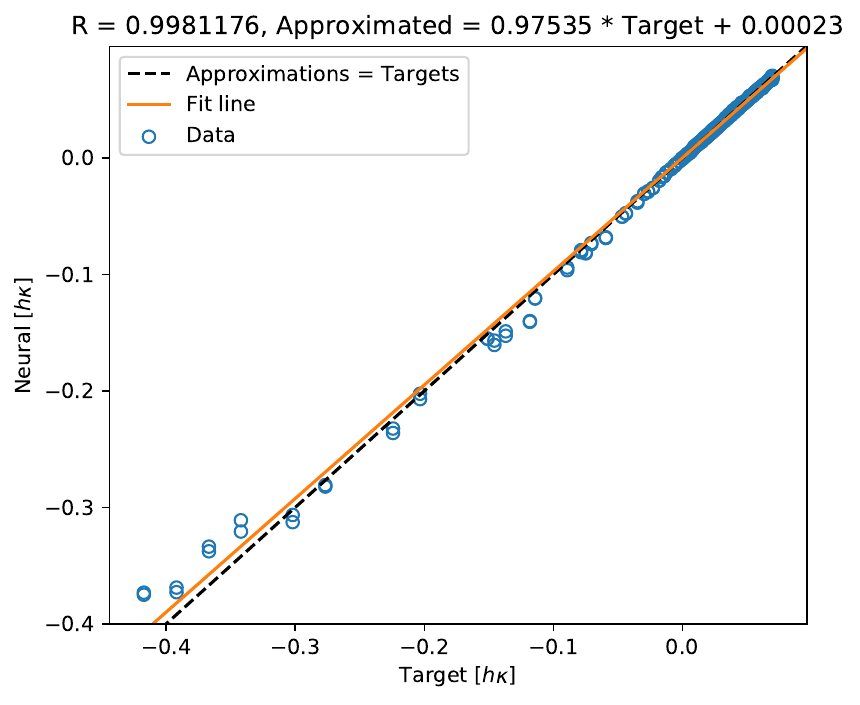}
		\caption{Neural, 20 iterations}
		\label{fig.flower.u.acute.lowerRes.nnet.iter20}
	\end{subfigure}
    
	\begin{subfigure}[b]{0.3\textwidth}
		\includegraphics[width=\textwidth]{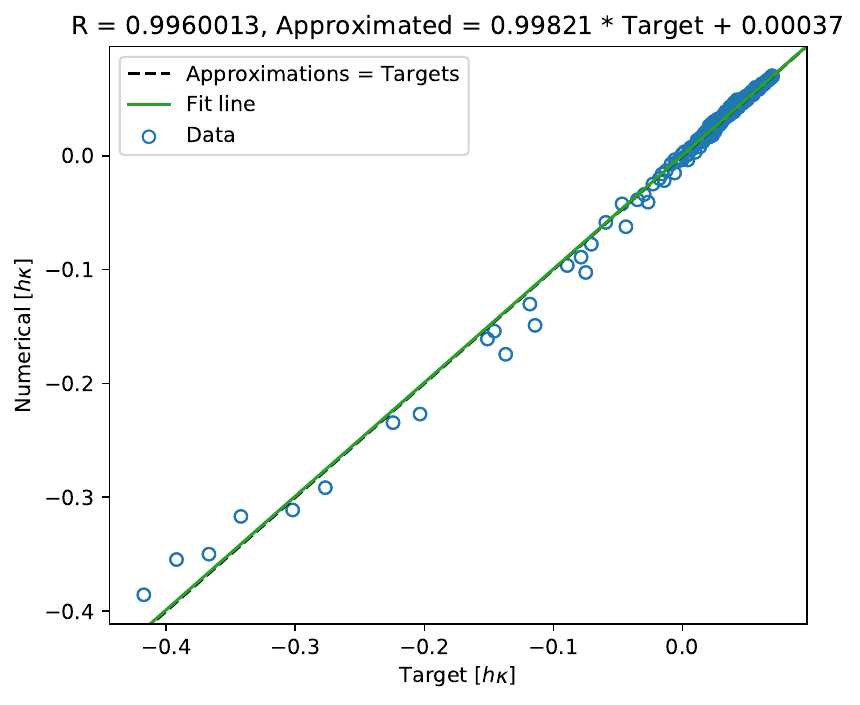}
		\caption{Numerical, 5 iterations}
		\label{fig.flower.u.acute.lowerRes.numerics.iter5}
	\end{subfigure}
	\begin{subfigure}[b]{0.3\textwidth}
		\includegraphics[width=\textwidth]{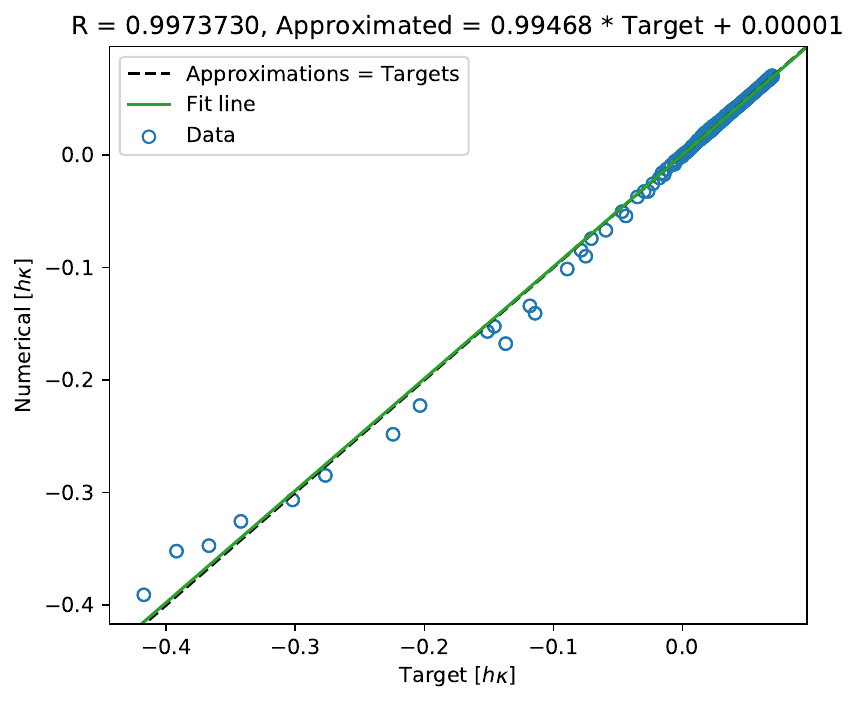}
		\caption{Numerical, 10 iterations}
		\label{fig.flower.u.acute.lowerRes.numerics.iter10}
	\end{subfigure}
	\begin{subfigure}[b]{0.3\textwidth}
		\includegraphics[width=\textwidth]{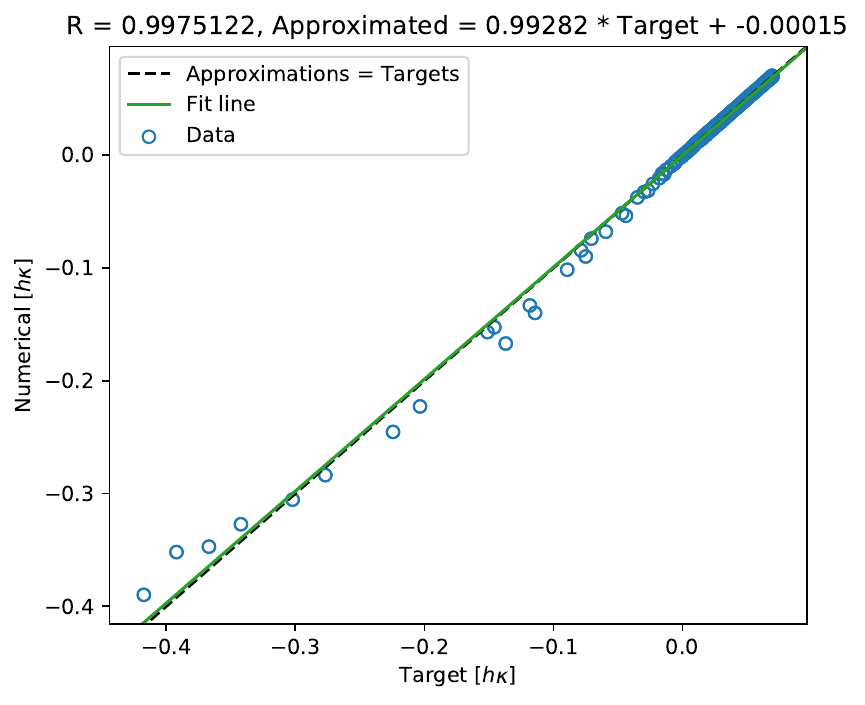}
		\caption{Numerical, 20 iterations}
		\label{fig.flower.u.acute.lowerRes.numerics.iter20}
	\end{subfigure}
	\caption{\small Correlation between expected and approximated $h\kappa$ using the neural network and the numerical method for the acute flower interface in a regular grid of $120 \times 120$ nodes.}
	\label{fig.flower.u.acute.lowRes.correlation}
\end{figure}

Next, we analyze the curvature approximations for the acute interface, $\Gamma_a$, in a low-resolution regular grid.  For this, we set the domain $\Omega \equiv [-0.232826, 0.232826]^2$ and discretize it into 120 equally spaced vertices along each Cartesian direction (see Figure \ref{fig.flower.u.acute}).  This discretization yields $h = 3.913043\times 10^{-3}$, which corresponds to a unit-square resolution of $256.56$ nodes per unit length.  Such a mesh size is compatible with the neural model trained for $\rho = 256$ and allows us to collect 624 samples along $\Gamma_a$. Figure \ref{fig.flower.u.acute.lowRes.correlation} shows the quality of the neural and numerical $h\kappa$ approximations.

Compared to the previous experiment, the steeper curvatures at the rapidly turning petal junctions in $\Gamma_a$ present major difficulties for both methods.  Our model, however, once again produces similar curvature values to the approximations calculated with the compound numerical approach.  In fact, as seen in Figure \ref{fig.flower.u.acute.lowRes.correlation}, our strategy exhibits higher correlation factors regardless of the number of iterations used for level-set redistancing.  Table \ref{tbl.flower.u.acute.lowRes.errors} supports these findings by showing that our deep learning framework is better than the numerical scheme when contrasting the error $L^1$ and $L^2$ norms.  As with the results provided in Table \ref{tbl.flower.u.smooth.lowRes.errors}, the two-step numerical approach still proves to be more effective at reducing the error $L^\infty$ norm.  Nevertheless, these maximum absolute errors are not significantly improved beyond 12\% of the corresponding neural network errors.

\begin{table}[!b]
	\centering
	\footnotesize
	\bgroup
	\def\arraystretch{1.1}%
	\begin{tabular}{|l|l|r|r|r|}
		\hline
		Iterations & Method & MAE & Max AE & MSE \\
		\hline \hline
		\multirow{2}{*}{5} & Neural & $2.731657\times 10^{-3}$ & $4.264017\times 10^{-2}$ & $2.897127\times 10^{-5}$ \\
 		& Numerical & $2.782704\times 10^{-3}$ & $3.743515\times 10^{-2}$ & $3.317569\times 10^{-5}$ \\
		\hline \hline
		\multirow{2}{*}{10} & Neural & $1.319063\times 10^{-3}$ & $4.306652\times 10^{-2}$ & $1.786897\times 10^{-5}$ \\
 		& Numerical & $1.438051\times 10^{-3}$ & $3.979733\times 10^{-2}$ & $2.163957\times 10^{-5}$ \\
		\hline \hline
		\multirow{2}{*}{20} & Neural & $1.165201\times 10^{-3}$ & $4.412831\times 10^{-2}$ & $1.736991\times 10^{-5}$ \\
 		& Numerical & $1.183391\times 10^{-3}$ & $4.001732\times 10^{-2}$ & $2.059011\times 10^{-5}$ \\
		\hline
	\end{tabular}
	\egroup
	\caption{\small Error analysis for the acute flower interface in a regular grid of $120 \times 120$ nodes.}
	\label{tbl.flower.u.acute.lowRes.errors}
\end{table}

\subsection{Medium-resolution regular grid}

We continue with our experiments by assessing the $h\kappa$ approximations for the smooth interface, $\Gamma_s$, in a medium-resolution regular grid.  Following the protocol outlined above, we start with a computational domain $\Omega \equiv [-0.207547, 0.207547]^2$ and discretize it into 111 uniform nodes along each Cartesian direction.  This yields a cell's width of $h = 3.773585\times 10^{-3}$, which is equivalent to a $266 \times 266$-node unit square and is compatible with our neural model for $\rho = 266$.  Also, the total number of samples we collect along the interface is 552.  Since this interface has the same shape as the smooth flower in a low-resolution mesh, we refer the reader to Figure \ref{fig.flower.u.smooth} for a visual representation of $\Gamma_s$.

\begin{figure}[t]
	\centering
	\begin{subfigure}[b]{0.3\textwidth}
		\includegraphics[width=\textwidth]{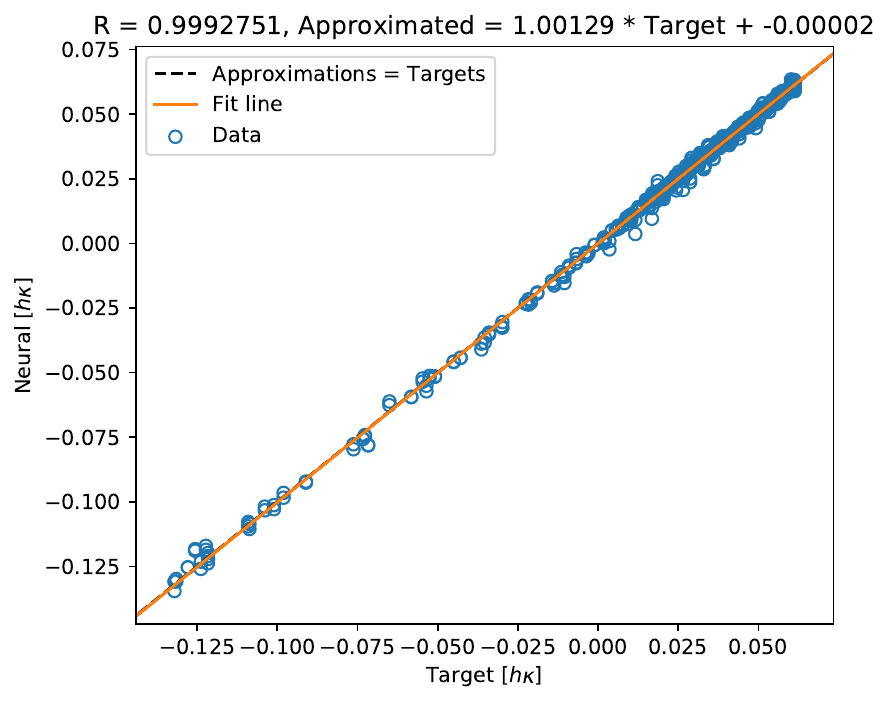}
        \caption{Neural, 5 iterations}
        \label{fig.flower.u.smooth.medRes.nnet.iter5}
    \end{subfigure}
	\begin{subfigure}[b]{0.3\textwidth}
		\includegraphics[width=\textwidth]{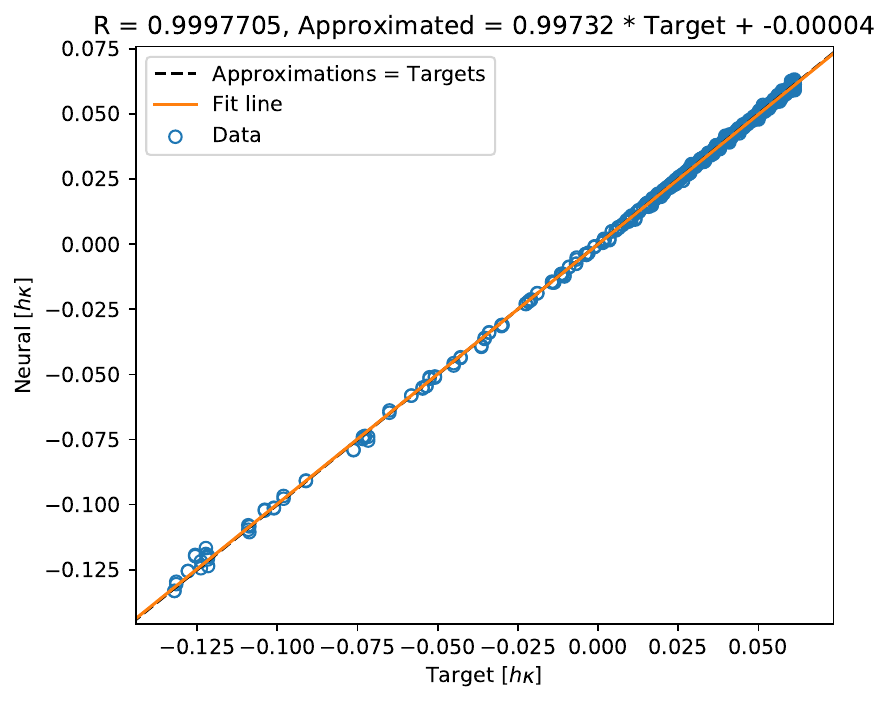}
		\caption{Neural, 10 iterations}
		\label{fig.flower.u.smooth.medRes.nnet.iter10}
	\end{subfigure}
	\begin{subfigure}[b]{0.3\textwidth}
		\includegraphics[width=\textwidth]{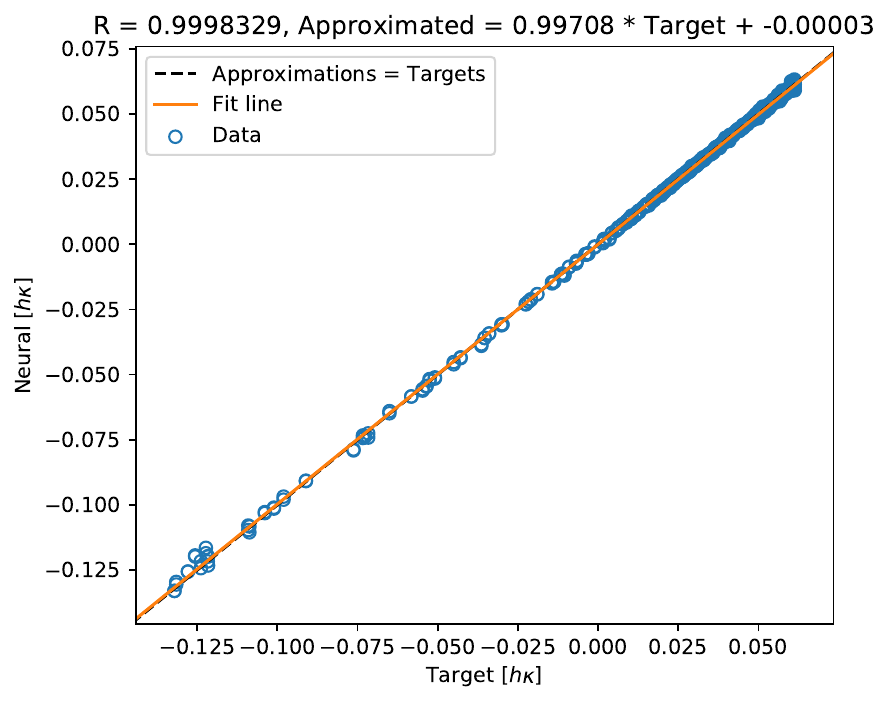}
		\caption{Neural, 20 iterations}
		\label{fig.flower.u.smooth.medRes.nnet.iter20}
	\end{subfigure}
    
	\begin{subfigure}[b]{0.3\textwidth}
		\includegraphics[width=\textwidth]{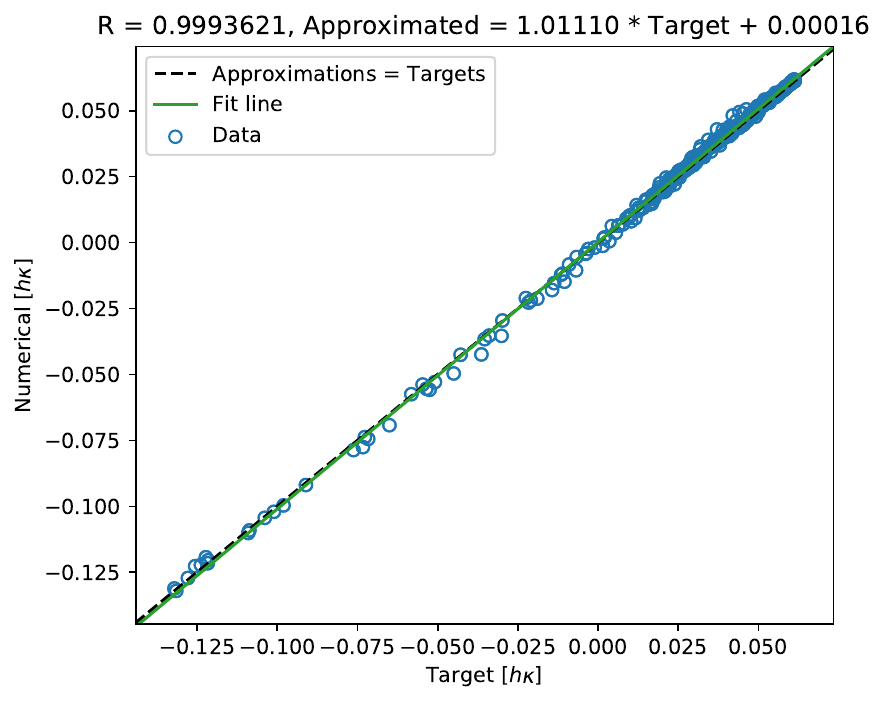}
		\caption{Numerical, 5 iterations}
		\label{fig.flower.u.smooth.medRes.numerics.iter5}
	\end{subfigure}
	\begin{subfigure}[b]{0.3\textwidth}
		\includegraphics[width=\textwidth]{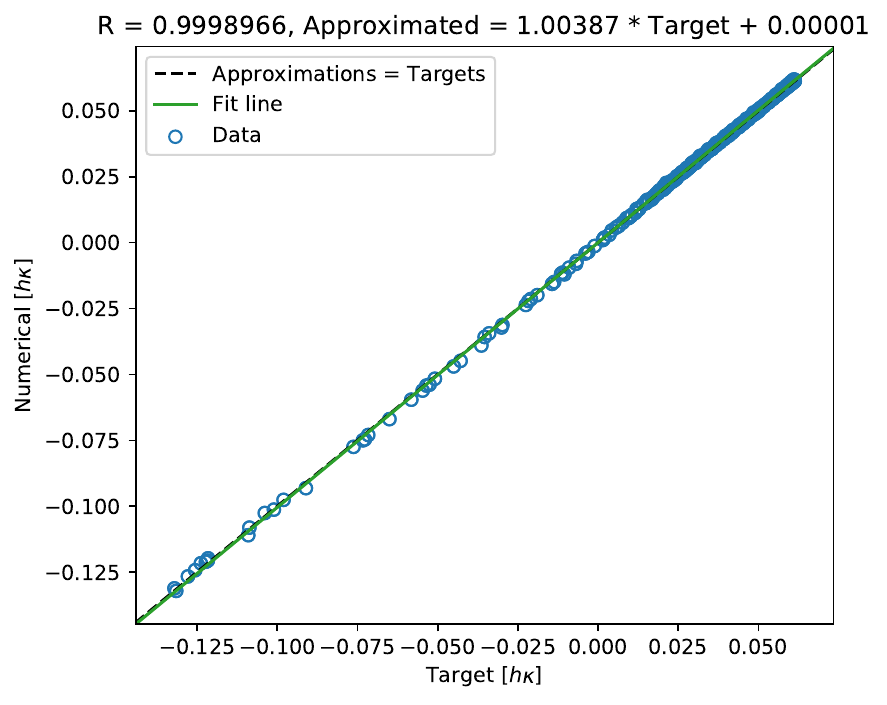}
		\caption{Numerical, 10 iterations}
		\label{fig.flower.u.smooth.medRes.numerics.iter10}
	\end{subfigure}
	\begin{subfigure}[b]{0.3\textwidth}
		\includegraphics[width=\textwidth]{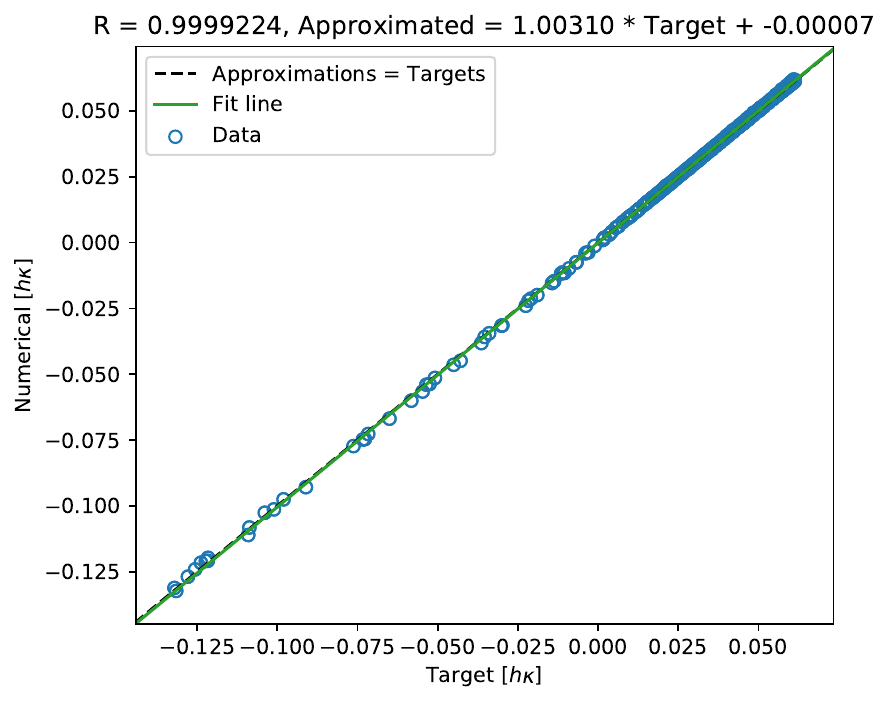}
		\caption{Numerical, 20 iterations}
		\label{fig.flower.u.smooth.medRes.numerics.iter20}
	\end{subfigure}
	\caption{\small Correlation between expected and approximated $h\kappa$ using the neural network and the numerical method for the smooth flower interface in a regular grid of $111 \times 111$ nodes.}
	\label{fig.flower.u.smooth.medRes.correlation}
\end{figure}

Figure \ref{fig.flower.u.smooth.medRes.correlation} illustrates the quality of the $h\kappa$ approximations for the current case study.  As in the previous analysis of $\Gamma_s$ (see Figure \ref{fig.flower.u.smooth.lowRes.correlation}), the neural network correlation factors are again in good agreement with the estimations from the numerical method.  As shown in Table \ref{tbl.flower.u.smooth.medRes.errors}, our proposed strategy produces MAEs that are not significantly different from those incurred by the compound numerical approach.  However, we note that in terms of the error $L^\infty$ norm, the finite-difference schemes followed by bilinear interpolation can reduce as much as 63\% of the maximum error experienced with the deep learning strategy.  We remark that such a notable improvement in the numerical method characterizes the estimation process for the smooth flower interface, $\Gamma_s$, regardless of the discretization considered in this section.

\begin{table}[!b]
	\centering
	\footnotesize
	\bgroup
	\def\arraystretch{1.1}%
	\begin{tabular}{|l|l|r|r|r|}
		\hline
		Iterations & Method & MAE & Max AE & MSE \\
		\hline \hline
		\multirow{2}{*}{5} & Neural & $1.111979\times 10^{-3}$ & $7.959664\times 10^{-3}$ & $2.584322\times 10^{-6}$ \\
 		& Numerical & $1.137535\times 10^{-3}$ & $6.195024\times 10^{-3}$ & $2.681028\times 10^{-6}$ \\
		\hline \hline
		\multirow{2}{*}{10} & Neural & $5.838719\times 10^{-4}$ & $6.307398\times 10^{-3}$ & $8.313230\times 10^{-7}$ \\
 		& Numerical & $4.180538\times 10^{-4}$ & $2.631675\times 10^{-3}$ & $4.035943\times 10^{-7}$ \\
		\hline \hline
		\multirow{2}{*}{20} & Neural & $4.642169\times 10^{-4}$ & $6.246050\times 10^{-3}$ & $6.122834\times 10^{-7}$ \\
 		& Numerical & $3.014267\times 10^{-4}$ & $2.294022\times 10^{-3}$ & $2.942238\times 10^{-7}$ \\
		\hline
	\end{tabular}
	\egroup
	\caption{\small Error analysis for the smooth flower interface in a regular grid of $111 \times 111$ nodes.}
	\label{tbl.flower.u.smooth.medRes.errors}
\end{table}

Next, we examine the accuracy of our deep learning approach at approximating $h\kappa$ for the acute interface in a medium-resolution regular grid.  Following the previous procedure, we begin by defining $\Omega \equiv [-0.232563, 0.232563]^2$ and  discretizing it into 124 uniform nodes along each Cartesian direction.  Thus, we get a mesh size $h = 3.781513\times 10^{-3}$, which makes the regular domain equivalent to a unit square with 265.44 nodes per side.  With these settings we collect 648 samples next to $\Gamma_a$ (see Figure \ref{fig.flower.u.acute}) and satisfy the resolution scope for the model trained for $\rho = 266$.

\begin{figure}[t]
	\centering
	\begin{subfigure}[b]{0.3\textwidth}
		\includegraphics[width=\textwidth]{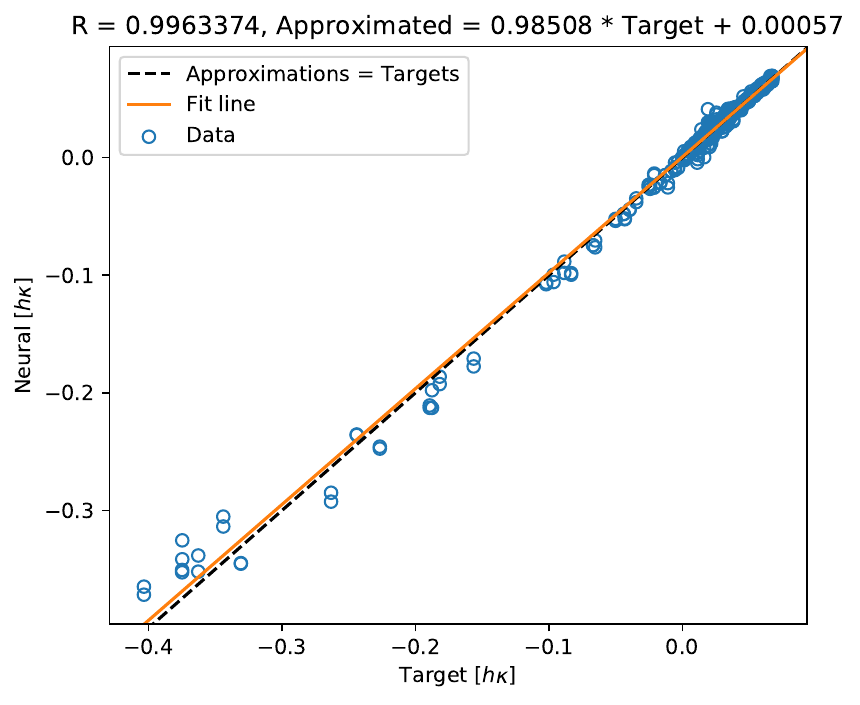}
        \caption{Neural, 5 iterations}
        \label{fig.flower.u.acute.medRes.nnet.iter5}
    \end{subfigure}
	\begin{subfigure}[b]{0.3\textwidth}
		\includegraphics[width=\textwidth]{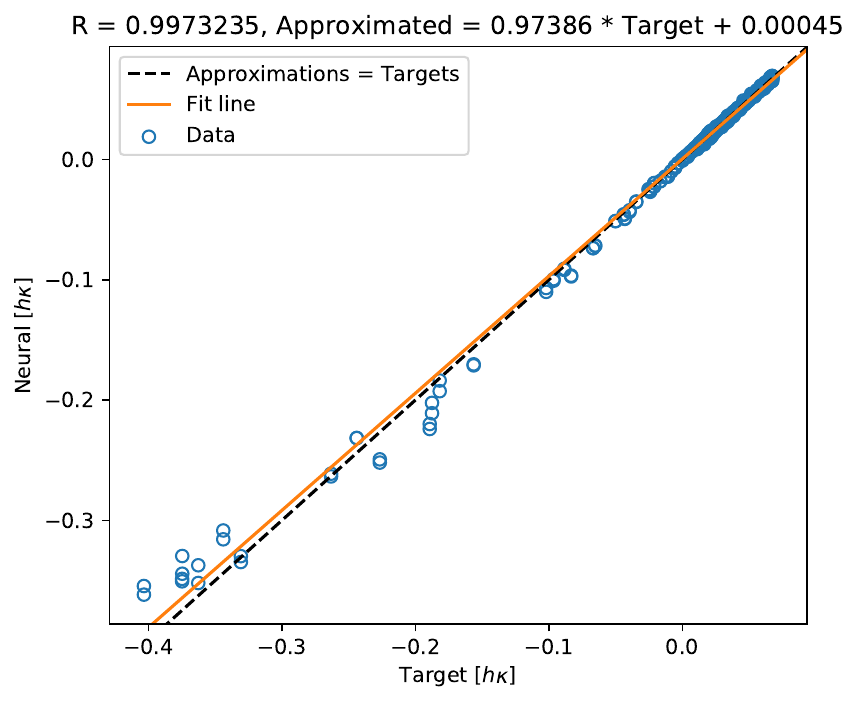}
		\caption{Neural, 10 iterations}
		\label{fig.flower.u.acute.medRes.nnet.iter10}
	\end{subfigure}
	\begin{subfigure}[b]{0.3\textwidth}
		\includegraphics[width=\textwidth]{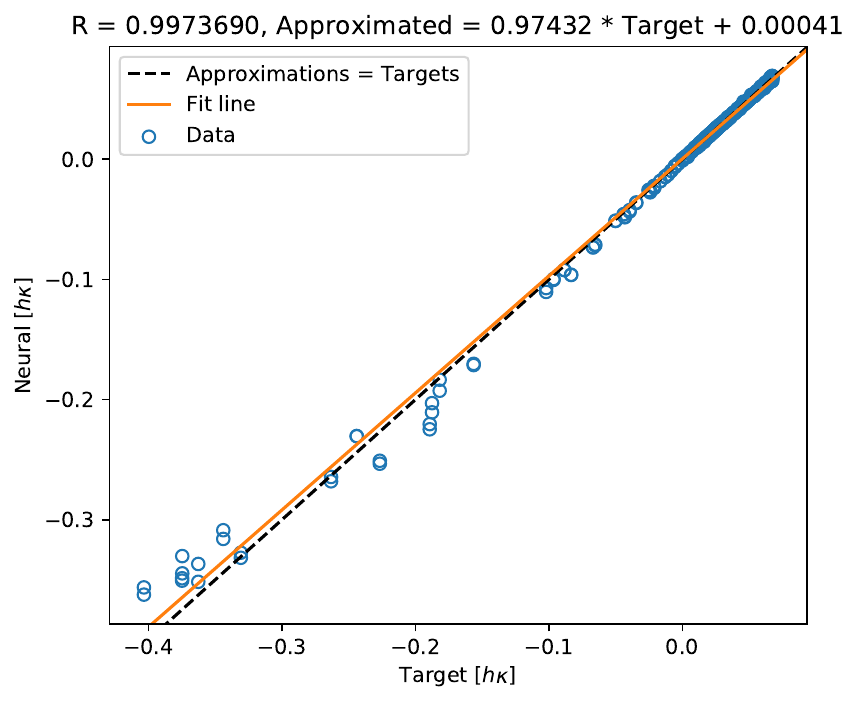}
		\caption{Neural, 20 iterations}
		\label{fig.flower.u.acute.medRes.nnet.iter20}
	\end{subfigure}
    
	\begin{subfigure}[b]{0.3\textwidth}
		\includegraphics[width=\textwidth]{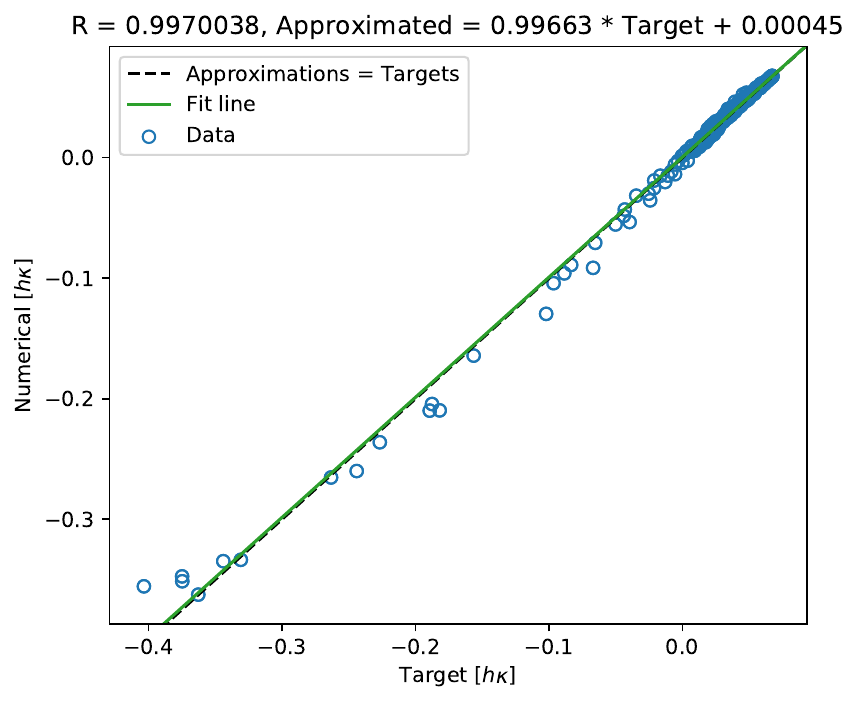}
		\caption{Numerical, 5 iterations}
		\label{fig.flower.u.acute.medRes.numerics.iter5}
	\end{subfigure}
	\begin{subfigure}[b]{0.3\textwidth}
		\includegraphics[width=\textwidth]{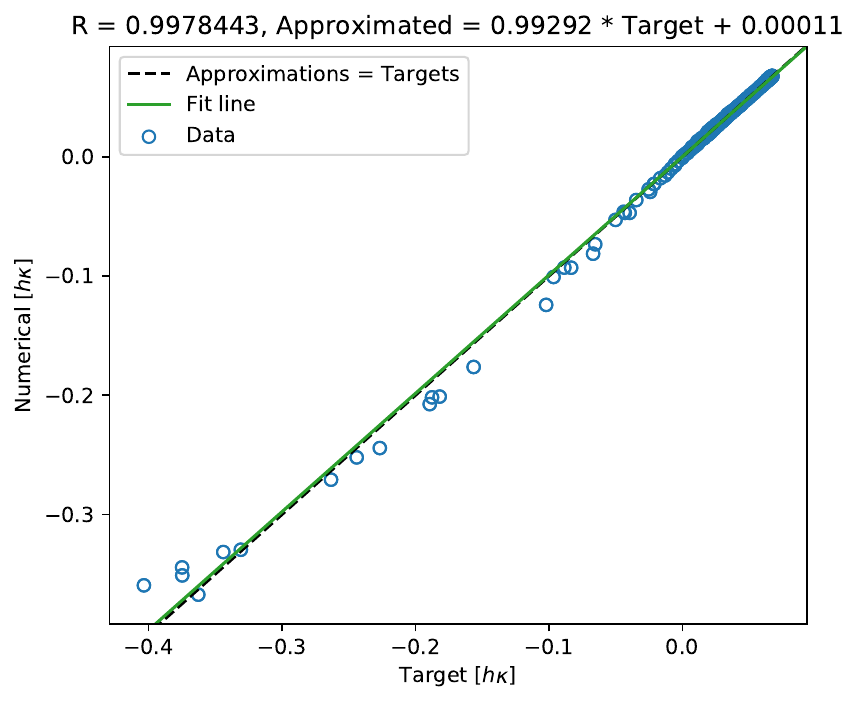}
		\caption{Numerical, 10 iterations}
		\label{fig.flower.u.acute.medRes.numerics.iter10}
	\end{subfigure}
	\begin{subfigure}[b]{0.3\textwidth}
		\includegraphics[width=\textwidth]{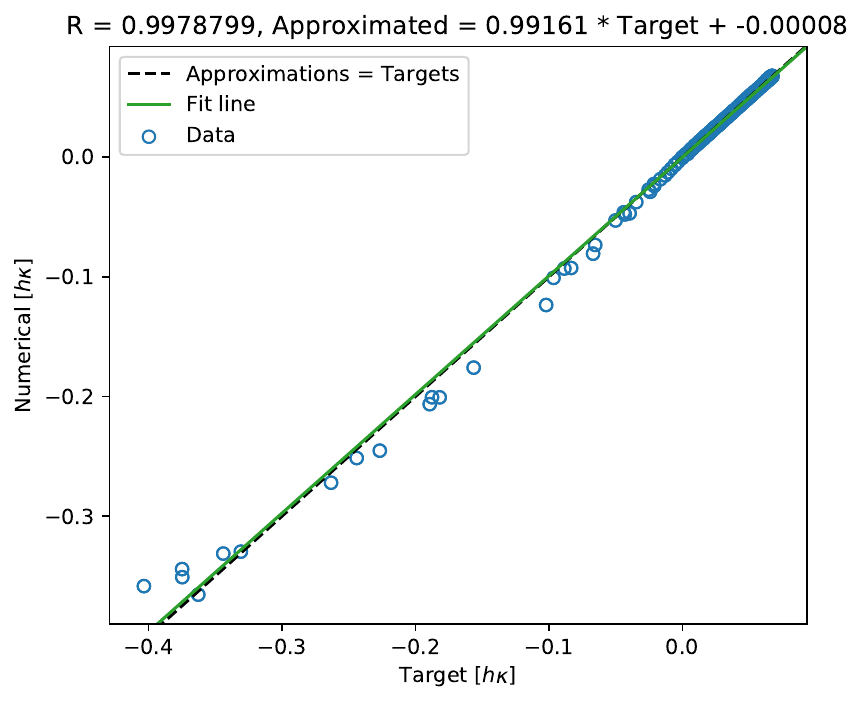}
		\caption{Numerical, 20 iterations}
		\label{fig.flower.u.acute.medRes.numerics.iter20}
	\end{subfigure}
	\caption{\small Correlation between expected and approximated $h\kappa$ using the neural network and the numerical method for the acute flower interface in a regular grid of $124 \times 124$ nodes.}
	\label{fig.flower.u.acute.medRes.correlation}
\end{figure}

Figure \ref{fig.flower.u.acute.medRes.correlation} depicts the quality of the curvature approximations for the current experiment.  A visual inspection reveals that our model outputs modestly accurate $h\kappa$ values despite the steep curvatures in $\Gamma_a$.  These estimations, albeit noisier than the numerical ones, have high correlation factors and improve as one increases the number of redistancing operations.  The statistical summary in Table \ref{tbl.flower.u.acute.medRes.errors} further corroborates that our neural network and the compound numerical method deliver practically the same level of precision in the $L^1$ and $L^2$ norms.  Likewise, the maximum absolute errors are not significantly distinct, and the numerical estimations experience only a minor improvement of up to 10\% over our strategy.  The latter confirms what is also clear in Table \ref{tbl.flower.u.acute.lowRes.errors}: the accuracy gap between both techniques diminishes when the interface features steep curvatures.  We repeatedly see this phenomenon throughout our experiments with uniform grids whenever $\Gamma_a$ is involved.

\begin{table}[!b]
	\centering
	\footnotesize
	\bgroup
	\def\arraystretch{1.1}%
	\begin{tabular}{|l|l|r|r|r|}
		\hline
		Iterations & Method & MAE & Max AE & MSE \\
		\hline \hline
		\multirow{2}{*}{5} & Neural & $2.851336\times 10^{-3}$ & $4.926337\times 10^{-2}$ & $3.416267\times 10^{-5}$ \\
 		& Numerical & $2.608772\times 10^{-3}$ & $4.774025\times 10^{-2}$ & $2.783543\times 10^{-5}$ \\
		\hline \hline
		\multirow{2}{*}{10} & Neural & $1.571223\times 10^{-3}$ & $4.877240\times 10^{-2}$ & $2.671261\times 10^{-5}$ \\
 		& Numerical & $1.378433\times 10^{-3}$ & $4.396094\times 10^{-2}$ & $1.993660\times 10^{-5}$ \\
		\hline \hline
		\multirow{2}{*}{20} & Neural & $1.378460\times 10^{-3}$ & $4.705131\times 10^{-2}$ & $2.622015\times 10^{-5}$ \\
 		& Numerical & $1.140213\times 10^{-3}$ & $4.493223\times 10^{-2}$ & $1.969833\times 10^{-5}$ \\
		\hline
	\end{tabular}
	\egroup
	\caption{\small Error analysis for the acute flower interface in a regular grid of $124 \times 124$ nodes.}
	\label{tbl.flower.u.acute.medRes.errors}
\end{table}

\subsection{High-resolution regular grid}

We now assess our deep learning strategy on the smooth interface, $\Gamma_s$, in a high-resolution regular grid (see Figure \ref{fig.flower.u.smooth}).  We begin by defining the computational domain $\Omega \equiv [-0.207339, 0.207339]^2$ and discretizing it into 114 equally spaced points along each Cartesian direction.  As a result, we collect 564 samples from a uniform grid with $h = 3.669724\times 10^{-3}$.  This cell's width is equivalent to a unit-square resolution of 273.5 nodes per side length and matches the training resolution of our neural network for $\rho = 276$.

\begin{figure}[t]
	\centering
	\begin{subfigure}[b]{0.3\textwidth}
		\includegraphics[width=\textwidth]{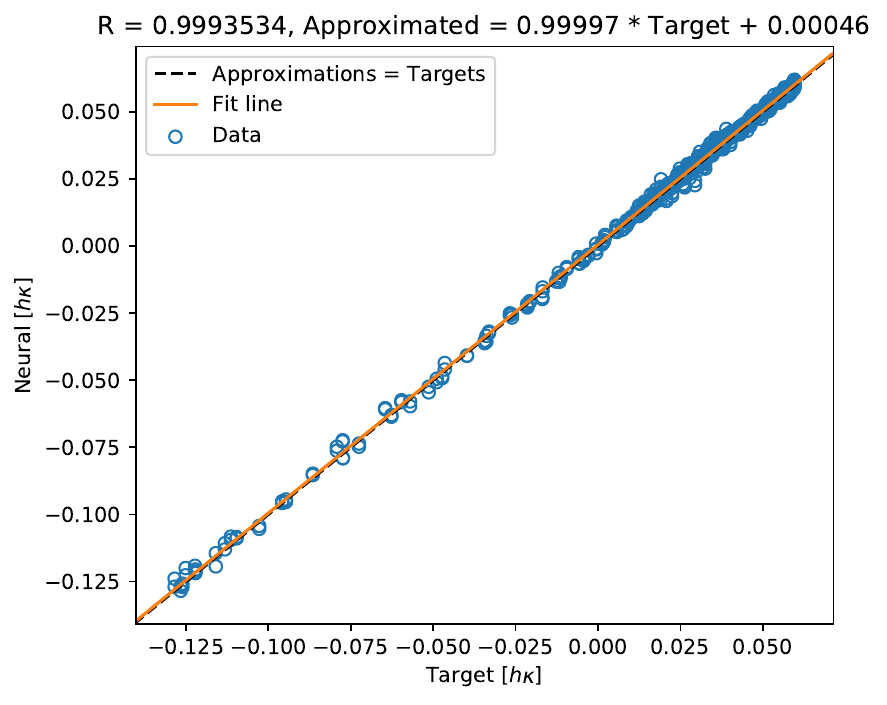}
        \caption{Neural, 5 iterations}
        \label{fig.flower.u.smooth.highRes.nnet.iter5}
    \end{subfigure}
	\begin{subfigure}[b]{0.3\textwidth}
		\includegraphics[width=\textwidth]{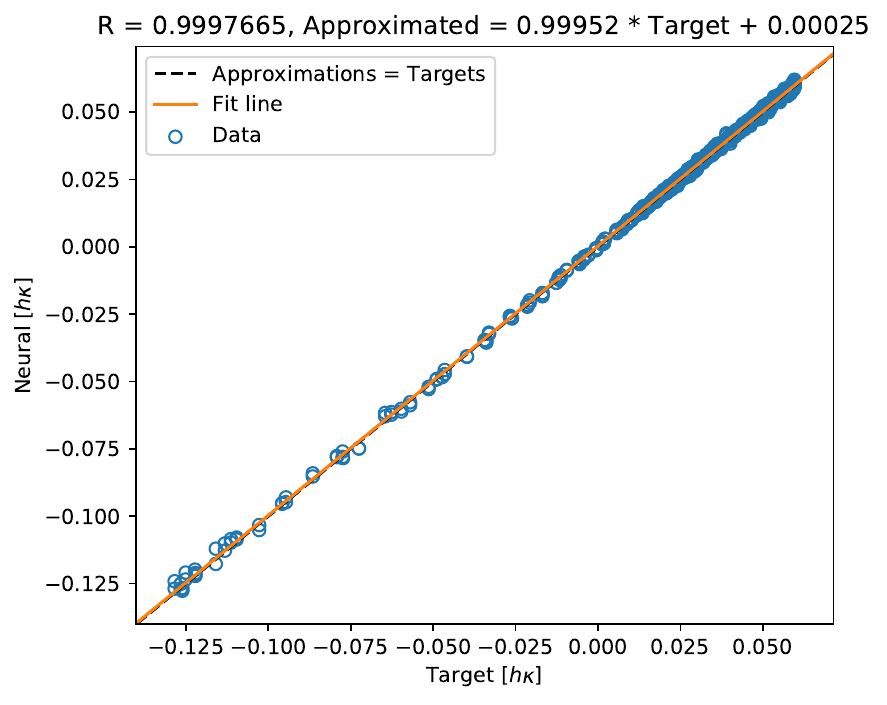}
		\caption{Neural, 10 iterations}
		\label{fig.flower.u.smooth.highRes.nnet.iter10}
	\end{subfigure}
	\begin{subfigure}[b]{0.3\textwidth}
		\includegraphics[width=\textwidth]{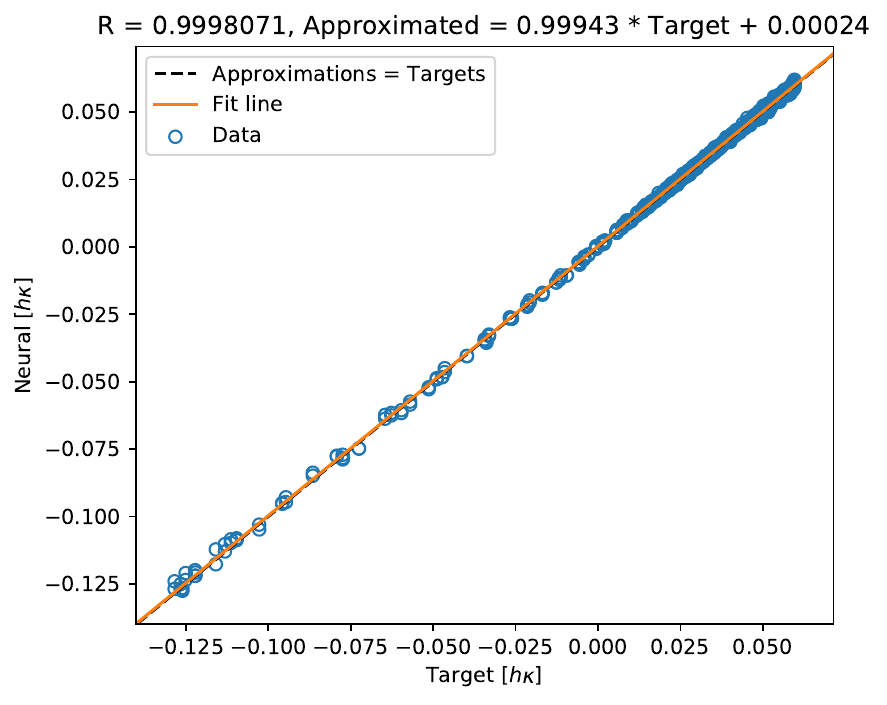}
		\caption{Neural, 20 iterations}
		\label{fig.flower.u.smooth.highRes.nnet.iter20}
	\end{subfigure}
    
	\begin{subfigure}[b]{0.3\textwidth}
		\includegraphics[width=\textwidth]{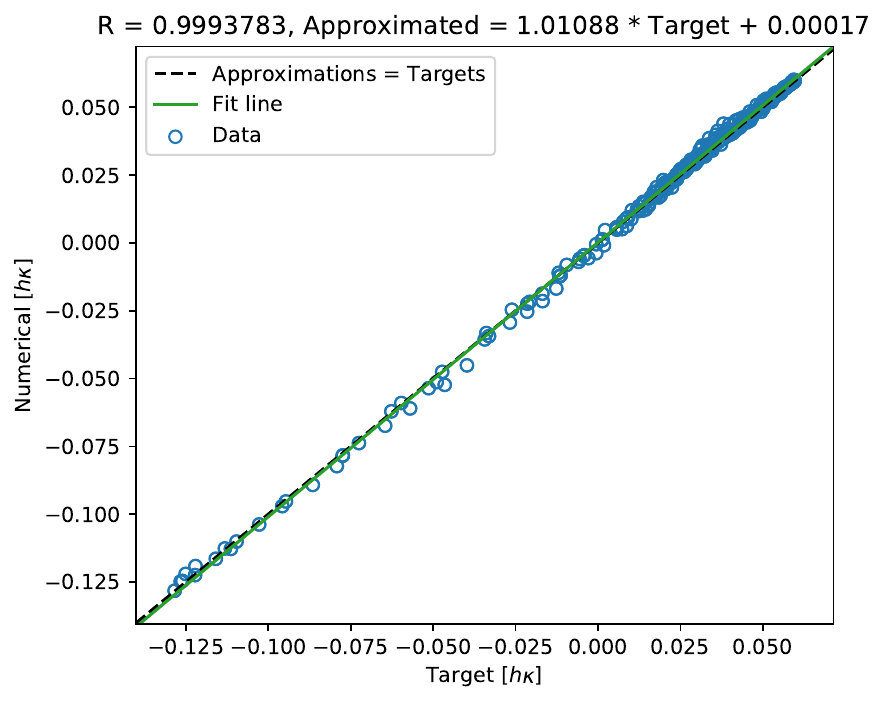}
		\caption{Numerical, 5 iterations}
		\label{fig.flower.u.smooth.highRes.numerics.iter5}
	\end{subfigure}
	\begin{subfigure}[b]{0.3\textwidth}
		\includegraphics[width=\textwidth]{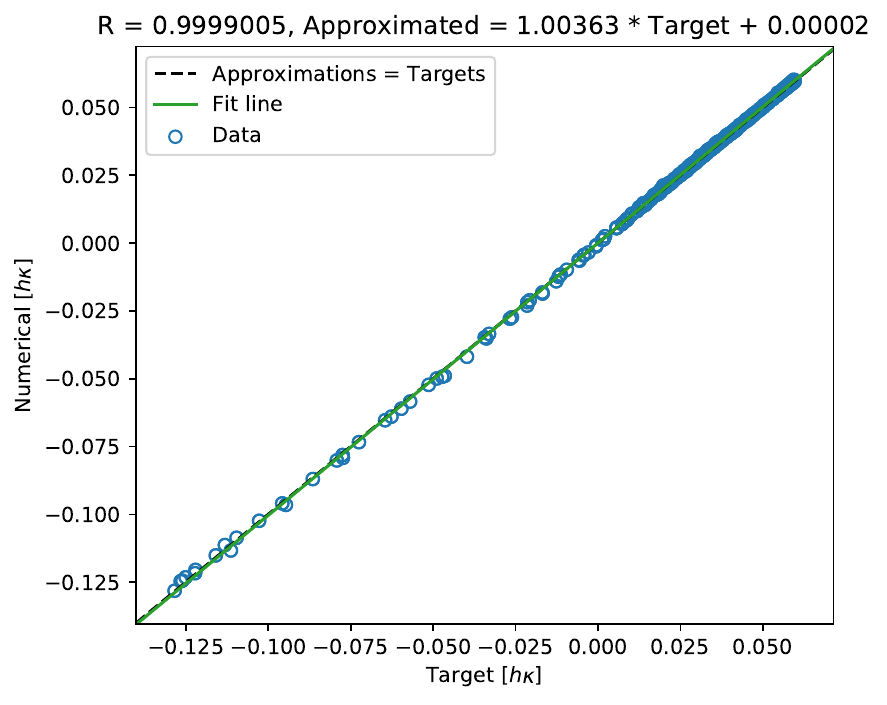}
		\caption{Numerical, 10 iterations}
		\label{fig.flower.u.smooth.highRes.numerics.iter10}
	\end{subfigure}
	\begin{subfigure}[b]{0.3\textwidth}
		\includegraphics[width=\textwidth]{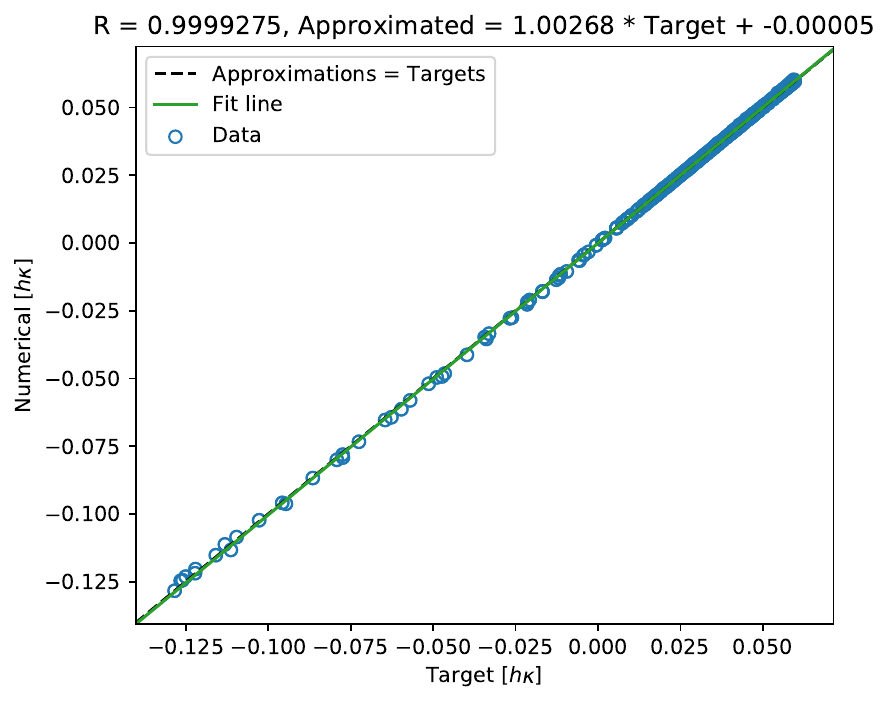}
		\caption{Numerical, 20 iterations}
		\label{fig.flower.u.smooth.highRes.numerics.iter20}
	\end{subfigure}
	\caption{\small Correlation between expected and approximated $h\kappa$ using the neural network and the numerical method for the smooth flower interface in a regular grid of $114 \times 114$ nodes.}
	\label{fig.flower.u.smooth.highRes.correlation}
\end{figure}

Figure \ref{fig.flower.u.smooth.highRes.correlation} contrasts the quality of the fit of our trained multilayer perceptron and how well the compound numerical method estimates curvature in this experiment.  Like in the previous cases, the neural accuracy is again comparable to the numerical precision.  One can prove this by noticing the high correlation factors across the plots in Figure \ref{fig.flower.u.smooth.highRes.correlation}.  However, the error summary provided in Table \ref{tbl.flower.u.smooth.highRes.errors} demonstrates that the accuracy gap between both methods stretches as one uses more iterations for level-set reinitialization.  Regarding the error $L^1$ and $L^2$ norms, for instance, when the number of iterations is small, both techniques estimate curvature at practically the same level.  On the other hand, when the number of iterations is large, the numerical framework approximates curvature with an MAE that is at most 44\% of the error incurred by the neural model.  A similar calculation shows that the error $L^\infty$ norm also improves up to 54\% in favor of the conventional approach.

\begin{table}[!b]
	\centering
	\footnotesize
	\bgroup
	\def\arraystretch{1.1}%
	\begin{tabular}{|l|l|r|r|r|}
		\hline
		Iterations & Method & MAE & Max AE & MSE \\
		\hline \hline
		\multirow{2}{*}{5} & Neural & $1.162374\times 10^{-3}$ & $6.806313\times 10^{-3}$ & $2.408337\times 10^{-6}$ \\
 		& Numerical & $1.117490\times 10^{-3}$ & $5.903280\times 10^{-3}$ & $2.507478\times 10^{-6}$ \\
		\hline \hline
		\multirow{2}{*}{10} & Neural & $6.950655\times 10^{-4}$ & $4.271494\times 10^{-3}$ & $8.508377\times 10^{-7}$ \\
 		& Numerical & $4.017341\times 10^{-4}$ & $2.461716\times 10^{-3}$ & $3.710802\times 10^{-7}$ \\
		\hline \hline
		\multirow{2}{*}{20} & Neural & $6.178167\times 10^{-4}$ & $4.365818\times 10^{-3}$ & $7.070641\times 10^{-7}$ \\
 		& Numerical & $2.748485\times 10^{-4}$ & $1.969042\times 10^{-3}$ & $2.593444\times 10^{-7}$ \\
		\hline
	\end{tabular}
	\egroup
	\caption{\small Error analysis for the smooth flower interface in a regular grid of $114 \times 114$ nodes.}
	\label{tbl.flower.u.smooth.highRes.errors}
\end{table}

\begin{figure}[t]
	\centering
	\begin{subfigure}[b]{0.3\textwidth}
		\includegraphics[width=\textwidth]{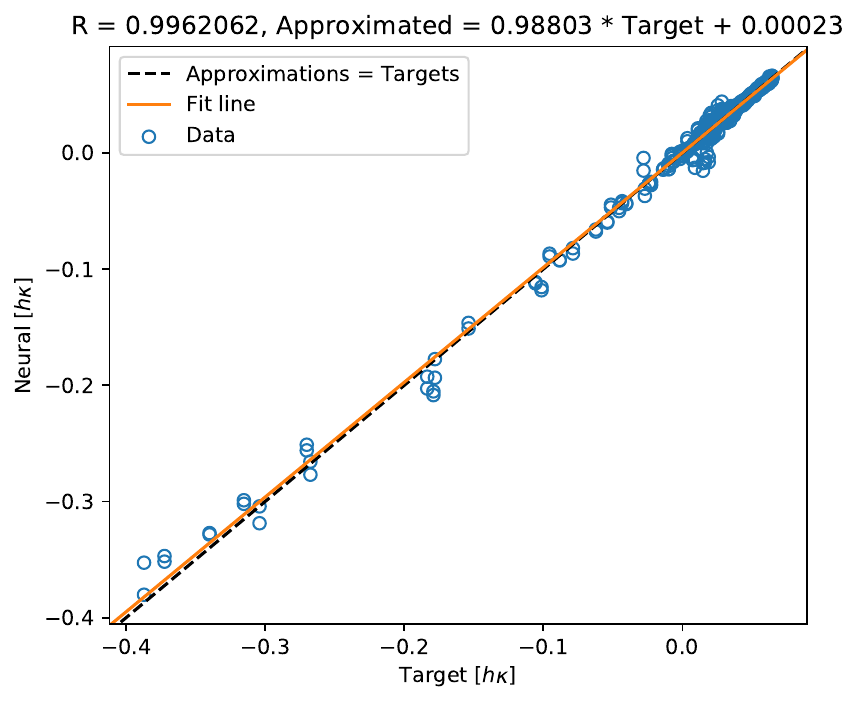}
        \caption{Neural, 5 iterations}
        \label{fig.flower.u.acute.highRes.nnet.iter5}
    \end{subfigure}
	\begin{subfigure}[b]{0.3\textwidth}
		\includegraphics[width=\textwidth]{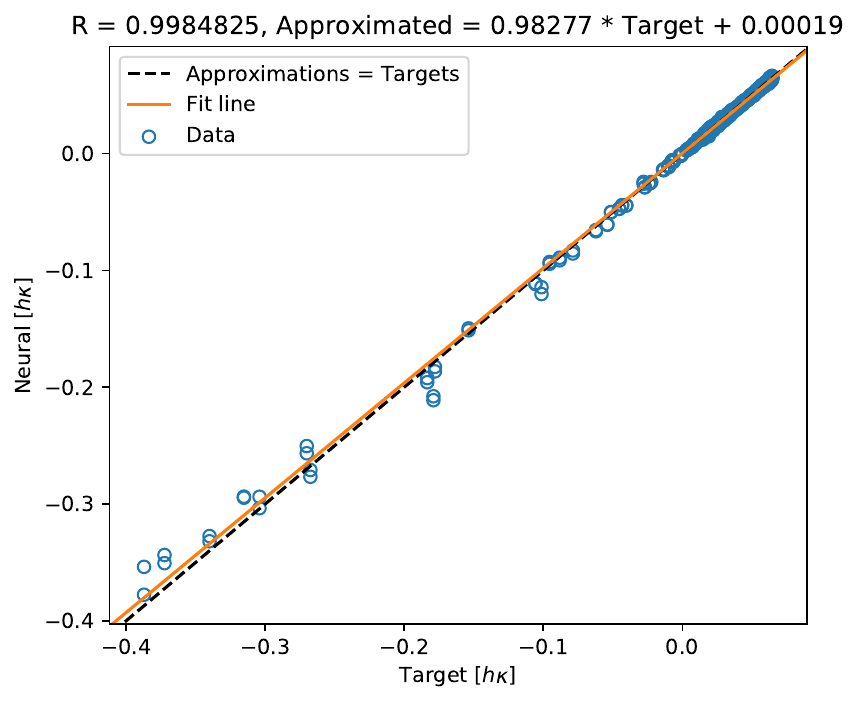}
		\caption{Neural, 10 iterations}
		\label{fig.flower.u.acute.highRes.nnet.iter10}
	\end{subfigure}
	\begin{subfigure}[b]{0.3\textwidth}
		\includegraphics[width=\textwidth]{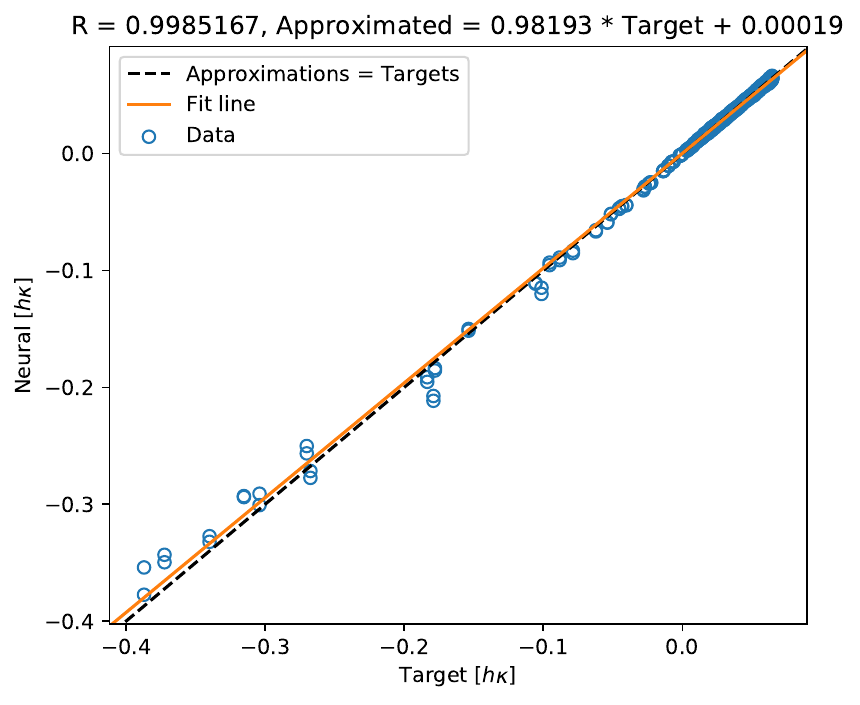}
		\caption{Neural, 20 iterations}
		\label{fig.flower.u.acute.highRes.nnet.iter20}
	\end{subfigure}
    
	\begin{subfigure}[b]{0.3\textwidth}
		\includegraphics[width=\textwidth]{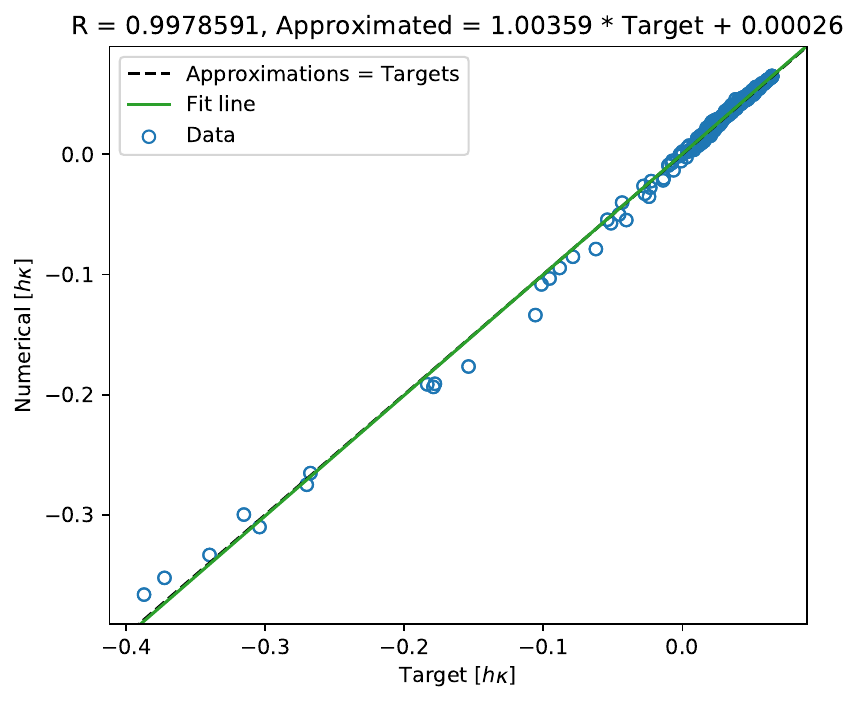}
		\caption{Numerical, 5 iterations}
		\label{fig.flower.u.acute.highRes.numerics.iter5}
	\end{subfigure}
	\begin{subfigure}[b]{0.3\textwidth}
		\includegraphics[width=\textwidth]{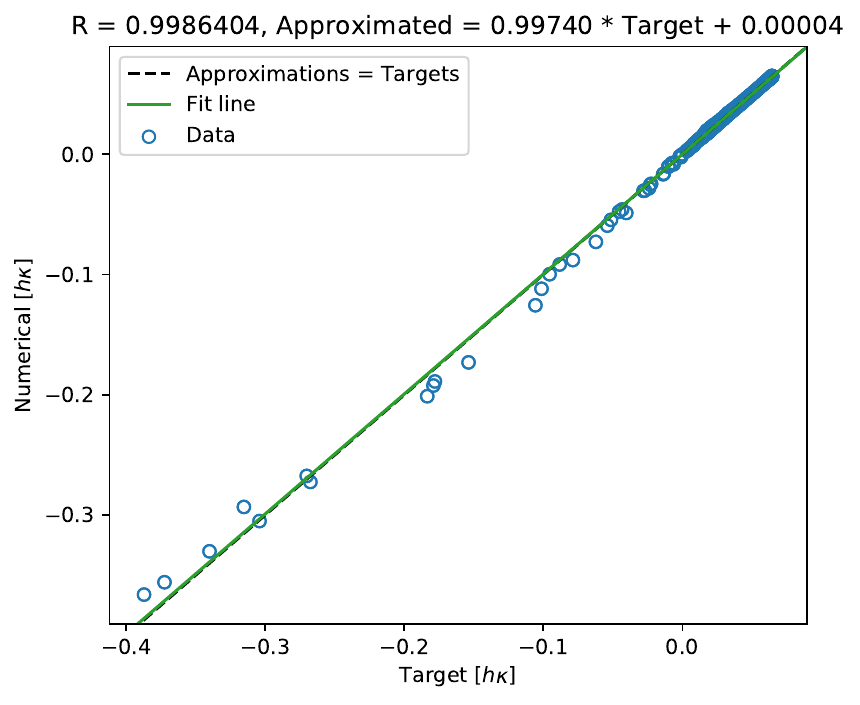}
		\caption{Numerical, 10 iterations}
		\label{fig.flower.u.acute.highRes.numerics.iter10}
	\end{subfigure}
	\begin{subfigure}[b]{0.3\textwidth}
		\includegraphics[width=\textwidth]{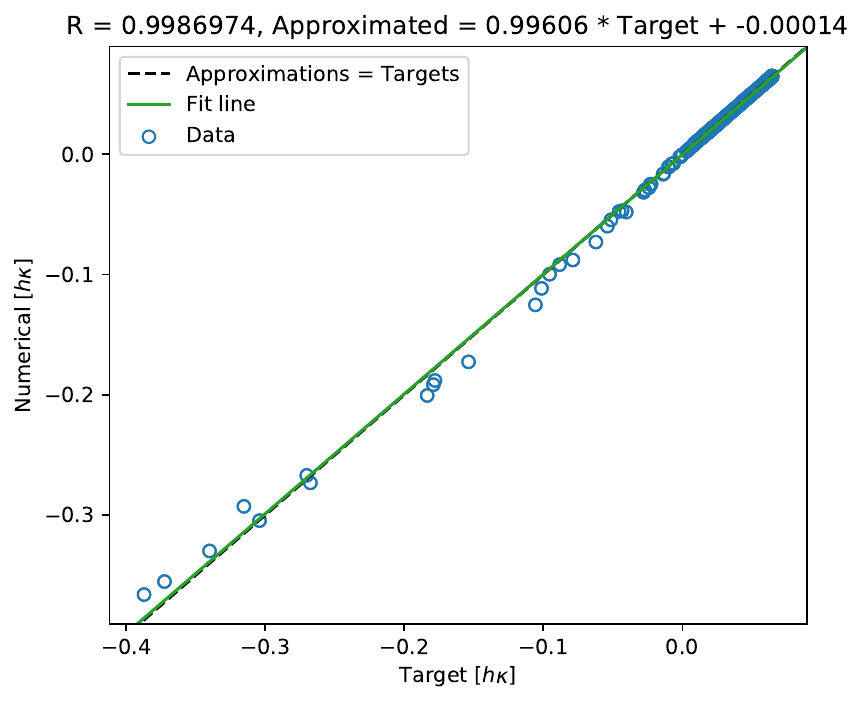}
		\caption{Numerical, 20 iterations}
		\label{fig.flower.u.acute.highRes.numerics.iter20}
	\end{subfigure}
	\caption{\small Correlation between expected and approximated $h\kappa$ using the neural network and the numerical method for the acute flower interface in a regular grid of $129 \times 129$ nodes.}
	\label{fig.flower.u.acute.highRes.correlation}
\end{figure}

The last curvature accuracy assessment for the flower interface in a high-resolution regular grid corresponds to $\Gamma_a$ (see Figure \ref{fig.flower.u.acute}).  As with the preceding analyses, we set the computational domain $\Omega \equiv [-0.232258, 0.232258]^2$ and discretize it into a $129 \times 129$ regular mesh.  This yields a cell's width of $h = 3.629032\times 10^{-3}$, which is equivalent to a unit square with 276.56 grid points per side.  This configuration allows us to apply the neural network trained for $\rho = 276$ and gather 672 samples along the interface.  Figure \ref{fig.flower.u.acute.highRes.correlation} shows the correlation plots between the target $h\kappa$ values and the neural and numerical approximations.

On average, our neural network generates curvature values with similar precision to the compound numerical approach.  We realize, however, that our neural approximations in Figure \ref{fig.flower.u.acute.highRes.correlation} are noisier than their numerical counterparts.  The issue is more obvious when one uses 5 iterations and as the target $h\kappa \rightarrow 0$.  In reality, this problem is a clear indicator that training a multilayer perceptron with only circular interfaces seems to not be enough to guarantee satisfactory results.  

We complement the previous findings with the statistics in Table \ref{tbl.flower.u.acute.highRes.errors}.  Despite the difficulties imposed by the steep curvatures in $\Gamma_a$, our deep learning strategy delivers $h\kappa$ estimations with a mean accuracy that does not differ significantly from the one in the numerical approach.  Regarding the error $L^\infty$ norm,  however, the compound numerical method is better than our framework by a proportion that is never larger than the improvement in $\Gamma_s$.  This can be seen if we compare the 34\% error reduction in this experiment, for example, with the 54\% encountered in Table \ref{tbl.flower.u.smooth.highRes.errors}.

\begin{table}[!b]
	\centering
	\footnotesize
	\bgroup
	\def\arraystretch{1.1}%
	\begin{tabular}{|l|l|r|r|r|}
		\hline
		Iterations & Method & MAE & Max AE & MSE \\
		\hline \hline
		\multirow{2}{*}{5} & Neural & $2.825795\times 10^{-3}$ & $3.444712\times 10^{-2}$ & $2.764972\times 10^{-5}$ \\
 		& Numerical & $2.281847\times 10^{-3}$ & $2.811109\times 10^{-2}$ & $1.590543\times 10^{-5}$ \\
		\hline \hline
		\multirow{2}{*}{10} & Neural & $1.213594\times 10^{-3}$ & $3.333340\times 10^{-2}$ & $1.179404\times 10^{-5}$ \\
 		& Numerical & $1.073979\times 10^{-3}$ & $2.198447\times 10^{-2}$ & $9.895744\times 10^{-6}$ \\
		\hline \hline
		\multirow{2}{*}{20} & Neural & $1.055267\times 10^{-3}$ & $3.302984\times 10^{-2}$ & $1.164381\times 10^{-5}$ \\
 		& Numerical & $8.570835\times 10^{-4}$ & $2.250886\times 10^{-2}$ & $9.531167\times 10^{-6}$ \\
		\hline
	\end{tabular}
	\egroup
	\caption{\small Error analysis for the acute flower interface in a regular grid of $129 \times 129$ nodes.}
	\label{tbl.flower.u.acute.highRes.errors}
\end{table}

\subsection{Embedding the irregular interface in an adaptive grid}
\label{subsec:adaptiveGrid}

\begin{figure}[t]
	\centering
	\includegraphics[width=0.6\textwidth]{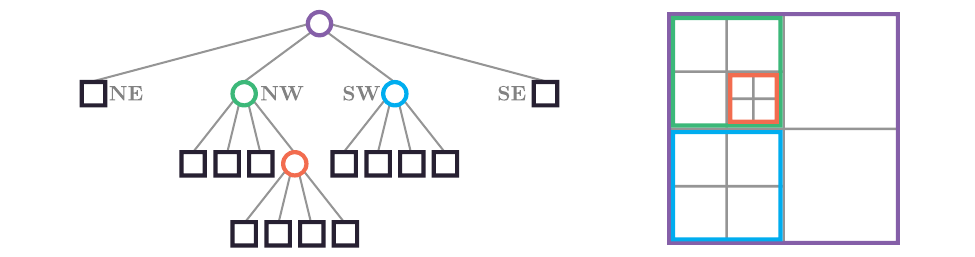}
	\caption{\small A quadtree and its cell subdivisions (adapted from \cite{BKOS00}).}
	\label{fig.quadtree}
\end{figure}

To close the results section, we assess the neural network $h\kappa$ approximations for $\Gamma_s$ and $\Gamma_a$ using an adaptive grid discretization.  A two-dimensional adaptive grid is a non-uniform quadtree mesh of vertices that covers the entire computational domain.  A quadtree is a rooted data structure where each tree cell, $\mathcal{C}$, either has 4 children (quadrants) or is a leaf \cite{BKOS00}.  The tree cells are organized into $L \ge 0$ levels, where the cell vertices (Cartesian nodes) store positional information and $\phi(x,y)$ values \cite{Strain1999}.  An example quadtree is shown in Figure \ref{fig.quadtree}. 

The process for building a quadtree consists of a recursive subdivision of cells that starts at the root ($l_0 \equiv \Omega$) and stops when the tree reaches a user-defined maximum level of refinement.  At every subdivision, we determine if a cell $\mathcal{C}$ needs to be split based on its distance to $\Gamma$ \cite{Strain1999}. This is expressed through the following criterion \cite{Min;Gibou:07:A-second-order-accur}:

\begin{equation}
\min_{v\in\mathcal{V}(\mathcal{C})}|\phi(v)| \leq \textrm{Lip}(\phi(v)) \times \textrm{diag}_\ell(\mathcal{C}),
\label{eq.quadtreeRefinementCriterion}
\end{equation}
where $\mathcal{V}$ is the set of $\mathcal{C}$'s vertices, $\textrm{diag}_\ell(\mathcal{C})$ defines the diagonal length of $\mathcal{C}$, and $\textrm{Lip}(\phi(v))$ is the Lipschitz constant of $\phi$ (here set to 1.2 for all cells).

\begin{figure}[t]
	\centering
	\begin{subfigure}[b]{0.4\textwidth}
		\includegraphics[width=\textwidth]{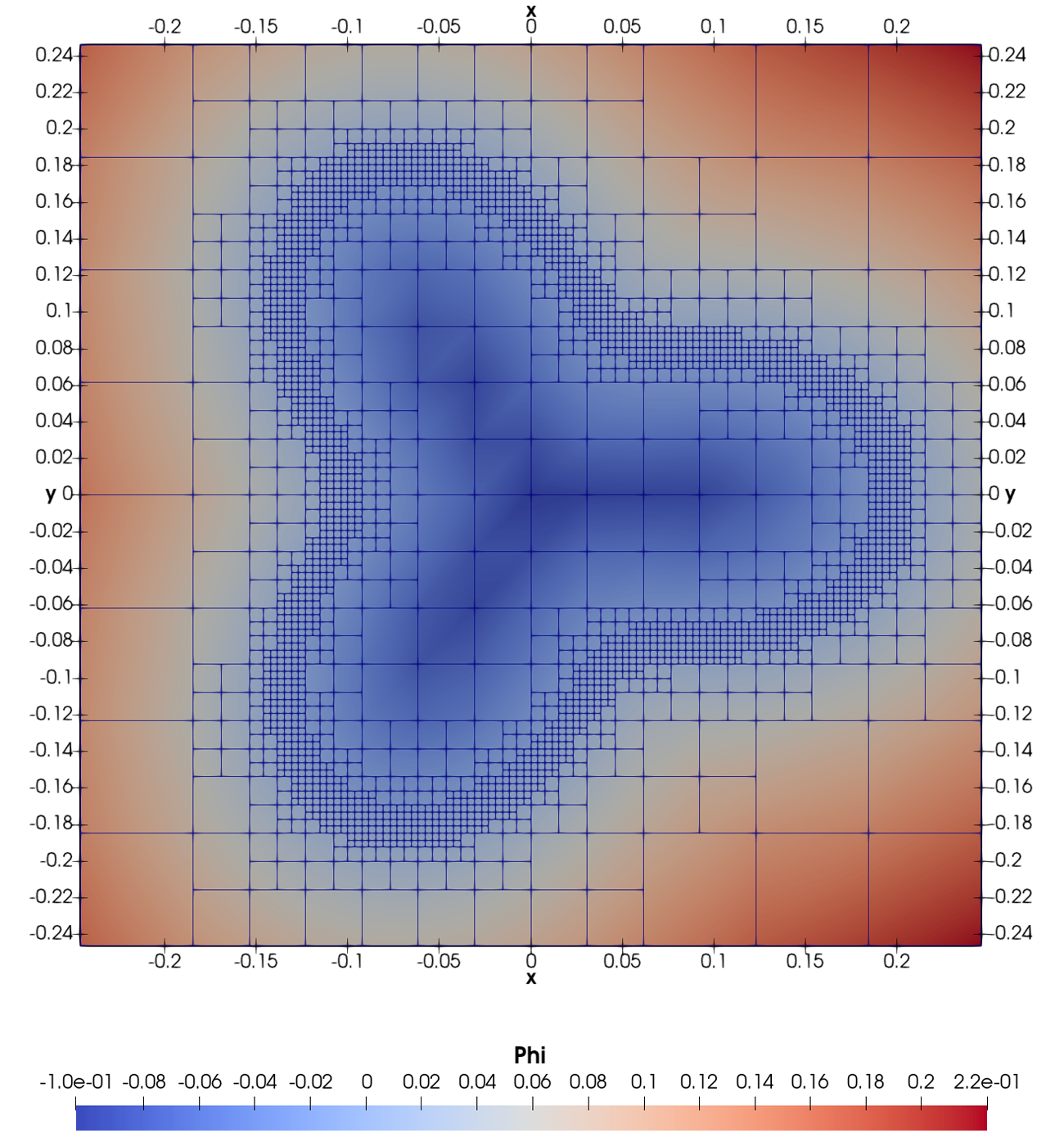}
        \caption{Quadtree discretization}
        \label{fig.flower.a.smooth.edges}
    \end{subfigure}\hfill%
	\begin{subfigure}[b]{0.4\textwidth}
		\includegraphics[width=\textwidth]{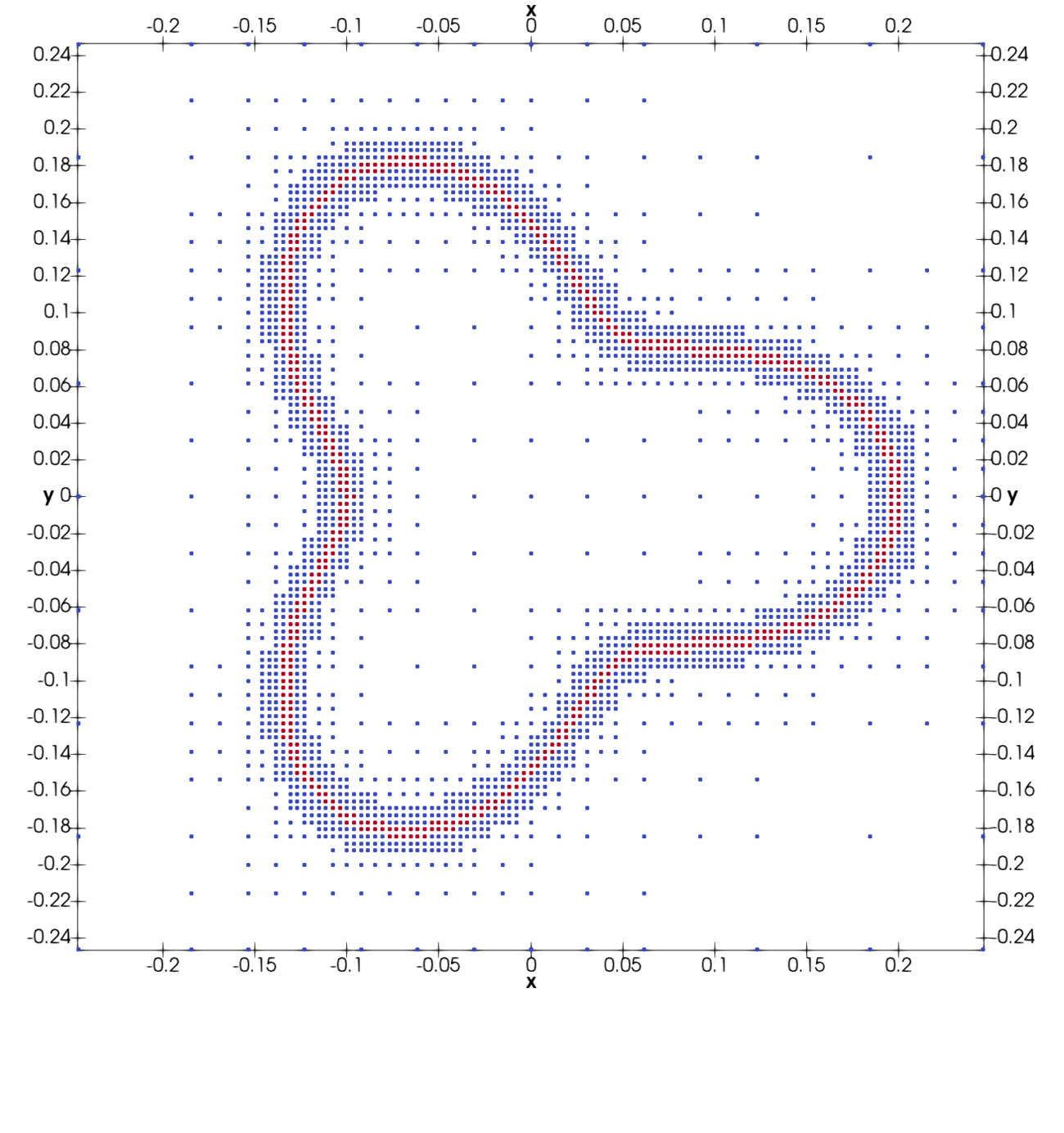}
		\caption{Sampled nodes along interface}
		\label{fig.flower.a.smooth.nodes}
	\end{subfigure}
	\caption{\small Smooth flower interface embedded in a non-uniform grid.  (For interpretation of the references to color in this figure, the reader is referred to the web version of this article.)}
	\label{fig.flower.a.smooth}
\end{figure}

First, we analyze the quality of the neural $h\kappa$ approximations for the smooth interface, $\Gamma_s$, in a nonuniform grid.  We define the two-dimensional computational domain $\Omega \equiv [-0.246154, 0.246154]^2$, discretize it using one quadtree with $l_{max} = 7$ (see Figure \ref{fig.flower.a.smooth.edges}), and get a minimum cell's width of $h = 3.846154\times 10^{-3}$.  This is equivalent to a local resolution of 261 uniform grid points per unit length and is compatible with the neural model trained for a $266 \times 266$ uniform mesh.  With these settings we collect 536 samples next to $\Gamma_s$, which we show in red in Figure \ref{fig.flower.a.smooth.nodes}.

\begin{figure}[t]
	\centering
	\begin{subfigure}[b]{0.3\textwidth}
		\includegraphics[width=\textwidth]{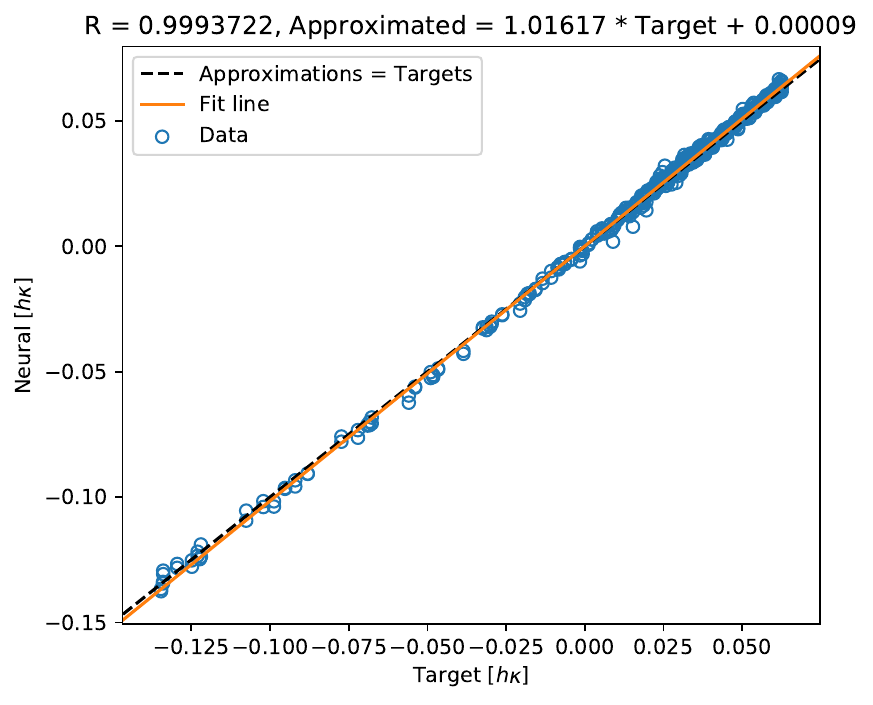}
        \caption{Neural, 5 iterations}
        \label{fig.flower.a.smooth.nnet.iter5}
    \end{subfigure}
	\begin{subfigure}[b]{0.3\textwidth}
		\includegraphics[width=\textwidth]{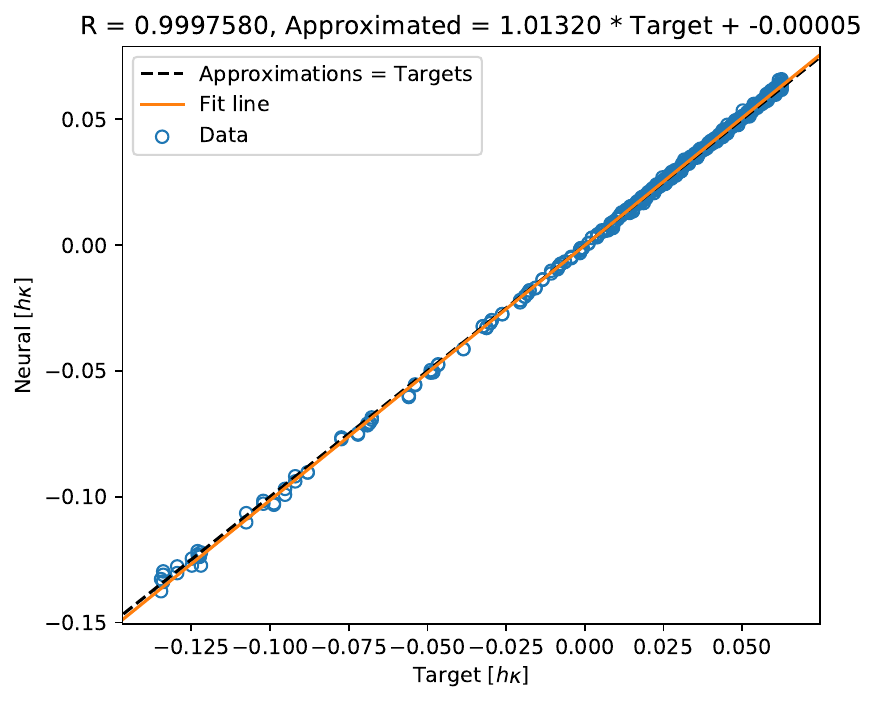}
		\caption{Neural, 10 iterations}
		\label{fig.flower.a.smooth.nnet.iter10}
	\end{subfigure}
	\begin{subfigure}[b]{0.3\textwidth}
		\includegraphics[width=\textwidth]{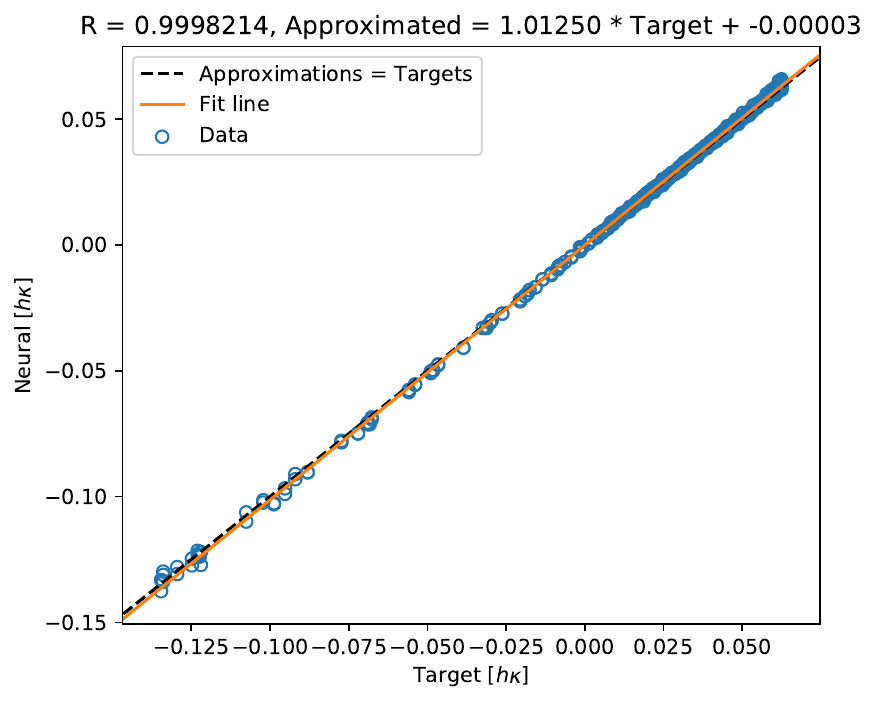}
		\caption{Neural, 20 iterations}
		\label{fig.flower.a.smooth.nnet.iter20}
	\end{subfigure}
    
	\begin{subfigure}[b]{0.3\textwidth}
		\includegraphics[width=\textwidth]{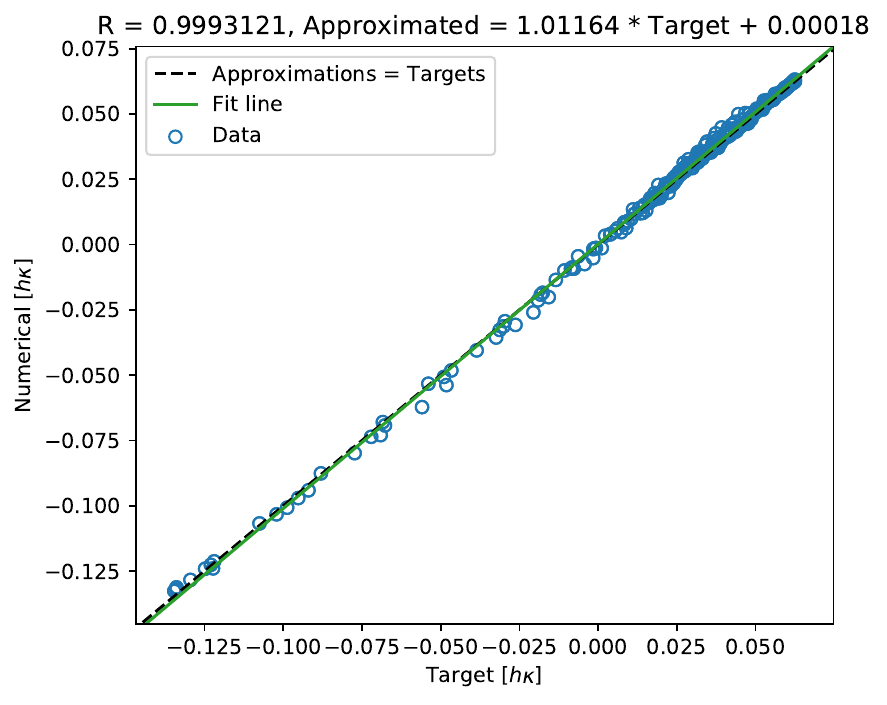}
		\caption{Numerical, 5 iterations}
		\label{fig.flower.a.smooth.numerics.iter5}
	\end{subfigure}
	\begin{subfigure}[b]{0.3\textwidth}
		\includegraphics[width=\textwidth]{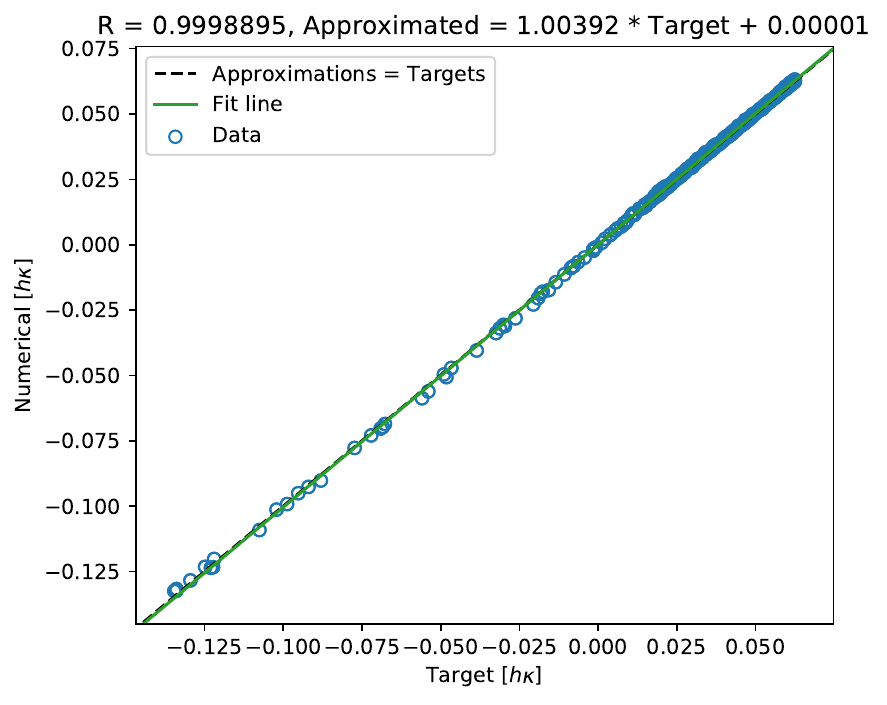}
		\caption{Numerical, 10 iterations}
		\label{fig.flower.a.smooth.numerics.iter10}
	\end{subfigure}
	\begin{subfigure}[b]{0.3\textwidth}
		\includegraphics[width=\textwidth]{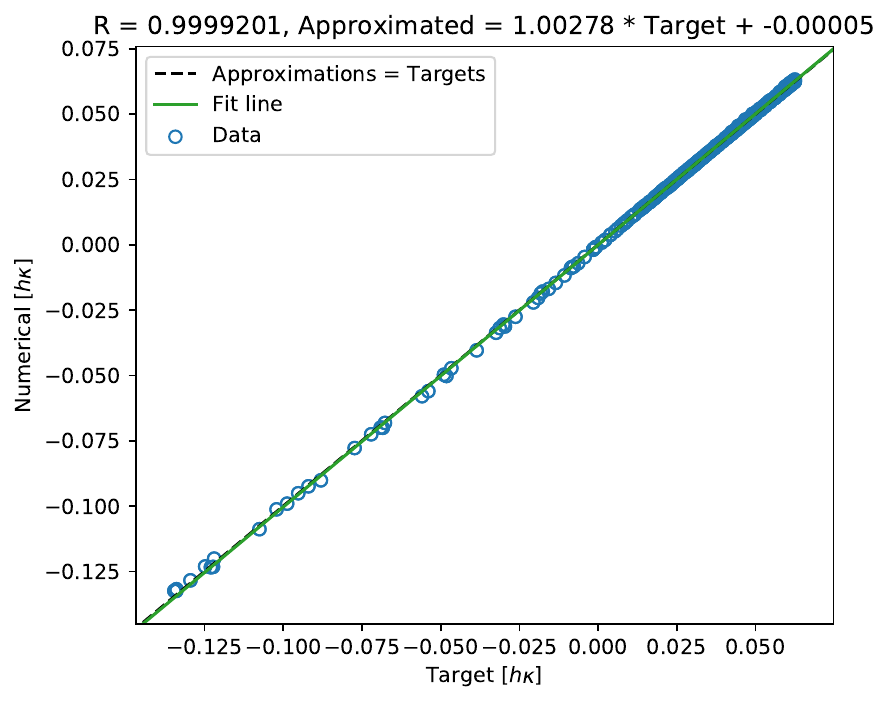}
		\caption{Numerical, 20 iterations}
		\label{fig.flower.a.smooth.numerics.iter20}
	\end{subfigure}
	\caption{\small Correlation between expected and approximated $h\kappa$ using the neural network and the numerical method for the smooth flower interface in an nonuniform grid.}
	\label{fig.flower.a.smooth.correlation}
\end{figure}

Figure \ref{fig.flower.a.smooth.correlation} illustrates the correlation between the expected $h\kappa$ values and their approximations when $\Gamma_s$ is embedded in a nonuniform grid.  In analogy to the regular grid test cases, the quality of the neural inferences is again comparable to the numerical estimations.  The error summary given in Table \ref{tbl.flower.a.smooth.errors} complements these visual results.  Once more, our multilayer perceptron can estimate mean curvature with practically the same average precision as the compound numerical method when one uses 5 iterations to reinitialize \eqref{eq.implicitFlower}.  As the number of iterations increases, however, the gap between the $L^1$ norms widens, and the corresponding error in the numerical method decreases by as much as 58\%.  Similar to the previous experiments, the finite-difference schemes followed by interpolation are also more robust in the error $L^\infty$ norm.  In this scenario, the numerical method exhibits an error reduction of as much as 47\% when compared with the neural maximum absolute error.

\begin{table}[!b]
	\centering
	\footnotesize
	\bgroup
	\def\arraystretch{1.1}%
	\begin{tabular}{|l|l|r|r|r|}
		\hline
		Iterations & Method & MAE & Max AE & MSE \\
		\hline \hline
		\multirow{2}{*}{5} & Neural & $1.325434\times 10^{-3}$ & $7.451553\times 10^{-3}$ & $2.979805\times 10^{-6}$ \\
 		& Numerical & $1.195590\times 10^{-3}$ & $6.273169\times 10^{-3}$ & $2.949928\times 10^{-6}$ \\
		\hline \hline
		\multirow{2}{*}{10} & Neural & $8.303893\times 10^{-4}$ & $5.367417\times 10^{-3}$ & $1.256426\times 10^{-6}$ \\
 		& Numerical & $4.289887\times 10^{-4}$ & $2.802827\times 10^{-3}$ & $4.350392\times 10^{-7}$ \\
		\hline \hline
		\multirow{2}{*}{20} & Neural & $7.147239\times 10^{-4}$ & $5.192075\times 10^{-3}$ & $9.882747\times 10^{-7}$ \\
 		& Numerical & $3.025235\times 10^{-4}$ & $2.215459\times 10^{-3}$ & $3.016474\times 10^{-7}$ \\
		\hline
	\end{tabular}
	\egroup
	\caption{\small Error analysis for the smooth flower interface in an adaptive grid.}
	\label{tbl.flower.a.smooth.errors}
\end{table}

\begin{figure}[t]
	\centering
	\begin{subfigure}[b]{0.4\textwidth}
		\includegraphics[width=\textwidth]{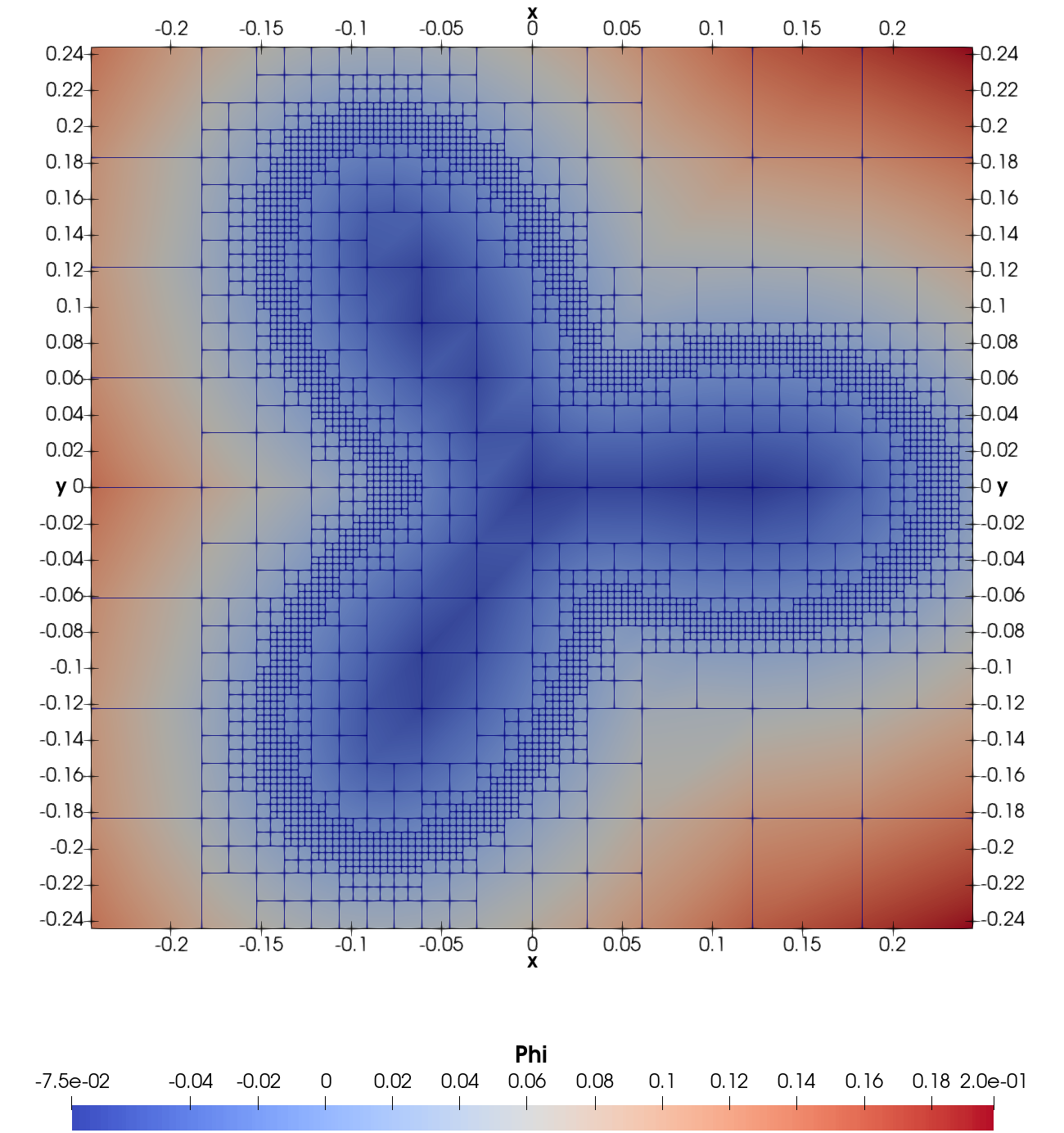}
        \caption{Quadtree discretization}
        \label{fig.flower.a.acute.edges}
    \end{subfigure}\hfill%
	\begin{subfigure}[b]{0.4\textwidth}
		\includegraphics[width=\textwidth]{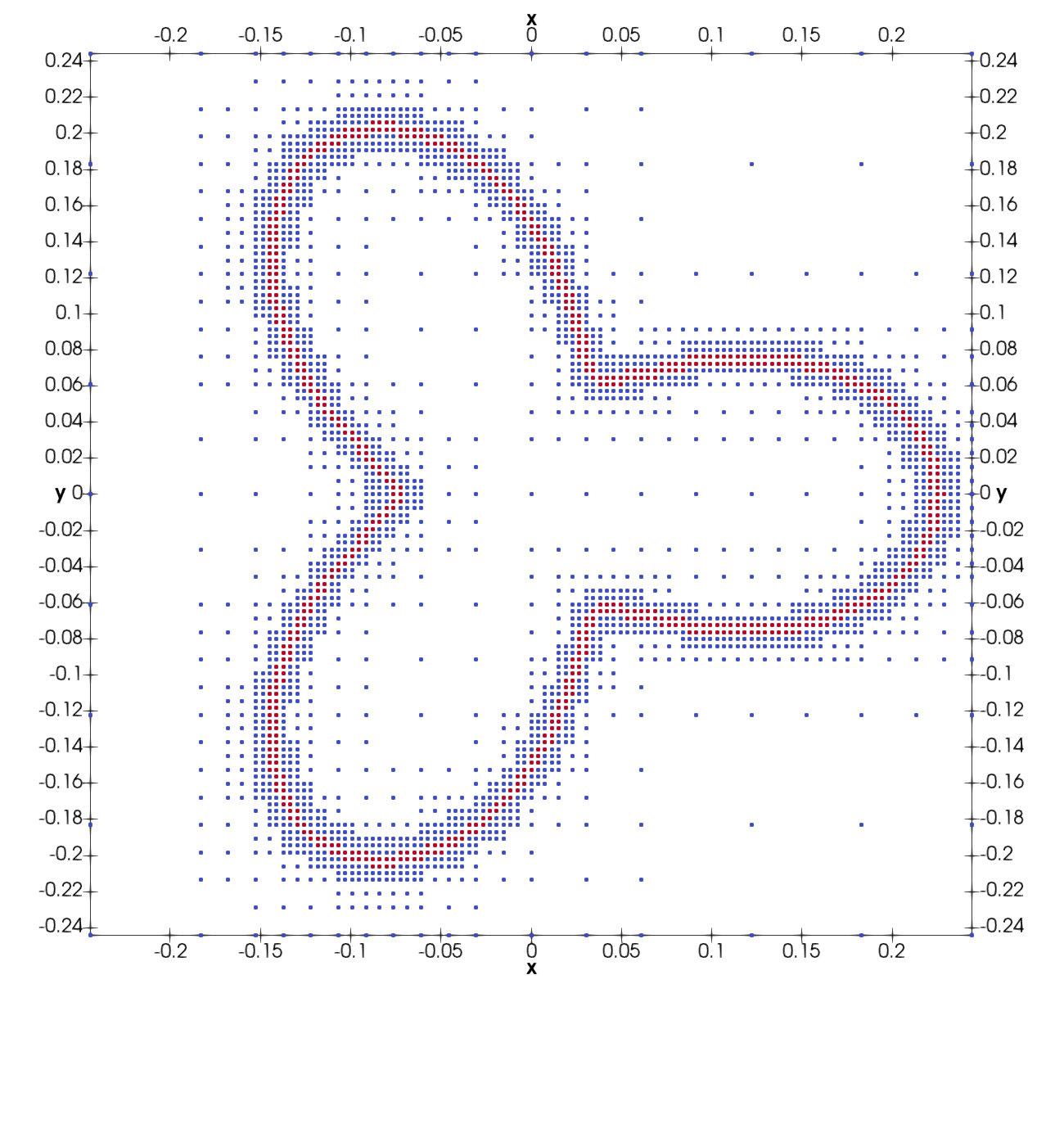}
		\caption{Sampled nodes along interface}
		\label{fig.flower.a.acute.nodes}
	\end{subfigure}
	\caption{\small Acute flower interface embedded in a nonuniform grid.  (For interpretation of the references to color in this figure, the reader is referred to the web version of this article.)}
	\label{fig.flower.a.acute}
\end{figure}

Last, we analyze the neural $h\kappa$ approximation for the acute interface when $\Gamma_a$ is embedded in a nonuniform grid.  We now set the computational domain to $\Omega \equiv [-0.244068, 0.244068]^2$ and discretize it using one quadtree with $l_{max} = 7$ (see Figure \ref{fig.flower.a.acute.edges}).  This yields a minimum cell's width of $h = 3.813559\times 10^{-3}$, which is equivalent to a unit square with 263.22 equally spaced nodes along each Cartesian direction.  We thus collect 644 samples and use the neural network trained for $\rho = 266$ to carry out the accuracy assessment.  Figure \ref{fig.flower.a.acute.nodes} shows the irregular interface, the sampled points, and the adaptive grid discretization.

\begin{figure}[t]
	\centering
	\begin{subfigure}[b]{0.3\textwidth}
		\includegraphics[width=\textwidth]{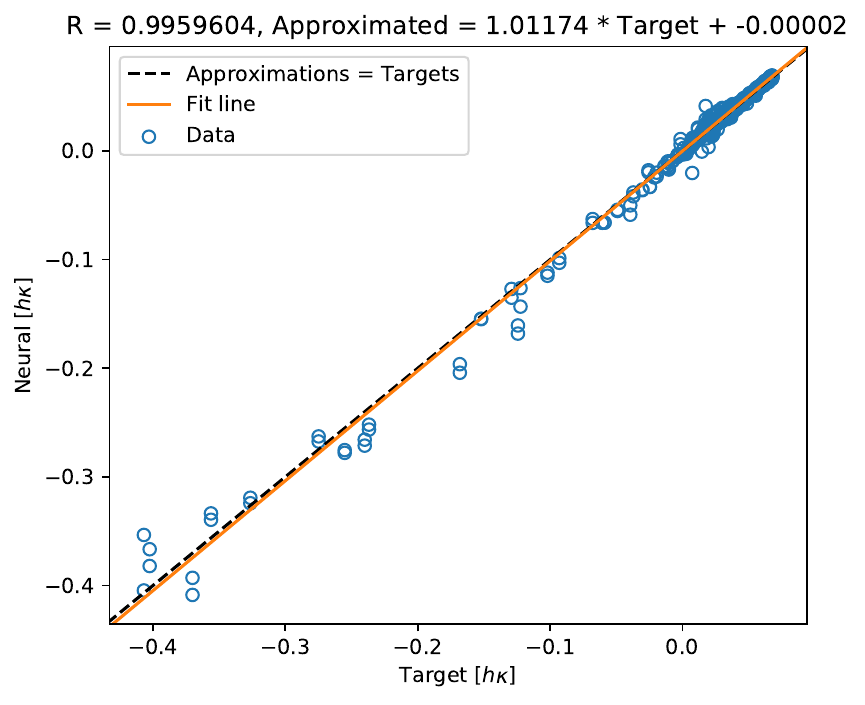}
        \caption{Neural, 5 iterations}
        \label{fig.flower.a.acute.nnet.iter5}
    \end{subfigure}
	\begin{subfigure}[b]{0.3\textwidth}
		\includegraphics[width=\textwidth]{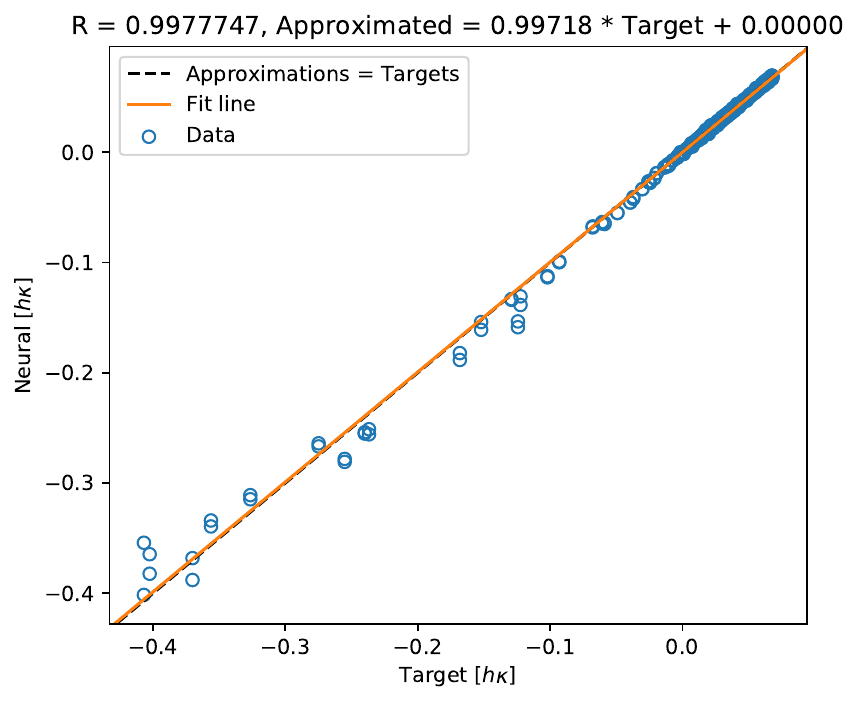}
		\caption{Neural, 10 iterations}
		\label{fig.flower.a.acute.nnet.iter10}
	\end{subfigure}
	\begin{subfigure}[b]{0.3\textwidth}
		\includegraphics[width=\textwidth]{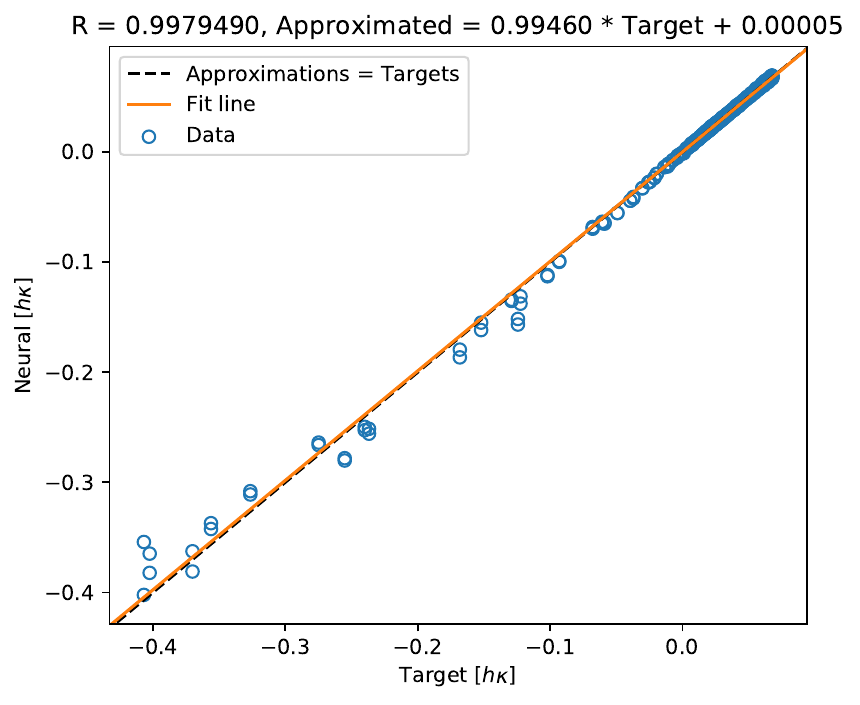}
		\caption{Neural, 20 iterations}
		\label{fig.flower.a.acute.nnet.iter20}
	\end{subfigure}
    
	\begin{subfigure}[b]{0.3\textwidth}
		\includegraphics[width=\textwidth]{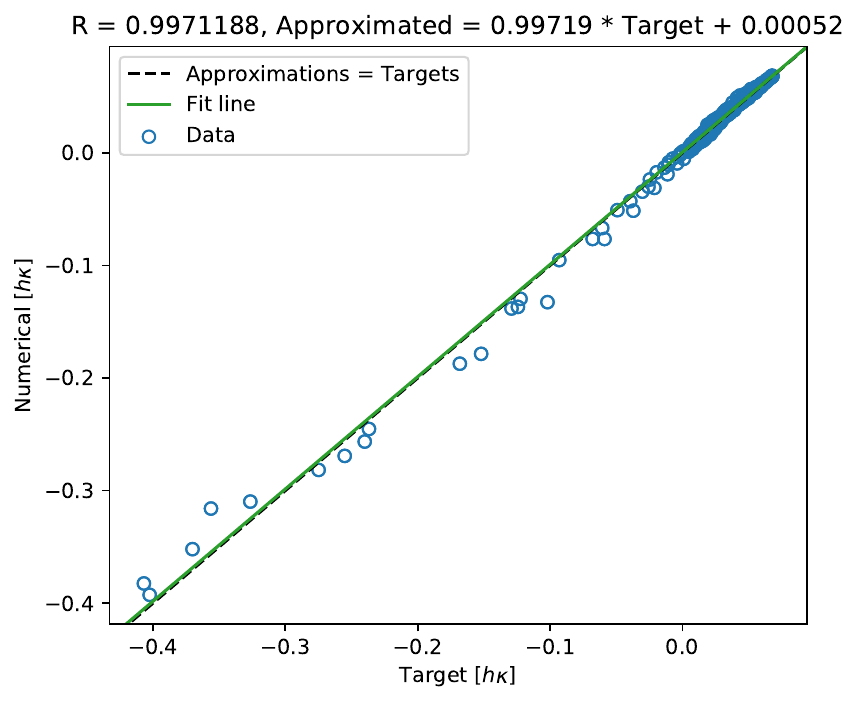}
		\caption{Numerical, 5 iterations}
		\label{fig.flower.a.acute.numerics.iter5}
	\end{subfigure}
	\begin{subfigure}[b]{0.3\textwidth}
		\includegraphics[width=\textwidth]{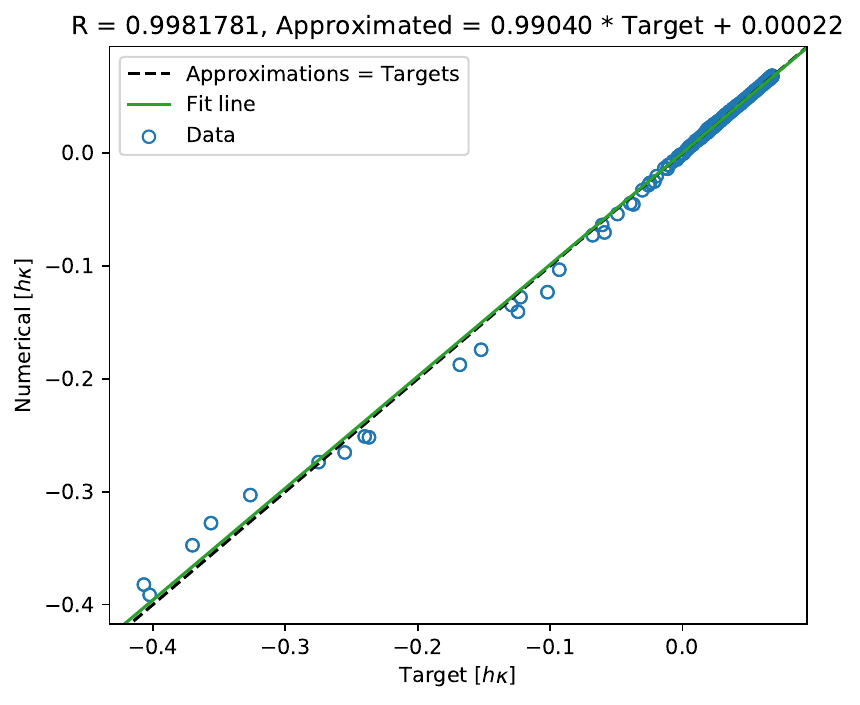}
		\caption{Numerical, 10 iterations}
		\label{fig.flower.a.acute.numerics.iter10}
	\end{subfigure}
	\begin{subfigure}[b]{0.3\textwidth}
		\includegraphics[width=\textwidth]{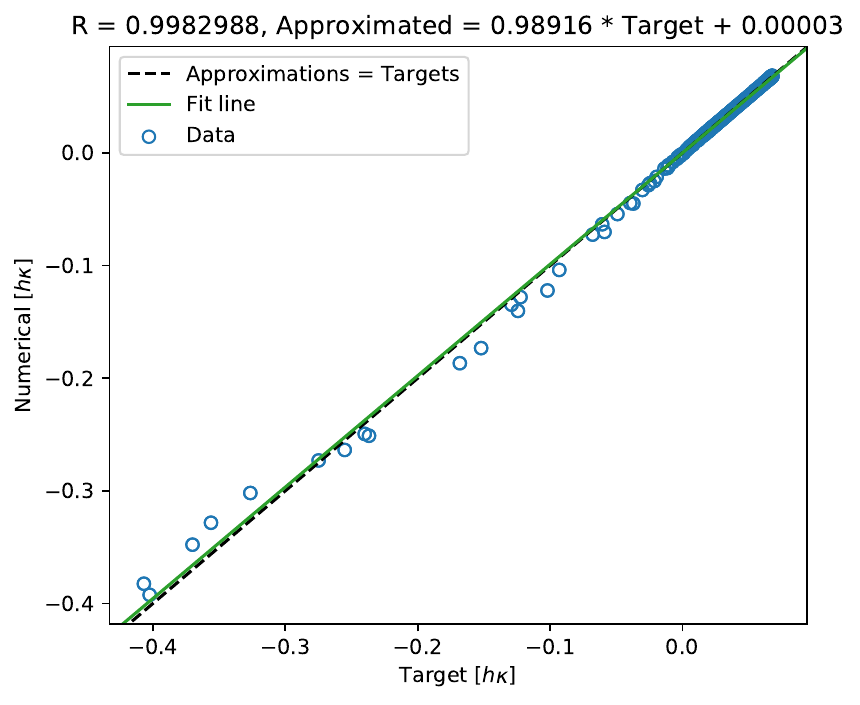}
		\caption{Numerical, 20 iterations}
		\label{fig.flower.a.acute.numerics.iter20}
	\end{subfigure}
	\caption{\small Correlation between expected and approximated $h\kappa$ using the neural network and the numerical method for the acute flower interface in a nonuniform grid.}
	\label{fig.flower.a.acute.correlation}
\end{figure}

Figure \ref{fig.flower.a.acute.correlation} compares the quality of the approximated $h\kappa$ values in the current study.  Despite the steep curvatures at the rapidly turning petal junctions, our deep learning strategy is again in good agreement with the results achieved by the compound numerical approach.  This can be seen in the high correlation factors across the plots in Figure \ref{fig.flower.a.acute.correlation}.   The statistical report in Table \ref{tbl.flower.a.acute.errors} supports these observations.  This summary shows that the neural network's $L^1$ and $L^2$ norms do not differ significantly from the corresponding errors in the numerical method.  As for the maximum absolute error, however, there is an important reduction of 47\% in favor of the conventional framework.  Interestingly, such a drop in the $L^\infty$ norm is atypical to the trend observed in the previous analyses of $\Gamma_a$.  This underlines the need for more fine-tuned training and for possibly considering a larger set of input patterns in $\mathcal{D}$.

\begin{table}[!b]
	\centering
	\footnotesize
	\bgroup
	\def\arraystretch{1.1}%
	\begin{tabular}{|l|l|r|r|r|}
		\hline
		Iterations & Method & MAE & Max AE & MSE \\
		\hline \hline
		\multirow{2}{*}{5} & Neural & $2.973764\times 10^{-3}$ & $5.340193\times 10^{-2}$ & $3.739998\times 10^{-5}$ \\
 		& Numerical & $2.701956\times 10^{-3}$ & $4.015226\times 10^{-2}$ & $2.567785\times 10^{-5}$ \\
		\hline \hline
		\multirow{2}{*}{10} & Neural & $1.481318\times 10^{-3}$ & $5.257659\times 10^{-2}$ & $1.965320\times 10^{-5}$ \\
 		& Numerical & $1.326347\times 10^{-3}$ & $2.827962\times 10^{-2}$ & $1.624310\times 10^{-5}$ \\
		\hline \hline
		\multirow{2}{*}{20} & Neural & $1.289340\times 10^{-3}$ & $5.250363\times 10^{-2}$ & $1.811191\times 10^{-5}$ \\
 		& Numerical & $1.082114\times 10^{-3}$ & $2.778140\times 10^{-2}$ & $1.528986\times 10^{-5}$ \\
		\hline
	\end{tabular}
	\egroup
	\caption{\small Error analysis for the acute flower interface in an adaptive grid.}
	\label{tbl.flower.a.acute.errors}
\end{table}


\section{Conclusions}
\label{sec:conclusions}

We have presented a novel deep learning approach\footnote{Our neural networks are available at \url{https://github.com/UCSB-CASL/LSCurvatureDL}.} for the estimation of mean curvature in the level-set method.  Our experiments on static, irregular interfaces, embedded in uniform and adaptive grids, show that our neural networks are competitive with conventional numerical methods at approximating $h\kappa$.

The case studies in section \ref{sec:results} indicate that, unlike the coupling between machine learning and VOF methods \cite{CurvatureML19, VOFCurvature3DML19}, level-set curvature neural networks trained only on circular interfaces do not always produce satisfactory results.  Our findings, however, show that deep learning has enormous potential for tackling the shortcomings in the level-set method.  As seen above, these drawbacks get particularly emphasized when dealing with coarse grids and steep curvatures.  Indeed, the accuracy analyses for the acute interface, $\Gamma_a$, corroborate that the neural and numerical $h\kappa$ estimations are practically equivalent in the $L^1$ and $L^2$ norms.  Furthermore, the very accurate numerical approximations in well-resolved regions suggest that a hybrid approach seems more reasonable for handling inferences in the whole curvature spectrum.

We also remark that our models cannot be transferred to domains where the grid resolutions differ significantly from the spatial discretization used in training \cite{BNK;2020;ML-for-FD}.  However, given some a priori knowledge about the application's computational domain, we claim that a dictionary of multilayer perceptrons can be constructed for various mesh sizes.  Then, one can use this map to approximate the dimensionless curvature of any interface.  A dictionary of neural networks presents two advantages over a universal neural model: (1) it is modular and faster to train, and (2) it is more efficient in both time and space.  As part of our future work, we plan to carry out a rigorous study to establish the resolution span over which one can apply a neural network with high confidence.  Similarly, given the statistics in Table \ref{tbl.nnetStats}, a meticulous convergence analysis is necessary to determine a sound relationship between mesh size and neural network precision and complexity.  In addition, we leave to future endeavors the integration of our deep learning strategy to a real simulation.  The latter will allow us to better assess our results not only with static but also with moving boundaries and level-set functions subject to advection and redistancing.

Another area of opportunity that deserves further exploration is the study of fast redistancing algorithms and machine learning for mean curvature estimation.  As we pointed out in section \ref{sec.buildingNNets}, level-set reinitialization has the benefit of removing unstructured noise that can otherwise undermine the performance of our neural networks.  However, solving the pseudotransient equation \eqref{eq.reinitialization} to steady state is expensive, so other reinitialization schemes, such as the fast marching method \cite{Sethian:96:A-Fast-Marching-Leve, Chacon;Vladimirsky:15:A-Parallel-Two-Scale, Chacon;Vladimirsky:12:Fast-two-scale-metho} and the fast sweeping method \cite{Zhao:04:A-Fast-Sweeping-Meth, Zhao:07:Parallel-Implementat, Detrixhe;Gibou;Min:13:A-parallel-fast-swee, Detrixhe;Gibou:14:Hybrid-Massively-Par}, are worth considering alongside machine learning in the pursuit of an efficient framework.

To conclude, we briefly mention our findings concerning a simple execution-time analysis of our neural models' performance.  Given a batch of 14,400 samples and a one-processor system, the compound numerical method requires about 83ms in C++ to calculate $h\kappa$ after using 20 iterations for reinitialization.  As for the multilayer perceptrons, these take an average of 351ms to output the corresponding estimations in Python.  Consequently, there is still much to be done to increase the neural networks' evaluation efficiency.  We can, for instance, migrate the neural networks to C++ and PETSc \cite{Balay;Brown;Buschelman;etal:12:PETSc-Web-page} or resort to model pruning \cite{ZG18} to favor sparsity and reduce the number of operations.



\section*{Acknowledgments}
Use was made of computational facilities purchased with funds from the National Science Foundation (CNS-1725797) and administered by the Center for Scientific Computing (CSC).  The CSC is supported by the California NanoSystems Institute and the Materials Research Science and Engineering Center (MRSEC; NSF DMR 1720256) at UC Santa Barbara.


\bibliographystyle{siamplain}
\bibliography{references}

\end{document}